%% file: thesis.tex
\documentclass[twoside]{ruthesis}
\usepackage{amssymb, times, latexsym}
\newcommand{\QED}{\hfill $\square$}
\newcommand{\PQED}{\hfill $\triangle$}
\def\proof{\noindent{\bf Proof. }}
\def\claim{\noindent{\bf Claim. }}
\def\poc{\noindent{\bf Proof of Claim. }}
\def\hint{\noindent{\bf Hint. }}

\newtheorem{Theorem}{Theorem}[section]
\newtheorem{Proposition}[Theorem]{Proposition}
\newtheorem{Lemma}[Theorem]{Lemma}
\newtheorem{Corollary}[Theorem]{Corollary}
\newtheorem{Definition}[Theorem]{Definition}
\newtheorem{Example}[Theorem]{Example}
\newtheorem{Problem}[Theorem]{Problem}
\newtheorem{Conjecture}[Theorem]{Conjecture}
\newtheorem{Exercise}[Theorem]{Exercise}

\newtheorem{Th}{Theorem}[chapter]
\newtheorem{Pro}[Th]{Proposition}
\newtheorem{Le}[Th]{Lemma}

\newtheorem{Ex}[Th]{Exercise}

\newtheorem{The}{Theorem}
\newtheorem{Def}[The]{Definition}

\newcommand{\lar}{\longrightarrow}
\newcommand{\depth}{\mbox{\rm depth}}
\newcommand{\length}{\ell}
\newcommand{\coker}{\mbox{\rm coker }}
\newcommand{\image}{\mbox{\rm image }}
\newcommand{\m}{\mathfrak{m}}

\newcommand{\height}{\mbox{\rm ht }}
\newcommand{\reg}{\mbox{\rm reg}}
\newcommand{\Ext}{\mbox{\rm Ext}}

\newcommand{\Hom}{\mbox{\rm Hom}}

\newcommand{\hdeg}{\mbox{\rm hdeg}}
\newcommand{\bdeg}{\mbox{\rm bdeg}}
\newcommand{\hdil}{\mbox{\rm hdil}}
\newcommand{\dil}{\mbox{\rm dil}}
\newcommand{\Deg}{\mbox{\rm Deg}}
\newcommand{\type}{\mbox{\rm type}}
\newcommand{\socle}{\mbox{\rm soc}}
\newcommand{\onto}{\to\!\!\!\!\!\to}
\newcommand{\grm}{\mbox{\rm gr}_{\m}}
\newcommand{\Ass}{\mbox{\rm Ass}}

\newcommand{\red}{\mbox{\rm r}}
\newcommand{\gm}{\Gamma_{\m}}
\newcommand{\gin}{\mbox{\rm gin}}
\newcommand{\init}{\mbox{\rm in}}
\newcommand{\Gl}{\mbox{\rm GL}}

\newcommand{\adeg}{\mbox{\rm adeg}}
\newcommand{\gdeg}{\mbox{\rm gdeg}}
\newcommand{\range}{\mbox{\rm range}}
\newcommand{\Mon}{\mbox{\it Mon}}
\newcommand{\embcod}{\mbox{\rm embcod}}
\newcommand{\embdim}{\mbox{\rm embdim}}
\newcommand{\Spec}{\mbox{\rm Spec}}
\newcommand{\rees}{\mbox{\it rees}}

\renewcommand{\deg}{\mbox{\rm deg}}
\renewcommand{\dim}{\mbox{\rm dim\,}}
\renewcommand{\max}{\mbox{\rm max}}
\renewcommand{\min}{\mbox{\rm min}}
\renewcommand{\sup}{\mbox{\rm sup}}
\renewcommand{\ker}{\mbox{\rm ker\,}}

\begin{document}

	\include{thesis-note}
    \include{thesis-title}

    \afterpreface 
    \include{thesis-intro}
    \include{thesis-chap1}

    \include{thesis-chap2}

    \include{thesis-chap3}

    \include{thesis-chap4}

    \appendix
    \include{thesis-appendix}

    \bibliography{thesis}
    \bibliographystyle{plain}
    \include{thesis-vita} 
\end{document}

%% file: thesis-note.tex
\chapter*{Note}
\pagestyle{empty}
This version differs slightly from the original.  Typographical errors have been corrected
without comment.  More substantial changes are indicated by footnotes.  The numbering of 
definitions, propositions and references is unchanged.

%% file: thesis-title.tex

    \phd
\title{\textbf{Cohomological Degrees, Dilworth Numbers and Linear Resolution}}
    \author{ Tor Kenneth Gunston }
    \campus{New Brunswick}
    \program{Mathematics}
    \director{Wolmer V. Vasconcelos}
    \approvals{4}
    \copyrightpage 
    \submissionmonth{October}
    \submissionyear{1998}
    \abstract{\input{thesis-abstract}}  
\beforepreface 

    \acknowledgements{
I am deeply indebted to my thesis advisor, Wolmer Vasconcelos, a warm and decent person, generous with his time and attention, tolerant of my faults, a Fountain of Algebra, and a wonderful {\em teacher} in the broadest -- and noblest -- sense of the word.

I am also grateful to the balance of my patient thesis committee, Joseph Brennan, Friedrich Knop and Charles Weibel.  Friedrich Knop and Charles Weibel offered excellent instruction in the courses basic to this dissertation: Commutative Algebra, Algebraic Geometry and Homological Algebra.  Charles Weibel devoted considerable energy to reading this thesis and suggesting numerous improvements.  Joseph Brennan was like a second advisor to me during his visit in the Fall of 1996 and the Spring of 1997.  It was Joe who sparked my interest in Dilworth numbers.  The ideas in Chapter 4 of this thesis were worked out with his considerable help.

My friend Luisa Doering, with whom I shared an office and an advisor, and with whom I learned Commutative Algebra, served as audience, critic and collaborator for many of the ideas in this dissertation.

This thesis would not have been completed without the support of Carol Hamer.

My Uncle Wayne and Aunt Barbara Perkins and my late grandparents, Doris and Howard Gunston, have given me generous financial support throughout my education.  I can not thank them enough.

Finally, I thank my family for laying a firm foundation.

}
    \dedication{\vspace{7cm}\begin{center}This is dedicated to my family: Mom, Dad, Hope, Emily and Christian.\end{center} }

    \abbreviations{\input{thesis-abbreviations}}


%% file: thesis-abstract.tex
This thesis is a study of various ways of measuring the size and complexity of finitely generated $R$-modules, where $(R,\m)$ is a Noetherian local ring, by attaching a number $\delta(M)$ to each finitely generated $R$-module $M$.  The classical example is $\delta(M)= \deg(M)$, the {\em multiplicity} or {\em degree} of $M$.  Here we investigate several variants of the degree function: the {\em homological Dilworth number} $\hdil(-)$ and the family of {\em cohomological degrees}, such as the homological degree, $\hdeg(-)$, and the extremal cohomological degree, $\bdeg(-)$.

A cohomological degree is a function $\Deg(-)$ which assigns to every isomorphism class of finitely generated $R$-modules a real number $\Deg(M)$ subject to the following conditions: (1) If $M$ is Cohen-Macaulay, then $\Deg(M) = \deg(M)$; (2) If $\depth(M)>0$ and $x$ is a generic hyperplane on $M$, then $\Deg(M) > \Deg(M/xM)$; and (3) $\Deg(M) = \length(H_{\m}^0(M)) + \Deg(M/H_{\m}^0(M))$.

The Dilworth number, $\dil(M)$, of a module $M$ is given by $\dil(M) = \sup\{\,\nu(N)\;|\; N \subseteq M\,\}$, where $\nu(N)$ denotes the size of a minimal generating set for the submodule $N$.  The Dilworth number is finite only in dimension $0$ and $1$.  The homological Dilworth number extends the Dilworth number to higher dimensions.  Its definition is similar to the definition of the homological degree.

Sally, Valla and others have established bounds for the number of generators of ideals (or modules) in terms of multiplicities and other numerical data, usually under the assumption that the ideal is Cohen-Macaulay.  We use cohomological degrees and the homological Dilworth number in place of the classical degree to extend some of these results from the Cohen-Macaulay case to the non-Cohen-Macaulay case.

We give particular attention to the following result: If $M$ is Cohen-Macaulay, then $\nu(M) \leq \deg(M)$, and equality holds if and only if $\grm(M)$, viewed as a module over a polynomial ring, has a linear resolution.  We give extensions of this theorem to the non-Cohen-Macaulay case (Theorems \ref{Yosh}, \ref{linresthm}, \ref{main} and Corollary \ref{dil-dim2}).  Our results are strongest for the extremal cohomological degree, $\bdeg(-)$, and this provides an avenue for connecting cohomological degrees with Castelnuovo-Mumford regularity.

%% file: thesis-abbreviations.tex

\begin{descriptionlist}{box width} 
\item[$\bf N$] the natural numbers

\item[$\bf R$] the real numbers
\item[$\bf Z$] the integers
\item[${\cal M}(R)$] finitely generated $R$-modules, graded if $R$ is graded; Section 1.1
\item[$M_{\langle i \rangle}$] the $i^{th}$ graded component of $M$; Section 1.1
\item[$H(M;n)$] the Hilbert function of $M$; Section 1.1
\item[$\deg(M)$] the multiplicity or degree of $M$; Section 1.1
\item[$\length(M)$] the length of $M$
\item[$H_M(t)$] the Hilbert series of $M$; Section 1.1
\item[$\grm(M)$] the associated graded module of $M$; Section 1.1
\item[$\chi^I(M;n)$] the Hilbert-Samuel function of $M$ with respect to $I$; Section 1.1
\item[$e(I;M)$] the multiplicity of $M$ with respect to $I$; Section 1.1
\item[$e({\bf x};M)$] the multiplicity symbol; Definition 1.2.2
\item[$\red_I(J)$] the reduction number of $J$ with respect to $I$; Definition 1.3.3
\item[$\red(J)$] the reduction number of $J$; Definition 1.3.3
\item[$\depth(M)$] the depth of $M$
\item[$\type(M)$] the type of $M$; Definition 1.4.1
\item[$\socle(M)$] the socle of $M$; Definition 1.4.3
\item[$\hat{M}$] the $\m$-adic completion of $M$
\item[$E_R(M)$] the injective envelope of the $R$-module $M$; Definition 1.5.1
\item[$M^{\vee}$] the Matlis dual of $M$; Section 1.5
\item[${\cal M}_0(R)$] the category of $R$-modules of finite length, graded if $R$ is graded; Section 1.5
\item[$\nu(M)$] the size of a minimal generating set of $M$
\item[$\gm(M)$] the zeroth local cohomology module of $M$; Definition 1.5.3
\item[$H_{m}^i(M)$] the $i^{th}$ local cohomology module of $M$; Definition 1.5.5
\item[$M_i$] $\Ext_S^{n-i}(M,S)$, where $S$ is a Gorenstein ring of dimension $n$; Section 1.5
\item[$\Spec(R)$] the set of prime ideals in $R$
\item[$\Ass(M)$] the set of primes associated to $M$
\item[$I(M)$] the invariant of St\"{u}ckrad and Vogel; Definition 1.6.4
\item[$M_{\geq r}$] the truncation $\oplus_{i\geq r}M_{\langle i \rangle}$ of the graded module $M$; Definition 1.7.3
\item[$\reg(M)$] the Castelnuovo-Mumford regularity of $M$; Definition 1.7.5
\item[$\succ_{rlex}$] the reverse lexicographic order; Section 1.9
\item[$\deg_{x_i}(m)$] the degree of the monomial $m$ in the variable $x_i$
\item[$\init(f)$] the initial term of $f$ with respect to the reverse lexicographic order; Section 1.9
\item[$\init(I)$] the initial ideal of $I$ with respect to the reverse lexicographic order; Section 1.9
\item[$\Gl(V)$] the group of inveritible linear transformations of the vector space $V$
\item[$\gin(I)$] the generic initial ideal of $I$ with respect to the reverse lexicographic order; Definition 1.9.3
\item[$I^{sat}$] the saturation of the ideal $I$; Section 1.9
\item[$S_h, I_h$] the restriction of the ring $S$ and ideal $I$ to the hyperplane $h$; Section 1.9
\item[$\hdeg$] the homological degree; Example 2.1.3
\item[$\bdeg_U$] the extremal cohomological degree with respect to $U$; Example 2.1.6
\item[$\embdim(R)$] the embedding dimension of $R$; Section 2.2
\item[$\embcod(R)$] the embedding codimension of $R$; Section 2.2
\item[$\bdeg$] the extremal cohomological degree; Definition 3.1.3
\item[${\cal E}(M)$] Definition 3.1.10
\item[$\rees(M)$] the Rees number of $M$; Definition 3.1.11
\item[$e^{(r,d)}$] see Section 3.4
\item[$\dil(M)$] the Dilworth number of $M$; Definition 4.2.1
\item[$D(M)$] the unique maximal submodule $N$ of $M$ satisfying $\nu(N) = \dil(M)$; Definition 4.2.8
\item[$\hdil$] the homological Dilworth number; Definition 4.4.1
\item[$\triangle$] the end of a portion of a proof
\item[$\square$] the end of a proof
\end{descriptionlist}

%% file: thesis-intro.tex
\chapter*{Introduction}
\addcontentsline{toc}{chapter}{Introduction}

This dissertation is about measuring modules.

\bigskip

Let $(R,\m)$ be a local, Noetherian ring, with residue field $k$, or a homogeneous $k$-algebra, with $k$ a field. Write ${\cal M}(R)$ for the category of finitely generated (and graded if $R$ is homogeneous) $R$-modules. 

It is of interest to assign a number, $\delta(M)$, to every $M$ in  ${\cal M}(R)$, which measures the complexity of the module $M$.  Some examples are the rank of $M$; the size of a minimal generating set for $M$, $\nu(M)$; the multiplicity of $M$; the arithmetic and geometric degrees, $\adeg(M)$ and $\gdeg(M)$, studied in Bayer-Mumford \cite{BM} and Sturmfels-Trung-Vogel \cite{STV}; the homological degree, $\hdeg(M)$, introduced in \cite{V96}, and other cohomological degrees.  For a survey, see the notes of Vas\-con\-celos \cite{VBarc}.
In the case where $M$ is Cohen-Macaulay, the multiplicity (or degree) of $M$, written $\deg(M)$, serves this purpose admirably.  It is defined as follows:  If $\dim(M)=0$, then the length of $M$, $\length(M)$, is finite and we let $\deg(M) = \length(M)$.  If $\dim(M)=d>0$ then $M$ has infinite length and we let 
$$\deg(M) = \lim_{n\to\infty}\frac{d!\cdot \length(M/\m^n M)}{n^d}.$$
$\deg(M)$ is always a positive integer unless $M=0$, in which case $\deg(M)=0$.

The multiplicity function has served particularly well in bounding numbers of generators of ideals and modules.  For example, Sally \cite{Sally76},\cite{Sally}, has shown that if $R$ is Cohen-Macaulay and $I$ is a Cohen-Macaulay ideal of height $c$, then the number of generators of $I$, $\nu(I)$, is bounded as follows:
\begin{equation}\label{sallybound} \nu(I) \leq \deg(R/I)^{c-1}\deg(R)+c-1.\end{equation}
Many similar bounds can be found in the literature (\cite{Becker, BER, ERV, Shalev, Trung, Valla}).  In particular, an improvement to (\ref{sallybound}) can be found in \cite{DGV}. On the other hand if $M$ is Cohen-Macaulay, then we have 
\begin{equation}\label{nuleqdeg} \nu(M) \leq \deg(M), \end{equation} and, as demonstrated by Brennan, Herzog and Ulrich in \cite{BHU}, equality holds if and only if $\grm(M)$ has a linear resolution (viewed as a module over a polynomial ring.)  Unfortunately, all of these results are false when the Cohen-Macaulay condition is removed.  In the words of Vasconcelos \cite[page 145]{VBarc}, ``$\deg(M)$ leaks all over as a predictor of properties of $M$'' when $M$ fails to be Cohen-Macaulay.  {\em Cohomological degree functions} were introduced in \cite{DGV} in an attempt to remedy this situation.

\begin{Def} Let $k$ be an infinite field.  A {\em cohomological degree} is a function, $\Deg(-)$, from the isomorphism classes of modules in ${\cal M}(R)$ to ${\bf R}$, satisfying the following conditions. \begin{enumerate} 
\item If $L = \gm(M)$ and $\overline{M} = M/L$, then 
$$ \Deg(M) = \Deg(\overline{M}) + \length(L).$$
\item If $\depth (M) > 0$ and $x$ is a generic hyperplane on $M$, then
$$ \Deg(M) \geq \Deg(M/xM).$$
\item If $M$ is Cohen-Macaulay, then 
$$\Deg(M) = \deg(M).$$ \end{enumerate} \end{Def}
($\gm(M)$ denotes the zeroth local cohomology module.  The condition that $k$ be infinite is required, of course, to make sense of the term 'generic'.)  The prototype for such a function is the homological degree, $\hdeg(-)$, introduced in \cite{V96}. 
	If $R$ is a Gorenstein ring of dimension $n$, then we define $\hdeg(M)$ inductively as follows.  If $\dim(M)=0$ then we let $\hdeg(M)=\length(M)$.  If $\dim(M)=d>0$, then we let 
$$\hdeg(M) = \deg(M) + \sum_{i=0}^{d-1}{d-1 \choose i}\hdeg(\Ext_R^{n-i}(M,R)).$$
(This definition is well founded since $\dim(\Ext_R^{n-i}(M,R))\leq i$ for all $i$.)  In \cite{V96} it is shown that $\hdeg(-)$ is a cohomological degree.

In Chapter 2 of this dissertation, we lay out the basic properties of cohomological degrees (many of which first appear in \cite{DGV}).  Included are several bounds on the number of generators of ideals in the spirit of Sally, but extended to the non-Cohen-Macaulay case (Section \ref{gen-props}).  Moreover, if $\Deg(-)$ is a cohomological degree and $M$ is a finitely generated module, we may generalize (\ref{nuleqdeg}) to 
$$\nu(M) \leq \Deg(M).$$
Yoshida \cite{Yoshida} has partially extended the result of Brennan-Herzog-Ulrich \cite{BHU} for the cohomological degree $\hdeg(-)$: 

\begin{The} If $\nu(M) = \hdeg(M)$ then $\grm(M)$ has a linear resolution. \end{The}

In Section \ref{Deg-linres} we show (following Yoshida's proof) that only the properties of cohomological degrees are needed.  Thus, we have:

\begin{The}\label{intro-Deglinres} Suppose $\Deg(-)$ is a cohomological degree and $M$ is a finitely generated module. If $\nu(M) = \Deg(M)$ then $\grm(M)$ has a linear resolution. \end{The}

However, unlike in the Cohen-Macaulay case, the converse does not hold.

In Chapter 3, we introduce a new cohomological degree, $\bdeg(-)$, which may be described axiomatically as follows.

\begin{Def} $\bdeg(-)$ is the unique function from the isomorphism classes of modules in ${\cal M}(R)$ to ${\bf N}$ satisfying the following properties.
\begin{enumerate}
\item $\bdeg(M) = \deg(M)$ if $M$ is Cohen-Macaulay.
\item If $L = \gm(M)$, then $\bdeg(M) = \bdeg(M/L) + \length(L)$.
\item If $\depth(M)>0$ and $x$ is a generic hyperplane on $M$, then $\bdeg(M) = \bdeg(M/xM)$. \end{enumerate} \end{Def}

It is not obvious that such a function exists.  Thus Section \ref{bdeg-axiom} is devoted to establishing that there is such a function (Theorem \ref{b-prop}).

$\bdeg(-)$ is sharper than $\hdeg(-)$, so the bounds in Section \ref{gen-props} are necessarily improved, but we also now have a converse to Theorem \ref{intro-Deglinres}.

\begin{The}\label{intro-bdeglinres} If $M$ is a finitely generated $R$ module then $\nu(M) = \bdeg(M)$ if and only if $\grm(M)$ has a linear resolution. \end{The}

A key step will be to show that $\bdeg(-)$ may only increase when we pass to the associated graded module.

\begin{The} $\bdeg(M) \leq \bdeg(\grm(M))$. \end{The}

In counterpoint to these degree-like functions, when $R$ is homogeneous, is the Castel\-nuovo-Mumford regularity of $M, \: \reg(M)$.  Since $R$ is homogeneous, it is the homomorphic image of a polynomial ring, $S$, and so we may view $M$ as an $S$-module.  $M$ has a graded minimal free resolution as an $S$-module:
$$0 \lar \bigoplus_j S[-j]^{\beta_{nj}} \lar \cdots \lar \bigoplus_j S[-j]^{\beta_{0j}} \lar M \lar 0.$$  
The regularity of $M$ is:
$$ \reg(M) := \max\{\, j-i \: | \: \beta_{ij}\neq 0 \, \}.$$
Regularity has a distinctly different flavor from the degree functions:  it is oblivious to rank ($\reg(M\oplus M) = \reg(M)$) but sensitive to shifts ($\reg(M[-37]) = \reg(M) + 37$), whereas a degree function, $\delta(-)$, tends to exhibit the opposite behavior;  it is sensitive to rank ($\delta(M\oplus M) = 2\delta(M)$) but oblivious to shifts ($\delta(M[-37]) = \delta(M)$).  This situation may be summarized by saying that degrees are ``ranky'' and regularity is ``shifty''.  Nonetheless, for algebras these two kinds of functions are not completely independent.  There is a growing literature (\cite{DGV,Hoa,HM,HT,MV1,MV2}) on the interplay between $\reg(R)$ and $\delta(R)$ for various $\delta(-)$.

For cohomological degrees, the starting point is the following theorem, taken from \cite{DGV}.

\begin{The} Suppose $R$ is a homogeneous $k$-algebra and $\Deg(-)$ is a cohomological degree on ${\cal M}(R)$. Then $$\reg(R) < \Deg(R).$$ \end{The}

We take up the question of how much $\Deg(R)$ may exceed $\reg(R)$.  Our first result is:

\begin{The} There exists a function $f(n,r)$ such that 
$$\hdeg(R) \leq f(\embdim(R),\reg(R)),$$
for all homogeneous $k$-algebras $R$ of dimension at most $2$. \end{The}

Here, $\embdim(R) = \nu(\m)$ is the embedding dimension of $R$.  We conjecture that the condition $\dim(R) \leq 2$ is unnecessary.

For $\bdeg(-)$ our results are stronger.  We give two bounds.  The first one is a consequence of Theorem \ref{intro-bdeglinres}.

\begin{The} If $M$ is a finitely generated graded module over a homogeneous $k$ algebra and $\depth(M)>0$, then
$$\bdeg(M) \leq H(M;r),$$
where $H(M;-)$ is the Hilbert function of $M$ and $r = \reg(M)$. \end{The}

Our second bound is less elegant but more concrete.  It is based on the fact that $\bdeg(R)$ does not change when we pass to the generic initial ideal.

\begin{The} Let $S=k[x_1,\ldots,x_n]$ be a polynomial ring over $k$ and let $I$ be a homogeneous ideal of $S$.  Then
$$\bdeg(S/I) = \bdeg(S/\gin(I)),$$
where $\gin(I)$ is the generic initial ideal of $I$. \end{The}

We give a combinatorial description of $\bdeg(S/\gin(I))$ and from this derive a sharp, explicit bound on $\bdeg(R)$, where $R$ is a homogeneous $k$-algebra.

Recall that any positive integer $e$ can be written in the form
$$e = {k(r) \choose r} + \cdots + {k(t) \choose t},$$
where $k(r) > \cdots > k(t) \geq t$, called the {\em $r^{th}$ Macaulay representation} of $e$.  For fixed positive integers $e$ and $r$, the numbers 
$k(r),\ldots,k(t)$ are uniquely determined.  Let
$$e^{(r,d)} = {k(r) + d \choose r} + \cdots + {k(t)+ d \choose t}.$$

\begin{The}Let $R$ be a homogeneous algebra over a field of characteristic $0$.  If $\deg(R) = e,\, \reg(R) = r,\, \dim R = d$, and $\depth (R) = g$, then $$\bdeg(R) \leq e + {n-g+r \choose r} - e^{(r,d-g)}. $$ \end{The}

In Chapter 4, we study another measurement on modules, which is also useful in bounding numbers of generators, the Dilworth number.

\begin{Def} If $M$ is in ${\cal M}(R)$, then we define the Dilworth number, $\dil(M)$, as follows:  $$\dil(M) = \max\{\,\nu(N)\;|\; N \subseteq M\,\}.$$
\end{Def}

Trivially, we have $\nu(M) \leq \dil(M)$, but unfortunately $\dil(M)$ is finite if and only if $M$ has dimension at most $1$.  Further, it simply does not hold that if $\nu(M) = \dil(M)$ then $\grm(M)$ has a linear resolution, even in dimension $0$.  Nevertheless, there is an analogy between the Dilworth number and cohomological degrees -- the homological degree in particular.  It is as follows.  If $R$ is Gorenstein of dimension $n$ and $M$ is an $R$-module of dimension $1$, then
$$\hdeg(M) = \deg(M) + \length(\Ext_R^n(M,R)).$$
By a theorem of Ikeda \cite{Ik2}, we have $$\dil(M) = \deg(M) + \dil(\Ext_R^n(M,R)).$$
Thus we are led to the definition of the homological Dilworth number, $\hdil(-)$.

\begin{Def}If $R$ is a Gorenstein ring of dimension $n$ and
$\dim M = d$, then the homological Dilworth number is defined recursively as follows,
$$ \hdil (M) = \left\{ \begin{array}{ll}
                      \deg (M) + \displaystyle\sum_{i=0}^{d-1}{d-1 \choose i}\hdil(\Ext_R^{n-i}(M,R)) & \mbox{if } d>0 \\
                      \dil(M) & \mbox{if } d=0 \end{array} \right. . $$ 
\end{Def}

We show that, like the homological degree, the homological Dilworth number may only decrease modulo a generic regular hyperplane section.

\begin{The} Let $M$ in ${\cal M}(R)$ have positive depth and let $x$ be a generic hyperplane on $M$.  Then $$\hdil(M/xM) \leq \hdil(M).$$ \end{The}

This allows us to prove Sally bounds analogous to those proved for cohomological degrees.  Since $\hdil(M) \leq \hdeg(M)$, the bounds are sharper for the homological Dilworth number than they are for the homological degree, though perhaps not for $\bdeg(-)$.  We also have another linear resolution theorem.

\begin{The} If $M$ is in ${\cal M}(R)$, then $\nu(M) \leq \hdil(M)$.  If $\nu(M)=\hdil(M)$ and $\m \gm(M) = 0$ then $\grm(M)$ has a linear resolution. \end{The}

In the final section, we compare $\hdil(-)$ to $\bdeg(-)$.  Our first result is:

\begin{The} If $M$ in ${\cal M}(R)$ has dimension at most $2$, then $\hdil(M) \leq \bdeg(M)$. \end{The}

As a corollary, we have:

\begin{The} If $M \in {\cal M}(R)$ has dimension at most $2$, then $\grm(M)$ has a linear resolution if and only if $\nu(M) = \hdil(M)$ and $\m \gm(M) = 0$. \end{The}

We show by example that this theorem is false in dimension $3$ and higher.  We also give examples to show that either inequality \begin{eqnarray*}
\bdeg(M) & > & \hdil(M) \\
\bdeg(M) & < & \hdil(M) \end{eqnarray*}
is possible when $\dim(M) \geq 3$.

%% file: thesis-chap1.tex
\chapter{Preliminaries}

In this chapter we provide a brief summary of the definitions and results in commutative algebra upon which we will rely in the subsequent chapters.  Our sources for most of this material are: for basic commutative algebra, \cite{Eisenbud} and \cite{Matsumura}; for the theory of Cohen-Macaulay rings and modules, \cite{BH}; for Buchsbaum rings and modules, \cite{SV}; and for homological algebra, \cite{Weibel}.  The reader is advised to skim this chapter briefly, coming back to it only as needed.

\section{Local and Homogeneous Rings; Hilbert Functions}

In this dissertation we will be concerned primarily with commutative rings $R$ (with identity, of course!) which are of either of the following types
\begin{itemize}
\item Noetherian local rings; or
\item homogeneous $k$-algebras, $k$ a field.
\end{itemize}
A homogeneous $k$-algebra (also called a standard graded algebra) is a positively graded ring $$R = R_{\langle 0 \rangle}\oplus R_{\langle 1 \rangle}\oplus R_{\langle 2 \rangle}\oplus \cdots,$$ where $R_{\langle i \rangle}$ denotes the $i^{th}$ graded component of $R$, such that $R_{\langle 0 \rangle} = k$, each $R_{\langle i \rangle}$ is a finite dimensional vector space over $k$, and $R$ is generated as a $k$-algebra in degree $1$:  $R=R_{\langle 0 \rangle}[R_{\langle 1 \rangle}]$.  Equivalently, $R$ is of the form $R=S/I$, where $S=k[x_1,\ldots,x_n]$ is a polynomial ring over $k$ and $I$ is a homogeneous ideal.

If $R$ is local, then we let $\m$ denote the maximal ideal and we let $k$ denote the residue class field $R/ \m$.  ${\cal M}(R)$ will denote the category of finitely generated $R$-modules.  If $R$ is homogeneous, then $\m$ will denote the irrelevant maximal ideal, $\m = R_{\langle 1 \rangle}\oplus R_{\langle 2 \rangle}\oplus\cdots$, and ${\cal M}(R)$ will denote the category of finitely generated graded $R$-modules.  For the proofs of the assertions in the following discussion, we refer the reader to \cite{BH}, \cite{Eisenbud} or \cite{Matsumura}.

If $R$ is homogeneous and $M = \oplus M_{\langle i \rangle}$ is in  ${\cal M}(R)$, then the Hilbert function of $M$ is $$H(M;n) = \dim_k M_{\langle n \rangle}.$$ The Hilbert function agrees with a polynomial for large enough values of $n$:
$$H(M;n) = \alpha \cdot n^{d-1}+ \mbox{ terms of lower order.}$$ Here, $d=\dim M$ and $(d-1)!\cdot \alpha$ is an integer, called the {\em multiplicity} or {\em degree} of $M$, written $\deg(M)$.  (If $\dim M = 0$, then we set $\deg(M) = \length(M)$.)

If $R$ is local and $M$ is in ${\cal M}(R)$, then the associated graded module
$$\grm(M) = \frac{M}{\m M} \oplus \frac{\m M}{\m^2 M} \oplus \frac{\m^2 M}{\m^3 M} \oplus \cdots$$ is a finitely generated graded module over the associated graded ring $\grm(R) = (R/\m) \oplus (\m / \m^2) \oplus \cdots$, which is a homogeneous $k$-algebra.  Thus $\grm(M)$ has a Hilbert function, and we define the degree of $M$ to be the degree of its associated graded module: $$\deg(M) = \deg(\grm(M)).$$

For any $\m$-primary ideal $I$, we may define the {\em Hilbert-Samuel} function of $M$ with respect to $I$:
$$\chi^I(M;n) = \length(I^nM/I^{t+1}M).$$
The Hilbert-Samuel function also agrees with a polynomial for $n$ large enough:
$$\chi^I(M;n) = \alpha\cdot n^{d-1} + \mbox{ terms of lower order, }$$
where $d=\dim M$ and $(d-1)!\cdot \alpha$ is an integer, called the {\em multiplicity of $M$ with respect to $I$}, and written $e(I;M)$.  Observe that $\deg(M) = e(\m;M)$.  For some interesting variations on this theme, we refer the reader to \cite{Doering}.

If $R$ is homogeneous and $M$ is in ${\cal M}(R)$, then the generating function  $$H_M(t) = \sum H(M;n)t^n$$ of the Hilbert function is called the Hilbert series.  The Hilbert series is a rational function of $t$, and can be written in
the form $$H_M(t) = \frac{Q_M(t)}{(1-t)^d},$$ where $d = \dim M$ and $Q_M(t) \in {\bf Z}[t,t^{-1}]$ satisfies $Q_M(1) = \deg(M)$.

\section{The Multiplicity Symbol}

Let $(R,\m)$ denote a Noetherian local ring or a homogeneous $k$-algebra.  For the proofs of the propositions in the section, see \cite[Section 4.6]{BH}

\begin{Definition} A sequence of elements ${\bf x}=x_1\ldots,x_t$ in $\m$ is called a {\em multiplicity system} for a module $M$ if $M/{\bf x}M$ has finite length. \end{Definition}

\begin{Definition} If ${\bf x}=x_1\ldots,x_t$ is a multiplicity system for $M$ then we define the {\em multiplicity symbol}, $e({\bf x};M)$, by induction on $t$, as follows.  If $t=0$, then $$e({\bf x};M) = e(\emptyset;M) = \length(M).$$  If $t>0$, then $$e({\bf x};M) = e(x_2,\ldots,x_t;M/x_1M)- e(x_2,\ldots,x_t;[0:x_1]_M).$$ \end{Definition}

\begin{Proposition} Let ${\bf x} = x_1,\ldots,x_t$ be a multiplicity system for $M$.  If $t>\dim(M)$, then $e({\bf x};M)=0$.  If $t=\dim(M)$, then the multiplicity symbol $e({\bf x};M)$ agrees with the multiplicity of the ideal generated by ${\bf x}$ with respect to $M$: $e({\bf x};M) = e(({\bf x});M).$ \end{Proposition}

\begin{Proposition} Given an exact sequence of modules $$0 \lar A \lar B \lar C \lar 0,$$ a sequence ${\bf x}$ of elements in $\m$ is a multiplicity system for $B$ if and only if it is a multiplicity system for $A$ and for $C$, and we have
$$e({\bf x};B) = e({\bf x};A) + e({\bf x};C).$$ \end{Proposition}

\section{Reductions}

Again, $(R,\m)$ is a Noetherian local ring or a homogeneous $k$-algebra.  If $R$ is homogeneous, then whenever we refer to an ideal of $R$, it is assumed that the ideal is homogeneous.

\begin{Definition} If $\mathfrak{a} \subseteq \mathfrak{b}$ are ideals of $R$ and $M$ is in ${\cal M}(R)$, then $\mathfrak{a}$ is said to be a {\em reduction of $\mathfrak{b}$ with respect to $M$} if $\mathfrak{a}\mathfrak{b}^t M = \mathfrak{b}^{t+1}M$ for $t \gg 0$.  We say that $\mathfrak{a}$ is a {\em reduction of $\mathfrak{b}$} if $\mathfrak{a}$ is a reduction of $\mathfrak{b}$ with respect to $R$.  If $\mathfrak{a}$ is a reduction of $\mathfrak{b}$ and no other reduction of $\mathfrak{b}$ is contained in $\mathfrak{a}$, then $\mathfrak{a}$ is said to be a {\em minimal reduction of $\mathfrak{b}$}. \end{Definition}

\begin{Proposition} Every ideal has a minimal reduction. \end{Proposition}

For a proof, see \cite{NR}.

\begin{Definition} If $\mathfrak{a} \subseteq \mathfrak{b}$ is a reduction, then the {\em reduction number of $\mathfrak{b}$ with respect to $\mathfrak{a}$}, written $r_{\mathfrak{a}}(\mathfrak{b})$, is given by
$$r_{\mathfrak{a}}(\mathfrak{b}) = \min\{\,t\;|\; \mathfrak{a} \mathfrak{b}^t = \mathfrak{b}^{t+1}\,\}.$$
The {\em reduction number of $\mathfrak{b}$}, denoted $r(\mathfrak{b})$, is given by
$$r(\mathfrak{b}) = \min \{\,r_{\mathfrak{a}}(\mathfrak{b})\;|\; \mathfrak{a} \mbox{ is a minimal reduction of } \mathfrak{b}\,\}.$$ \end{Definition}

For a proof of the following proposition, we refer the reader to \cite[Lemma 4.5.5.]{BH} or \cite[Theorem 14.13]{Matsumura}.

\begin{Proposition} If $\mathfrak{a} \subseteq \mathfrak{b}$ is a reduction with respect to $M$ and $\mathfrak{b}$ is $\m$-primary, then $e(\mathfrak{a};M) = e(\mathfrak{b};M)$. \end{Proposition}

The following proposition says, roughly, that any sequence of sufficiently generic linear combinations of the generators of an $\m$-primary ideal form a minimal reduction of the ideal.

\begin{Proposition} Suppose $k$ is infinite and $\dim M = d$.  Let $I$ be an $\m$-primary ideal, and suppose that $I=(y_1,\ldots,y_m)$, where $m=\nu(I)$.  Then there is a Zariski-open subset $U$ of the affine space $k^{md}$ such that if $\alpha_{1,1},\ldots,\alpha_{m,d} \in R$ and $(\overline{\alpha}_{1,1},\ldots,\overline{\alpha}_{m,d}) \in U$, where $\overline{\alpha}_{i,j}$ denotes the residue class of $\alpha_{i,j}$ in $R/\m = k$, and $$x_j = \sum_1^m \alpha_{i,j}y_i, \qquad \mbox{ for } 1 \leq j \leq d$$ then $({\bf x}) = (x_1,\ldots,x_d)$ is a minimal reduction of $I$ with respect to $M$. In particular, ${\bf x}$ is a system of parameters for $M$. \end{Proposition}

A proof is contained in \cite[Theorem 14.14]{Matsumura}.

\begin{Corollary} Suppose the field $k$ is infinite.  If $M$ is in ${\cal M}(R)$, there is a multiplicity system ${\bf x}$ for $M$ such that $\deg(M) = e(({\bf x});M)$. \end{Corollary}

\begin{Definition} If ${\bf x}$ is a multiplicity system for $M$ such that $e({\bf x};M) = \deg(M)$, then ${\bf x}$ is called a {\em degree system} for $M$. \end{Definition}

\begin{Proposition} \label{lastone?} If $x$ is part of a multiplicity system for $M$, then $\deg(M/x^rM) = r \cdot \deg(M)$, for all $r \geq 1$. \end{Proposition}

\proof This follows quickly from, for example, Lech's formula \cite[Theorem 14.12]{Matsumura}. \QED

\section{Gorenstein Rings}

Let $(R,\m)$ be a Noetherian local ring or a homogeneous $k$-algebra, and let $M$ be in ${\cal M}(R)$.  Then we have (\cite[Theorem 1.2.8]{BH}) $$\depth(M) = \min\{\, t \; | \; \Ext_R^t(k,M) \neq 0\,\}.$$  The modules $\Ext_R^t(k,M)$ are all annihilated by $\m$, so it makes sense to speak of the dimension of $\Ext_R^t(k,M)$ as a vector space over $k$.  This leads to the definition of the {\em type} of a module.

\begin{Definition}\label{type-def} If $\depth(M) = t$ then the {\em type} of $M$ is the integer $$\type(M) = \dim_k \Ext_R^t(k,M).$$ \end{Definition}

Type is preserved modulo a regular element.

\begin{Proposition}\label{typelemma} Let $x \in \m$ be regular on $M$.  Then $$\type(M) = \type(M/xM).$$
\end{Proposition}

\proof Let $\depth(M)=t$.  We will show $\Ext_R^t(k,M) \cong \Ext_R^{t-1}(k,M/xM).$
From the exact sequence
$$0 \lar M \stackrel{x}{\lar} M \lar M/xM \lar 0$$
we obtain the exact sequence
$$0\lar \Ext_R^{t-1}(k,M/xM) \lar  \Ext_R^t(k,M)\stackrel{x}{\lar}\Ext_R^t(k,M).$$
Since $x$ annihilates $\Ext_R^t(k,M)$, this gives the desired isomorphism. \QED

\begin{Definition} The {\em socle} of $M$ is given by $$\socle(M) = \Hom_R(k,M) = [0:\m]_M.$$ \end{Definition}

Thus, the socle is nonzero precisely when the depth is zero.

\begin{Definition} $(R,\m)$ is called {\em Gorenstein} if it is Cohen-Macaulay and $\type(R)=1$. \end{Definition}

Thus, by Proposition \ref{typelemma}, if $R$ is Cohen-Macaulay and ${\bf x}$ is a maximal regular sequence in $\m$, then $R$ is Gorenstein if and only if $\socle(R/{\bf x}R) \cong k$.  In particular, if $S = k[x_1,\ldots,x_n]$ or if $S$ is a regular local ring, then $S$ is Gorenstein.

We will frequently be moved to require that $(R,\m)$ be the homomorphic image of a Gorenstein ring, in order to use the duality described in the next section.  If $R$ is homogeneous, this is automatically so, since $R$ is the quotient of a polynomial ring by a homogeneous ideal, $R=S/I$.  If $R$ is local, it is not necessarily the homomorphic image of a Gorenstein ring, but by Cohen's structure theorem (\cite[Theorem A.21]{BH} or \cite[Theorem 29.4]{Matsumura}), its $\m$-adic completion $\hat{R}$ is.

\section{Local Cohomology and Local Duality}\label{intro-local-duality}

In this section, we preserve the notation of the previous sections.  We refer the reader to \cite[Sections 3.2 and 3.5]{BH} for the proofs of the assertions in this section.

\begin{Definition} If $M$ is in ${\cal M}(R)$ and $M \hookrightarrow E$ is an embedding of $M$ into an injective module $E$ such that for every embedding $M \hookrightarrow I$ of $M$ into an injective module $I$, there is an embedding $E \hookrightarrow I$ of $E$ into $I$ such that the following diagram commutes
$$ \begin{array}{ccc}
M & \rightarrow & E \\
\downarrow & \swarrow & \\
I & & \end{array}$$
then $E$ is said to be the {\em injective envelope} of $M$. \end{Definition}

Injective envelopes always exist, and are unique up to isomorphisms which respect the embedding $M \hookrightarrow E$.  We let $E_R(M)$ denote the injective envelope of $M$.  For $M$ in ${\cal M}(R)$, we write $M^{\vee} = \Hom_R(M,E_R(k))$.  $M^{\vee}$ is called the {\em Matlis dual} of $M$.  For the following proposition, which is a summary of the facts about Matlis duality that we will require, we let ${\cal M}_0(R)$ denote the category of $R$-modules of finite length.

\begin{Proposition}\label{matlis-duality} {\em (Matlis Duality)} $(-)^{\vee}$ is an exact functor on ${\cal M}_0(R)$ satisfying, for all $M$ in ${\cal M}_0(R)$, 
\begin{enumerate}
\item $\length(M) = \length(M^{\vee})$,
\item $\nu(M) = \type(M^{\vee})$,
\item $\type(M) = \nu(M^{\vee})$, and
\item The natural homomorphism $M \rightarrow (M^{\vee})^{\vee}$ is an isomorphism.
\end{enumerate}
\end{Proposition}

Here, $\length(-)$ denotes the length function and $\nu(-)$ counts the number of generators in a minimal generating set.

\begin{Definition} For $M$ in ${\cal M}(R)$, we define the submodule $\gm(M)$ to be 
$$\gm(M) = \{\, z \in M \; | \; \m^t z = 0 \mbox{ for } t \gg 0\,\}.$$
\end{Definition}

\begin{Proposition} $\gm(-)$ is an additive, left exact functor on ${\cal M}(R)$. \end{Proposition}

\begin{Definition} The $i^{th}$ right derived functor of $\gm(-)$ is written $H_{\m}^i(-)$.  $H_{\m}^i(M)$ is called the $i^{th}$ {\em local cohomology module} of $M$, for $M$ in ${\cal M}(R)$. \end{Definition}

\begin{Proposition}\label{depth-dim} Local cohomology detects depth and dimension:
\begin{eqnarray*} \depth(M) & = & \min \{\, i \; | \; H_{\m}^i(M) \neq 0 \,\} \\
\dim M & = & \max \{\, i \; | \; H_{\m}^i(M) \neq 0 \,\} \end{eqnarray*}
\end{Proposition}

Suppose $R$ is the homomorphic image of a Gorenstein ring, $S$.  For all $M$ in ${\cal M}(R)$, we write
$$M_i = \Ext_S^{n-i}(M,S),$$
where $n = \dim S$.  $M_i$ is an $R$-module, and is independent of the choice of $S$.

\begin{Proposition} \label{local-duality}{\em (Local Duality)} Suppose $R$ is the homomorphic image of a Gorenstein ring.  Then for $M$ in ${\cal M}(R)$, we have
$$H_{\m}^i(M) \cong (M_i)^{\vee},$$ for all $i$.   In particular, $M_i = 0$ if $i>\dim M$. \end{Proposition}

We remark that the module $M_i$ need not be in ${\cal M}_0(R)$, so Proposition \ref{matlis-duality} does not necessarily apply.
Combining Propositions \ref{depth-dim} and \ref{local-duality}, we obtain the following corollary.

\begin{Corollary} \label{smiley} If $R$ is the homomorphic image of a Gorenstein ring and $M$ is in ${\cal M}(R)$, then $M$ is Cohen-Macaulay if and only if $M_i = 0$ for $i<\dim(M)$. \end{Corollary}

\begin{Proposition}\label{h-Spec} If $R$ is the homomorphic image of a Gorenstein ring and if $M$ is in ${\cal M}(R)$, then $\dim M_i \leq i.$ Furthermore, for all $P \in \Spec(R)$ of dimension $i$, we have that $P \in \Ass(M)$ if and only if $P \in \Ass(M_i)$. \end{Proposition}

\section{Buchsbaum Modules}

In this section, $(R,\m)$ will denote a Noetherian local ring or a homogeneous $k$-algebra.  The proof of the assertions can be found in \cite{SV}.

\begin{Definition} A module $M$ in ${\cal M}(R)$ is {\em Buchsbaum} if $\dim M = 0$, or if $\dim M > 0$ and there exists an integer $I(M) \geq 0$ such that for every system of parameters ${\bf x}$ for $M$, we have $$\length(M/{\bf x}M) - e(({\bf x});M) = I(M).$$ If $\dim M = 0$, then we set $I(M)=0$. \end{Definition}

$I(M)$ is sometimes referred to as the invariant of St\"{u}ckrad and Vogel.  Note that $M$ in ${\cal M}(R)$ is Cohen-Macaulay if and only if it is Buchsbaum with $I(M)=0$.

\begin{Proposition} \label{Buch-implies-qBuch} If $M$ in ${\cal M}(R)$ is Buchsbaum, then $\m H_{\m}^i(M)=0$ for $i < \dim M$.  \end{Proposition}

The converse to Proposition \ref{Buch-implies-qBuch} is not true.  If $\m H_{\m}^i(M)=0$ for $i < \dim M$, then we say that $M$ is {\em quasi-Buchsbaum}.

\begin{Proposition} If $M$ in ${\cal M}(R)$ is Buchsbaum of dimension $d$, then
$$I(M) = \sum_{i=0}^{d-1}{d-1 \choose i} \dim_k H_{\m}^i(M).$$
\end{Proposition}

\begin{Definition} If $\length(H_{\m}^i(M)) < \infty$ for $i < d = \dim M$, then we say that $M$ has {\em finite length cohomology}, and we set $$I(M) = \sum_{i=0}^{d-1} {d-1 \choose i} \length(H_{\m}^i(M)).$$
\end{Definition}

\section{Castelnuovo-Mumford Regularity}

In this section, $R$ will denote a homogeneous $k$-algebra of the form $R=S/I$, where $S=k[x_1,\ldots,x_n]$ is a polynomial ring in $n$ variables over the field $k$ and $I$ is a homogeneous ideal.  For the proofs of the results in this section, we refer the reader to \cite{Eisenbud} and \cite{EG}.

For $M$ in ${\cal M}(R)$, we let $M[i]$ denote the module $M$ shifted by $i$:
$$M[i]_{\langle j \rangle} = M_{\langle i + j \rangle}.$$
Every $M$ in ${\cal M}(R)$ can be viewed as a graded $S$-module, and as such has a minimal graded free resolution over $S$:
$$0 \lar \bigoplus_{j \geq s+n}\left(S[-j]\right)^{\beta_{nj}} \lar \cdots \lar \bigoplus_{j \geq s+1}(S[-j])^{\beta_{1j}} \lar \bigoplus_{j\geq s}(S[-j])^{\beta_{0j}}\lar M \lar 0.$$
The $\beta_{ij}$ are called the {\em graded Betti numbers} of $M$.  We let $\beta_i = \sum_j \beta_{ij}$ and call $\beta_i$ the $i^{th}$ {\em Betti number} of $M$.  Clearly, $\beta_{ij}$ and $\beta_i$ depend on the choice of $S$, but the definitions and propositions which follow are independent of the choice of $S$.

\begin{Proposition}\label{notsohard} If $M$ in ${\cal M}(R)$ has graded Betti numbers $\beta_{ij}$, then the Hilbert series of $M$ is given by
$$H_M(t) = \frac{\sum_{i,j}(-1)^i\beta_{ij}t^j}{(1-t)^n}.$$ \end{Proposition}

For a proof, see \cite[Lemma 4.1.13]{BH}.

\begin{Definition} We say that $M$ has a {\em linear resolution} if $\beta_{ij} = 0$ for $j\neq i$; that is, if the minimal graded free resolution of $M$ is of the form
$$0 \lar (S[-n])^{\beta_n} \lar \cdots \lar (S[-1])^{\beta_1} \lar S^{\beta_0} \lar M \lar 0.$$ \end{Definition}

\begin{Definition} If $M=\oplus_i M_{\langle i \rangle}$ is in ${\cal M}(R)$, then for $r \in {\bf Z}$, we let $M_{\geq r}$ denote the graded submodule of $M$ given by $$M_{\geq r} = \bigoplus_{i \geq r} M_{\langle i \rangle}.$$\end{Definition}

\begin{Proposition} \label{trunc-prop} For any $M$ in ${\cal M}(R)$, $M_{\geq r}[r]$ has a linear resolution for $r \gg 0$.   Furthermore, if $M_{\geq r}[r]$ has a linear resolution, then $M_{\geq s}[s]$ has a linear resolution for all $s \geq r$.\end{Proposition}

\begin{Definition} The {\em Castelnuovo-Mumford regularity} of $M$ is the integer $\reg(M)$, given by $$\reg(M) = \min\{\, r \;|\; M_{\geq r}[r] \mbox{ has a linear resolution}\,\}.$$ \end{Definition}

If $\length(M) < \infty$, then $\reg(M) = \max\{ \, i \;|\; M_{\langle i \rangle} \neq 0\,\}$.

\begin{Proposition} Let $M$ be in ${\cal M}(R)$ and let $\beta_{ij}$ be the graded Betti numbers of $M$.  Then \begin{eqnarray*}
\reg(M) & = & \max\{\, j-i \;|\; \beta_{ij}\neq 0\,\} \\
& = & \max\{\, j-i \;|\; H_{\m}^i(M)_{\langle j \rangle} \neq 0 \,\} \end{eqnarray*} \end{Proposition}

\begin{Proposition} \label{reg-ses} Given a short exact sequence of modules
$$0 \lar A \lar B \lar C \lar 0,$$ we have the following inequalities.
\begin{enumerate}
\item $\reg(A) \leq \max\{\,\reg(B),\,\reg(C)+1\,\}$
\item $\reg(B) \leq \max\{\,\reg(A),\,\reg(C)\,\}$
\item $\reg(C) \leq \max\{\,\reg(A)-1,\,\reg(B)\,\}$
\end{enumerate}
\end{Proposition}

\begin{Proposition}\label{reg-inductive} If $M$ is in ${\cal M}(R)$ and $x \in R_{\langle 1 \rangle}$ is a linear form such that $\length([0:x]_M) < \infty$, then
$$\reg(M) = \max\{\,\reg([0:x]),\,\reg(M/xM)\,\}.$$
In particular, if $x \in R_{\langle 1 \rangle}$ is regular on $M$, then $\reg(M) = \reg(M/xM)$.
 \end{Proposition}

\section{Linear Cohen-Macaulay and Linear Buchsbaum Modules}\label{linCM-linBuchs}

In this section, $(R,\m)$ will denote a Noetherian local ring or a homogeneous $k$-algebra.

\begin{Proposition}\label{BHUtriv} If $M$ is Cohen-Macaulay, then $\nu(M) \leq \deg(M)$. \end{Proposition}

\begin{Definition} If $M$ is Cohen-Macaulay and $\nu(M)=\deg(M)$, then we say that $M$ is a {\em linear Cohen-Macaulay module}. \end{Definition}

\begin{Proposition}\label{BHU} Suppose the field $k$ is infinite.  If $M$ is Cohen-Macaulay, then $M$ is a linear Cohen-Macaulay module if and only if $\grm(M)$ has a linear resolution. \end{Proposition}

The proofs of Propositions \ref{BHUtriv} and \ref{BHU} can be found in \cite{BHU}.

\begin{Proposition}\label{FLC-bound} If $M$ has finite length cohomology, then $\nu(M) \leq \deg(M) + I(M).$ \end{Proposition}

\begin{Definition} If $M$ has finite length cohomology and $\nu(M) = \deg(M) + I(M)$, then we say that $M$ is a {\em linear Buchsbaum module}. \end{Definition}

\begin{Proposition}\label{linBuchs} Suppose the field $k$ is infinite.  If $M$ has finite length cohomology, then $M$ is a linear Buchsbaum module if and only if $\grm(M)$ is Buchsbaum and has a linear resolution.  Furthermore, if these equivalent conditions hold, then $M$ is Buchsbaum. \end{Proposition}

The proofs of Propositions \ref{FLC-bound} and \ref{linBuchs} can be found in \cite{Y2}.

\section{Generic Initial Ideals}\label{GIN-intro}

In this section, $S=k[x_1,\ldots,x_n]$ will denote a polynomial ring in $n$ variables over the field $k$ and $I$ will denote a homogeneous ideal.  We begin by describing the {\em reverse lexicographic order}.

Let $m$ be a monomial in $S$.  Then  we let $\deg_{x_i}(m) = \alpha_i$ and $\deg(m) = \sum \alpha_i$, where $m=\beta\cdot x_1^{\alpha_1}\cdots x_n^{\alpha_n}$, with $0 \neq \beta \in k$ and $0 \leq \alpha_i \in {\bf Z}$ for $1\leq i \leq n$.  If $m_1$ and $m_2$ are monomials in $S$, then we write $m_1 \succ_{rlex} m_2$ if $\deg(m_1) > \deg(m_2)$ or if $\deg(m_1) = \deg(m_2)$ and there is some $2 \leq j \leq n$ such that $\deg_{x_i}(m_1) = \deg_{x_i}(m_2)$ for $i>j$ and $\deg_{x_j}(m_1) < \deg_{x_j}(m_2)$.

If $f \in S$, then $f$ may be written uniquely as a sum of monomials $$f = m_1 + \cdots + m_t,$$ where $m_1 \succ_{rlex} m_2 \succ_{rlex} \cdots \succ_{rlex} m_t$.  The first monomial, $m_1$, is called the {\em leading term} of $f$, and is written $\init(f)$.  If $I$ is a homogeneous ideal, then the initial ideal of $I$, $\init(I)$, is given by $$\init(I) = ( \{\, \init(f)\;|\; f \in I\,\}).$$
For a proof of the following proposition, see \cite[Theorem 4.2.3 and Corollary 4.2.4]{BH}.

\begin{Proposition}\label{Hilbert-init} If $I \subseteq S$ is a homogeneous ideal, then $S/I$ and $S/\init(I)$ have the same Hilbert function.  Thus $\dim(S/I) = \dim(S/\init(I))$ and $\deg(S/I) = \deg(S/\init(I))$. \end{Proposition}

Let $V = S_{\langle 1 \rangle}$ be the vector space of linear forms.  The group of invertible linear transformations of $V$, $\Gl(V)$, acts on $S_{\langle 1 \rangle}$, and hence on all of $S$:  If $\sigma \in \Gl(V)$ and $f = f(x_1,\ldots,x_n) \in S$, then $\sigma \cdot f = f(\sigma\cdot x_1,\ldots,\sigma \cdot x_n) \in S$.  If $I$ is an ideal of $S$, let $\sigma \cdot I = \{\, \sigma\cdot f \;|\; f \in I\,\}$.  This is again a homogeneous ideal.  Thus, fixing an a homogeneous ideal $I$, we have a map
\begin{eqnarray*}
\Gl(S_1) & \lar & \{ \, \mbox{monomial ideals of } S \, \} \\
\sigma & \longmapsto & \init(\sigma \cdot I). \end{eqnarray*} 

\begin{Proposition} \label{Galligo} Suppose the field $k$ is infinite and fix a homogeneous ideal $I$ of $S$.  Then there is a nonempty Zariski-open subset $U$ of $\Gl(V)$ such that the map described above is constant on $U$; that is, for $\sigma, \tau \in U$ we have $\init(\sigma\cdot I) = \init(\tau \cdot I)$. \end{Proposition}

A proof can be found in \cite[Theorem 15.18]{Eisenbud}.

\begin{Definition} Let $I \subseteq S$ be a homogeneous ideal and suppose the field $k$ is infinite.  The unique monomial ideal described in Proposition \ref{Galligo} is called the {\em generic initial ideal}, and is denoted $\gin(I)$. \end{Definition}

For the remainder of this section, we require that the field $k$ have characteristic $0$.

\begin{Definition} A monomial ideal $I\subseteq S$ is said to be {\em Borel-fixed} if for all monomials $m \in I$, if $\deg_{x_i}(m) > 0$ for some $2 \leq i \leq n$, then $(x_{i-1}/x_i)m \in I$ as well. \end{Definition}

The following proposition is \cite[Theorem 15.20]{Eisenbud}.

\begin{Proposition} Let $I$ be a homogeneous ideal.  Then $\gin(I)$ is Borel-fixed. \end{Proposition}

The {\em saturation} of an ideal $I$ is defined by
$$I^{sat} = \bigcup_{t\geq 0}[I:\m^t];$$ thus $\gm(S/I) = I^{sat}/I$.
For proofs of the remaining propositions in the section, we refer the reader to the notes of Green \cite{Green}.

\begin{Proposition} \label{BF-summary} Let $I\subseteq S=k[x_1,\ldots,x_n]$ be a Borel-fixed ideal and let $R = S/I$. \begin{enumerate}
\item $I^{sat}=\cup_{t\geq 0}[I:x_n^t] = [I:x_n^{\infty}]$.
\item If $x$ is a generic hyperplane section, $R/xR \cong R/x_nR$ by a linear change of variables.
\item $\reg(S/I) = \max \{ \, \deg(m) \; | \; m \mbox{ is a generator of } I \, \} -1$.
\end{enumerate} \end{Proposition}

If $0 \neq h \in S_{\langle 1 \rangle}$, then we let $S_h = S/(h)$ and $I_h = (I+(h))/(h) \subseteq S_h$.  Thus $S_h$ is a polynomial ring in $n-1$ variables and $I_h$ is a homogeneous ideal.  Observe that the generic initial ideal of $I_h$ is well defined (that is, it does not depend on which linear forms in $S_h$ we decide to distinguish as the variables of $S_h$) since $\gin(I_h)$ is computed after making a generic linear change of variables.

\begin{Proposition}\label{gin-gen-hyp} There is a nonempty Zariski-open subset $U$ of $S_{\langle 1 \rangle}$ such that $$S_h/\gin(I_h) \cong S/(\gin(I) + (x_n))$$
for all $h \in U$. \end{Proposition}

\begin{Proposition}\label{gin-summary} Let $I\subseteq S=k[x_1,\ldots,x_n]$ be a homogeneous ideal.

\begin{enumerate}
\item $\gin(I^{sat}) = \gin(I)^{sat}$.
\item  $\depth (S/I) = \depth (S/\gin(I))$.
\item  $\reg (S/I) = \reg(S/\gin(I))$.
\end{enumerate} \end{Proposition}

\section{Superficial Elements}

\begin{Definition} An element $x \in \m \setminus \m^2$ is called {\em superficial} on $M$ if the map $$\grm(M) \stackrel{x^*}{\lar} \grm(M),$$ where $x^*$ is the initial form of $x$ in $\grm(R)$, has kernel of finite length.\end{Definition}

In other words, $x$ is superficial on $M$ if $$[\m^{k+1}M:x]\cap \m^{k-1}M=\m^{k}M$$ for $k \gg 0$.  As an application of the Artin-Rees theorem we obtain the following result.

\begin{Proposition}\label{superlemma} If $x$ is superficial on $M$, then
$$xM \cap \m^kM = x\m^{k-1}M$$ for $k \gg 0$. \end{Proposition}

\proof Since $x$ is superficial, we have 
\begin{equation}\label{sf}x\m^{k-1}M \cap \m^{k+1}M = x\m^kM\end{equation} 
for $k \gg 0$, and, by the Artin-Rees theorem, $$xM \cap \m^{k+i}M = \m^i(xM\cap \m^kM)$$ for $k \gg 0$.  Choose $k$ large enough to satisfy both equations.  Write $xM \cap \m^kM = xN$.  Then, for $i \geq 0$, we have
$$xM \cap \m^{2k+i}M  =  \m^{k+i}(xM\cap \m^kM) 
 =  x\m^{k+i}N \subseteq x\m^{k+i}M.$$
Thus
\begin{eqnarray*} 
xM \cap \m^{2k+i}M & = & x\m^{k+i}M\cap \m^{2k+i}M \\
& = & x\m^{k+i}M \cap \m^{k+i+2}M \cap \m^{k+i+3}M \cap \cdots \cap \m^{2k+i}M \\
& = & x\m^{k+i+1}M \cap \m^{k+i+3}M \cap \cdots \cap \m^{2k+i}M \\
& = & x\m^{k+i+2}M \cap \cdots \cap \m^{2k+i}M \\
& = & x\m^{2k+i-1}M, \end{eqnarray*}
by repeatedly applying equation (\ref{sf}). \QED

\bigskip

We say that ${\bf x} = x_1,\ldots,x_d$ is a {\em superficial sequence} for $M$ if each $x_i$ is a superficial sequence on $M$ and $x_1,\ldots,\hat{x}_i,\ldots,x_d$ is superficial on $M/x_iM$.

\begin{Proposition} Suppose $k$ is infinite and $\dim M = d$.  Suppose that $\m=(y_1,\ldots,y_m)$, where $m=\nu(\m)$.  Then there is a Zariski-open subset $U$ of the affine space $k^{md}$ such that if $\alpha_{1,1},\ldots,\alpha_{m,d} \in R$ and $(\overline{\alpha}_{1,1},\ldots,\overline{\alpha}_{m,d}) \in U$, where $\overline{\alpha}_{i,j}$ denotes the residue class of $\alpha_{i,j}$ in $R/\m = k$, and $$x_j = \sum_1^m \alpha_{i,j}y_i, \qquad \mbox{ for } 1 \leq j \leq d$$ then $({\bf x}) = (x_1,\ldots,x_d)$ is a superficial sequence with respect to $M$.\end{Proposition}

%% file: thesis-chap2.tex
\chapter{Cohomological Degrees}
\section{Definitions and Examples}

Let $(R,\m)$ be either a local, Noetherian ring, with residue field $k$, or else a homogeneous $k$-algebra. In either case, we assume that the field $k$ is infinite.

\begin{Definition} If $R$ is local, a {\em notion of genericity} on ${\cal M}(R)$ is a function
$$U(-): \{ \mbox{ isomorphism classes of }{\cal M}(R) \} \lar \{ \mbox{ non-empty subsets of }\m \setminus \m^2 \}$$ subject to the following conditions for each $M$.
\begin{enumerate}
\item If $f-g \in \m^2$ then $f \in U(M)$ if and only if $g \in U(M)$.
\item The set $\overline{U(M)} \subseteq \m / \m^2$ contains a nonempty Zariski-open subset.
\item If $\depth (M) > 0$, and $f \in U(M)$, then $f$ is regular on $M$.
\end{enumerate}
If $R$ is homogeneous, then a {\em notion of genericity} is a function
$$U(-): \{ \mbox{ isomorphism classes of }{\cal M}(R) \} \lar \{  \mbox{ non-empty subsets of }R_{\langle 1\rangle} \}$$ subject to the following conditions for each $M$.
\begin{enumerate}
\item The set $U(M) \subseteq R_{\langle1\rangle}$ contains a nonempty Zariski-open subset.
\item If $\depth (M) > 0$, and $f \in U(M)$, then $f$ is regular on $M$.
\end{enumerate} \end{Definition}

The set $U(M)$ should be thought of as the set of all ``generic hyperplanes'' on $M$.  For the following definition, we fix a notion of genericity $U$.

\begin{Definition}\label{Deg} A {\em cohomological degree} on ${\cal M}(R)$ is a numerical function 
$$\Deg(-):\{ \mbox{ isomorphism classes of }{\cal M}(R) \}\lar \bf{R}$$ which satisfies the following conditions for all $M$ in ${\cal M}(R)$.
\begin{enumerate}
\item If $L=\gm(M)$ and $\overline{M}=M/L$, then $$\Deg(M) = \Deg(\overline{M})+\length(L).$$
\item If $x\in U(M)$ and $\depth(M)>0$, then $$\Deg(M) \geq \Deg(M/xM).$$
\item {\em (The calibration rule)} If $M$ is Cohen-Macaulay, then $$\Deg(M) = \deg(M).$$ \end{enumerate}
\end{Definition}

Implicit in every cohomological degree, $\Deg(-)$, is a choice of a notion of genericity, $U(-)$.  When a cohomological degree has been specified, we will usually suppress the use of the symbol $U$; instead, we will write ``$x$ is generic on $M$'' for ``$x\in U(M)$''.  Similarly, a {\em generic system of parameters, ${\bf x} = x_1,\ldots,x_d$ for $M$} is a system of parameters for $M$ such that $x_1 \in U(M)$ and, for $i=2,\ldots,d, \; x_i\in U(M/(x_1,\ldots,x_{i-1})M)$.  When we wish to emphasize the relationship between $\Deg(-)$ and $U(-)$, we say that $\Deg(-)$ is a {\em cohomological degree with respect to $U(-)$.}

\begin{Example} The Homological Degree \end{Example}

The homological degree is the prototype for a cohomological degree.  Indeed, its definition \cite{V96} predates the definition of cohomological degree \cite{DGV}.

Suppose $R$ is the homomorphic image of a Gorenstein ring $S$ of dimension $n$.  If $R$ is homogeneous, then we take $S$ to be a polynomial ring, so there is no conflict of notation.  Write $M_i = \Ext_S^{n-i}(M,S)$.  If $M$ is a finitely generated $R$-module of dimension $d$, then we define the homological degree, $\hdeg(-)$, as follows,
\begin{equation}\label{hdeg-def} \hdeg(M) = \deg(M) + \sum_{i=0}^{d-1}{d-1 \choose i}\hdeg(M_i),\end{equation}
if $d>0$, and $\hdeg(M) = \length(M)$ if $d=0$.

If $R$ is not the homomorphic image of a Gorenstein ring, then let $\hat{R}$ be the $\m$-adic completion of $R$ and define
$$\hdeg(M) = \hdeg(M\otimes \hat{R}).$$

Equation (\ref{hdeg-def}) may at first appear strange, but in fact it is quite natural:  Suppose that $M$ has finite length cohomology.  Then Equation (\ref{hdeg-def}) becomes, by local duality (Proposition \ref{local-duality}), 
$$ \hdeg(M) = \deg(M) + \sum_{i=0}^{d-1}{d-1 \choose i}\length (H_{\m}^{i}(M)),$$ which is the well-known invariant of St\"{u}ckrad and Vogel.

\begin{Proposition}\label{hdeg} If $\dim M = 1$, let $U(M) = \{\, x\in \m \; | \; x \mbox{ is a degree system on } M\,\}$.  For modules $M$ of dimension $d>1$, let $$U(M) = \{\, x\in \m \; | \; x \mbox{ is part of a degree system on $M$ and $\length([0:x]_M)<\infty$} \,\} \cap \left( \bigcap_1^{d-1}U(M_i) \right).$$
Then $U(-)$ is a notion of genericity, and for $x \in U(M)$ and $r \geq 1$ we have
$$\hdeg(M/x^rM) \leq r \cdot \hdeg(M).$$
\end{Proposition}

\proof See \cite{V96} or the proof of Proposition \ref{inductive}. \QED

\begin{Corollary} $\hdeg(-)$ is a cohomological degree with respect to the notion of genericity described in Proposition \ref{hdeg}. \end{Corollary}

\proof Properties (1) and (3) of Definition \ref{Deg} follow easily from Proposition \ref{local-duality} and Corollary \ref{smiley}.  Condition (2) is just Proposition \ref{hdeg} with $r=1$. \QED

\bigskip

We remark that if $R$ is homogeneous, we may actually compute $\hdeg(-)$ on a computer algebra system with the capabilities of {\sc Macaulay} \cite{MAC}, and further that the definition of $\hdeg(-)$ does not require that the field be infinite.

\begin{Example} The Extremal Cohomological Degree \end{Example} Fix a notion of genericity $U$.  Define $\bdeg_U(-)$ as follows:
$$ \bdeg_U(M) = \min \{ \,\Deg(M)\; |\; \Deg(-) \mbox{ is a cohomological degree with respect to } U \}.$$  It is easy to verify, assuming there are any cohomological degree functions with respect to $U$, that $\bdeg_U(-)$ is a cohomological degree.  Note that, if $L = \gm(M)$, then  
\begin{equation}\label{bdeg-characterization} \label{bdegU}\bdeg_U(M)  = 
\left\{ \begin{array}{lcl}
\length(L) + \bdeg_U(M/L) & \mbox{ if } & L \neq 0 \\
\max \{ \,\bdeg_U(M/xM)\; |\; x \in U(M) \,\} & \mbox{ if } & L=0  \end{array}\right. .\end{equation}
This description is not very explicit -- it is not clear how to compute $\bdeg_U(-)$ from (\ref{bdegU}), since there are infinitely many hyperplanes to consider.  In Chapter 3, we will see how to obtain an explicit description of $\bdeg_U(-)$ -- one that is more or less amenable to explicit computations -- when the notion of genericity is selected appropriately.

\section{Elementary Properties of Cohomological Degrees}\label{Deg-elem}

Throughout this section, $(R,\m)$ is a Noetherian local ring with infinite residue field $k$, or a homogeneous algebra over an infinite field $k$, with irrelevant maximal ideal $\m$, and $\Deg(-)$ is a cohomological degree on ${\cal M}(R)$.

\begin{Proposition}\label{basics} For each $M$ in ${\cal M}(R)$ the following hold.
\begin{enumerate} \item $\deg(M) \leq \Deg(M)$, with equality if and only if $M$ is Cohen-Macaulay.
\item If $L\subseteq M$ is any submodule of finite length, then $\Deg(M) = \Deg(M/L) + \length(L)$.
\end{enumerate}
\end{Proposition}

\proof Trivial. \QED

\begin{Proposition} \label{M/xM} Let $M$ be in ${\cal M}(R)$ and let $L=\gm(M)$.  If $x$ is a generic hyperplane on $M$, then
$$\Deg(M) \geq \Deg(M/xM) + \length(L) - \length(L/xL).$$
In particular, $\Deg(M/xM) \leq \Deg(M)$.
\end{Proposition}

\proof See \cite[Proposition 2.3]{DGV}. \QED

\begin{Proposition}\label{regDeg} Suppose $R$ is homogeneous.  Then
$$\reg(R) < \Deg(R).$$
\end{Proposition}

\proof See \cite[Proposition 2.4]{DGV}. \QED

\begin{Proposition} \label{inequalities} For every cohomological degree $\Deg(-)$, the following inequalities hold.

\begin{enumerate}
\item $\nu(M) \leq \Deg(M)$ for all $M$ in ${\cal M}(R)$.
\item $\type(M) \leq \Deg(M)$ for all $M$ in ${\cal M}(R)$.
\item $\embcod(R) + 1 \leq \Deg(R)$.
\end{enumerate}
\end{Proposition}

Here, $\embcod(R)$ denotes the {\em embedding codimension} of $R$, and is given by $\embcod(R) = \embdim(R)-\dim(R) = \nu(\m) - \dim(R)$.

\bigskip

\proof 

(1) is proved in \cite[Proposition 2.1]{DGV}. \PQED 

\bigskip

(2):  We induct on $t = \depth(M)$.  If $t=0$, then $$\type(M) = \dim_k \socle(M) \leq \length(\gm(M)) \leq \Deg(M).$$  If $t>0$, let $x \in \m$ be a generic hyperplane on $M$.  By Propositions \ref{typelemma} and \ref{M/xM} and by induction, we have $\type(M) = \type(M/xM) \leq \Deg(M/xM) \leq \Deg(M)$. \PQED

\bigskip

(3):  Let ${\bf x} = x_1,\ldots,x_d$ be a generic system of parameters of $R$.  Then $$\Deg(R) \geq \length(R/{\bf x}R) \geq \embcod(R/{\bf x}R) +1 = \embcod(R)+1.$$ \QED

\bigskip

In light of Proposition \ref{inequalities}, one may ask:  What happens if equality holds in (1), (2) or (3)?  The equality $\embcod(R) + 1 = \Deg(R)$ is dealt with in Section \ref{bdeg-reg}, and the equality $\nu(M) = \Deg(M)$ is the subject of Section \ref{Deg-linres}.  The equality $\type(M) = \Deg(M)$ is more elementary and can be dealt with immediately.

\begin{Proposition} If $\type(M) = \Deg(M)$, then $M$ is a linear Cohen-Macaulay module. \end{Proposition}

\proof Recall from Section \ref{linCM-linBuchs} that a linear Cohen-Macaulay module is a module which is Cohen-Macaulay and satisfies $\nu(M) = \deg(M)$.  To show that $M$ is Cohen-Macaulay, we induct on $t = \depth(M)$.  If $t=0$, then $$\type(M) = \dim_k \socle(M) \leq \length(\gm(M)),$$ but $\length(\gm(M)) < \Deg(M)$ unless $\dim M = 0$, in which case $M$ is Cohen-Macaulay.  If $t> 0$, let $x$ be a generic hyperplane on $M$.  We have $$\type(M) = \type(M/xM) \leq \Deg(M/xM) \leq \Deg(M) = \type(M),$$ whence $\type(M/xM) = \Deg(M/xM)$.  By induction, $M/xM$ is Cohen-Macaulay, and hence so is $M$.  

To show that $M$ is linear, let ${\bf x} = x_1,\ldots,x_d$ be a generic system of parameters of $M$. Then $\type(M) = \deg(M)$ implies that $\type(M/{\bf x}M) = \deg(M/{\bf x}M)$.  But this says that $M/{\bf x}M = M/\m M$, from which it follows that $\nu(M) = \deg(M)$. \QED

\section{Bounding Numbers of Generators and Hilbert Functions}\label{gen-props}

Here we reproduce, without proof, some results from \cite{DGV} which serve to motivate the study of cohomological degrees, and in particular the extremal cohomological degree which is studied in detail in Section \ref{bdeg}.  We preserve the assumptions about $R$ and $\Deg(-)$ from Section \ref{Deg-elem}.

\begin{Proposition} \label{gen-a}(\cite[Theorem 3.1]{DGV}) Let $R$ be Cohen-Macaulay of dimension $d$, and let $I$ be an ideal of codimension $g>0$.  If $\depth(R/I) = r$, then
\begin{eqnarray*} 
\nu(I) & \leq & \deg(R) + (g-1)\Deg(R/I) + (d-g-r)(\Deg(R/I)-\deg(R/I))\\
& = & \deg(R) + (g-1)\deg(R/I) + (d-r-1)(\Deg(R/I)-\deg(R/I)). \end{eqnarray*}
\end{Proposition}

\begin{Proposition} \label{gen-b}(\cite[Theorem 4.6]{DGV}) Let $I$ be an $\m$-primary ideal of $(R,\m)$.  Then
$$\nu(I) \leq \Deg(R){s+d-2 \choose d-1}+{s+d-2 \choose d-2},$$
where $d = \dim(R) \geq 1$ and $s$ is the index of nilpotency of $R/I$. \end{Proposition}

Observe that substituting $I = \m^s$ into Proposition \ref{gen-b} yields a bound for the Hilbert function of $R$.

\begin{Proposition} \label{gen-c}(\cite[Corollary 4.7]{DGV})  The reduction number of $\m$ is bounded as follows.
$$ \red(\m) \leq d \cdot \Deg(R)-2d+1,$$
where $d = \dim(R) \geq 1.$ \end{Proposition}

\section{Linear Resolution}\label{Deg-linres}

In this section we preserve the assumptions about $R$ and $\Deg(-)$ of Section \ref{Deg-elem}.  Our goal is to examine the consequences of the equality $\nu(M) = \Deg(M)$.  Namely, we have:

\begin{Theorem} \label{Yosh} If $\nu(M) = \Deg(M)$, then $\grm(M)$ has a linear resolution. \end{Theorem}

This theorem was first proved by Yoshida in \cite{Yoshida} for the homological degree.  Here we will show that the conclusion follows from only the properties of cohomological degrees.  Our proof follows the original proof of Yoshida, and proceeds by observing that modules which satisfy the equality $\nu(M) = \Deg(M)$ have much in common with Buchsbaum modules.  Indeed, Yoshida defines a {\em weakly linear Buchsbaum module}\footnote{This should be parsed ``(weakly) (linear Buchsbaum) module'' and not ``(weakly linear) (Buchsbaum) module'', since such modules are not necessarily Buchsbaum.} as one which satisfies the equality $\nu(M) = \hdeg(M)$.  However, we will use the less unwieldy term {\em Yoshida module}.  Note that Theorem \ref{Yosh} is a generalization of (half of) Proposition \ref{BHU}.

\begin{Proposition} \label{Yoshprop} Let $\Deg(-)$ be a cohomological degree function and suppose that $\nu(M) = \Deg(M)$. \begin{enumerate}
\item $M \cong N \oplus L$, where $\m L = 0$ and $\nu(N) = \Deg(N)$ and $\depth(N) > 0$.
\item If $x$ is a generic hyperplane then $\nu(M/xM) = \Deg(M/xM)$. \end{enumerate} \end{Proposition}

\proof 

(1):  Let $L=\gm(M)$ and $N=M/L$.  Then we have 
\begin{eqnarray*}
\nu(L)+\nu(N) & \leq & \length(L) + \Deg(N) = \Deg(M) \\
& = & \nu(M) \leq \nu(L)+\nu(N). \end{eqnarray*}
Hence $\length(L) = \nu(L)$, which is equivalent to $\m L = 0$, and $\Deg(N)=\nu(N)$.  Thus it only remains to show that the exact sequence
$$0 \lar L \lar M \lar N \lar 0$$
splits.  Tensoring with $R/\m$ yields
$$ \overline{L} \lar \overline{M} \lar \overline{N} \lar 0.$$
But $L = \overline{L}$, and, by counting dimensions, we have that the map $L \lar \overline{M}$ is an inclusion.  The composition of maps $M \lar \overline{M} \lar L$ gives the desired splitting. \PQED

(2):  Squeeze:
$$\nu(M/xM) = \nu(M) = \Deg(M) \geq \Deg(M/xM) \geq \nu(M/xM).$$ \QED


\begin{Lemma}\label{QM} Suppose $\nu(M) = \Deg(M)$, and let $Q = (x_1,\ldots,x_d)$, where $x_1,\ldots,x_d$ is a generic system of parameters for $M$.  Then $\m M = QM$. \end{Lemma}

\proof The result is clear in dimension $0$. Let $\overline{M} = M/QM$. By Proposition \ref{Yoshprop}(2), we have $\nu(\overline{M}) = \Deg(\overline{M})$, and so we are done. \QED

\bigskip

The following lemma, and its proof, is the analogue of \cite[Lemma 1.14]{SV}.

\begin{Lemma}\label{SV-lemma} Let $I$ be an ideal of $R$ such that $\nu(M) = \Deg(M/IM)$ and let $Q=({\bf x})$, where ${\bf x}=x_1,\ldots,x_r$ is a generic partial system of parameters for $M/IM$.  Then
$$ [IM:\m]\cap Q^kM \subseteq IQ^{k-1}M \qquad \mbox{for }k\geq 1.$$ \end{Lemma}

\proof Assume the lemma is false and that $I$ is maximal possible such that, for some $Q=({\bf x})$, where ${\bf x}$ is a partial system of parameters for $M/IM$,
$$ [IM:\m]\cap Q^kM \not\subseteq IQ^{k-1}M .$$
For this $I$, let $k \geq 1$ be minimal with respect to this property.  Let $V = [IM:\m]$.
\bigskip

\noindent {\bf Claim:} $[V:x]= V$ for generic $x\in\m$.

\noindent {\bf Proof of Claim.}  Since $\nu(M/IM) = \Deg(M/IM)$, we have, by Proposition \ref{Yoshprop}(1),
$$ M/IM \cong M' \oplus L,$$
where $\depth (M') > 0$ and $\m L = 0$.  Say $z \in [V:x]$.  So $xz \in V$ and $x\m z \subseteq IM$.  Let $\overline{z}$ be the image of $z$ in $M/IM$.  Then $x\m \overline{z} = 0$, and hence $x\overline{z} \in L$.  Since $x$ is generic, this implies $\overline{z} \in L$, so $\m \overline{z} = 0$, {\em i.e.} $z \in V$. \PQED
\bigskip

Note that the claim implies that $[V:x^t] = V$ for generic $x$ and $t \geq 1$.

If $r = 0$ then $Q=0$ and the lemma is trivial.  So $r > 0$.  Let $a \in (V \cap Q^kM)\setminus IQ^{k-1}M.$  Write
$$a = \sum_{i=1}^r x_i u_i, \qquad \mbox{where }\quad  u_i \in (x_1,\ldots x_i)^{k-1}M.$$
Let 
$$b = \sum_{i=1}^{r-1}x_i u_i.$$
Now, $xb + xx_ru_r = xa \in IM$, for all $x \in \m$, so
$$ b \in [(I+(x_r))M : \m]\cap Q_1^{k-1},$$
where $Q_1 = (x_1,\ldots,x_{r-1})$.  Since $\Deg(M/(IM+x_rM) = \Deg(M/IM)$, by the maximality of $I$ we must have 
$$b \in (I+(x_r))Q_1^{k-1}M.$$
Thus \begin{eqnarray*}
a = b+x_ru_r & \in & (IQ_1^{k-1} + x_rQ_1^{k-1} + x_rQ^{k-1})M \cap V \\
& = & IQ^{k-1}M + (x_rQ^{k-1}M)\cap V \\
& = & IQ_1^{k-1}M + x_r(Q^{k-1}M \cap [V:x_r]) \\
& = & IQ_1^{k-1}M + x_r(Q^{k-1}M\cap V). \end{eqnarray*}

If $k=1$, then $$a \in IM + x_rV = IM$$ and if $k \geq 2 $, then by the minimality of $k$, $$Q^{k-1}M \cap V \subseteq IQ^{k-2}M,$$ and so
$$ a \in IQ_1^{k-1}M + x_rIQ^{k-2}M = IQ^{k-1}M.$$
In either case we have a contradiction, and so we are done. \QED

\begin{Lemma}\label{exact}  Suppose $\depth(M) > 0$ and $\nu(M) = \Deg(M)$.  Let $x$ be a generic hyperplane section on $M$ and let $\overline{M} = M/xM$.  Then we have an exact sequence
$$0 \lar \grm(M) \stackrel{x^*}{\lar} \grm(M) \lar \grm( \overline{M} ) \lar 0,$$
where $x^*$ is the initial form of $x$ in $\grm(R)$.\end{Lemma}

\proof Let $Q = (x_2,\ldots,x_r)$, where $x,x_2,\ldots,x_r$ is a generic system of parameters for $M$.  By Lemma \ref{SV-lemma} with $I=(x)$, we have
$$ xM \cap Q^kM \subseteq [xM : \m]\cap Q^kM \subseteq xQ^{k-1}M.$$
Applying Lemma \ref{QM}, we have $$
xM\cap \m^kM =  xM \cap (x,Q)^kM 
 =  (xM \cap Q^kM) + xQ^{k-1}M,$$
whence $$xM \cap \m^kM  =  x \m^{k-1}M.$$
It follows that $[\m^{k+1}M:x] = \m^kM$, so $x^{*}$ is regular on $\grm(M)$.  Furthermore,
$$\begin{array}{rcccl}
\left[ \grm(\overline{M}) \right]_k & = & \frac{\m^kM}{(xM\cap \m^k)+ \m^{k+1}M} & & \\
& = & \frac{\m^kM}{x\m^{k-1}M + \m^{k+1}M} & = & \left[ \frac{\grm(M)}{x^{*}\grm(M)} \right]_k . 
\end{array}$$
This completes the proof. \QED

\bigskip 

\noindent {\bf Proof of Theorem \ref{Yosh}.}   We proceed by induction on the dimension of $M$.  If $\dim M = 0$, the result is trivial, since then $M$ is a vector space, $\m M = 0$.  Let $L = H_{\m}^0(M)$ and $N=M/L$.  By Lemma \ref{Yoshprop}, $M \cong N\oplus L$, $L$ has a linear resolution and $N$ satisfies the condition $\nu(N) = \Deg(N)$, so we may assume $\depth(M) > 0.$

Let $x$ be a generic, regular hyperplane on $M$.   Then $x^{*}$ is a regular hyperplane on $\grm{(M)}$.  Since $(\grm{M})/(x^{*}\grm{M}) \cong \grm{(M/xM)}$ has a linear resolution by induction, it follows that $\grm (M)$ does as well, by Proposition \ref{reg-inductive}\QED

\bigskip

Unfortunately, the converse to Theorem \ref{Yosh} fails for some cohomological degrees. We illustrate this with the homological degree.  Recall that a Yoshida module is a module which satisfies $\nu(M) = \hdeg(M)$.

\begin{Lemma} If $M$ is a finitely generated module of dimension $d\geq 2$ and $M$ has no associated primes of dimension $\leq 1$, then
$$\hdeg(\m M) = \hdeg(M) + (d-1)\cdot\nu(M).$$ \end{Lemma}

\proof We may assume $(R,\m)$ is the homomorphic image of a Gorenstein ring.  Since $M$ has no associated primes of dimension $\leq 1$, we have that $M_0=0$ and that $M_1$ is a module of finite length, by Proposition \ref{h-Spec}.  Let $N = \m M$.  The exact sequence $$0 \lar N \lar M \lar M/N \lar 0$$ yields 
$$0 \lar M_1 \lar N_1 \lar \left( \frac{M}{\m M} \right)^{\vee} \lar 0.$$
So
$$ \hdeg(M) = \deg(M) + \sum_{i=2}^{d-1}{d-1 \choose i}\hdeg(M_i) + (d-1)\ell(M_1) $$
and 
\begin{eqnarray*}
\hdeg(N) & = & \deg(N) + \sum_{i=2}^{d-1}{d-1 \choose i}\hdeg(N_i) + (d-1)\ell(N_1)\\
& = & \deg(M) + \sum_{i=2}^{d-1}{d-1 \choose i}\hdeg(M_i) + (d-1)(\ell(M_1) + \ell(M/\m M)^{\vee}) \\
& = & \hdeg(M) + (d-1)\nu(M). \end{eqnarray*} \QED

\begin{Proposition} If $M$ is a Yoshida module of dimension $d \geq 2$, then $\m^tM$ is not Yoshida
for some $t$. \end{Proposition}

\proof Clearly, we may assume that $\depth(M) >0$.  Since $M$ is Yoshida, we have by \cite[Proposition 2.2]{Yoshida} that $M$ has no associated primes of dimension $\leq 1$.  Now, by the above lemma we have that $$\hdeg (\m M) =\hdeg(M) + (d-1)\nu(M) = d\cdot\nu(M).$$  Therefore, if $\m^t M$ is Yoshida for all $t$, we have $\nu(\m^t M) = \nu(M)\cdot d^t$, which contradicts a famous theorem \cite{Cookie}. \QED

\section{The Homological Degree and Castelnuovo-Mumford regularity}

Throughout this section, $(R,\m)$ will denote a homogeneous $k$-algebra, where $k$ is infinite.  We choose a presentation, $R=S/I$, where $S=k[x_1,\ldots,x_n],\; n = \embdim(R),$ and $I$ is a homogeneous ideal.  We let $r=\reg(R)$.

If $k$ is infinite and $\Deg(-)$ is a cohomological degree on ${\cal M}(R)$, then by Proposition \ref{regDeg} we have $r < \Deg(R)$.  In this section, we are interested in finding a complement to Proposition \ref{regDeg}: an inequality of the form $\Deg(R) \leq f(r)$, where $f$ is some function of $r$.  In addition to depending on the cohomological degree $\Deg(-),\; f$ will necessarily incorporate the embedding dimension, $n$, since if $R=S/(x_1,\ldots,x_n)^2$ we have $\reg(R)=1$ but $\Deg(R)=\length(R)=n+1$.

Observe that both $d=\dim(R)$ and $e=\deg(R)$ may be bounded in terms of $n$ and $r$:  Clearly $d \leq n$.  By Propositions \ref{Hilbert-init} and \ref{gin-summary}, we have that $\deg(R)=\deg(S/\gin(I))$ and $\reg(R) = \reg(S/\gin(I))$, and for a given $n$ and $r$ there are only finitely many possible generic initial ideals.  Therefore, we may harmlessly incorporate $d$ and $e$ into any bound for $\Deg(R)$. Hence our problem is as follows.

\begin{Problem}\label{Deg-reg-prob} Given a cohomological degree $\Deg(-)$ on ${\cal M}(R)$, find a function $f$ such that $$\Deg(R) \leq f(r,n,d,e),$$
where $r=\reg(R),\; n=\embdim(R),\; d=\dim(R)$ and $e=\deg(R)$. \end{Problem}

With this in mind, we turn our attention to what is, at this point, our only {\em explicitly} defined (but see Chapter 3) cohomological degree, $\hdeg(-)$.

\begin{Conjecture}\label{hdeg-reg} There is a function $f(r,n,d,e)$ such that $$\hdeg(R) \leq f(r,n,d,e).$$ \end{Conjecture}

Conjecture \ref{hdeg-reg} would follow immediately from the following.

\begin{Conjecture} \label{hdeg-gin} $\hdeg(R)\leq \hdeg(S/\gin(I))$. \end{Conjecture}

\begin{Theorem} \label{hdeg-dim2} Conjecture \ref{hdeg-reg} is true for $d \leq 2$. \end{Theorem}

Before proving Theorem \ref{hdeg-dim2}, we need to establish some tools.  The first is a proposition taken from \cite[Theorem 2.13]{V96}.

\begin{Proposition}\label{hdeg-prop} Let $$0 \lar A \lar B \lar C \lar 0$$ be an exact
sequence of finitely generated modules with $\length(C) < \infty$.  Then
$$\hdeg(B) \leq \hdeg(A) + \length(C).$$
\end{Proposition}

The following proposition is based upon the idea of a theorem of Herzog-K\"{u}hl \cite{HK} and of Huneke-Miller \cite{HuM};  see also \cite[Theorem 4.1.15]{BH}.  In Appendix \ref{App}, we will develop some further applications of Proposition \ref{linres}.

\begin{Proposition}\label{linres} Let $M$ be a finitely generated graded module over $S=k[x_1,...,x_n]$ such that $M[j]$ has a linear resolution for some $j$, $\,\dim M=d$ and $\depth (M) \geq d-1$.  Then \begin{eqnarray*} \beta_{n-d+1} & = & \mu - e \\ \beta_{n-d} & = & (n-d+1)\mu - (n-d)e, \end{eqnarray*}
where $\beta_i$ is the $i^{th}$ Betti number of $M$, $\mu = \nu(M)$ is the minimal number of generators of $M$, and $e = deg(M)$ is the degree of $M$.
\end{Proposition}

\proof Without loss of generality, we may assume that $M$ is generated in degree 0.  The Hilbert series of $M$ is given by \begin{displaymath}
H_M(t) = \frac{S_M(t)}{(1-t)^n} = \frac{Q_M(t)}{(1-t)^d},\end{displaymath}
where $Q(1) = e(M)$.  Thus we have the system of equations \begin{eqnarray*}
S_M(1) & = & 0 \\
S_M'(1) & = & 0 \\
S_M''(1) & = & 0 \\
 & \vdots & \\
S_M^{(n-d-1)}(1) & = & 0 \\
S_M^{(n-d)}(1) & = & (-1)^{n-d}(d-n)!e 
\end{eqnarray*}

On the other hand, from Proposition \ref{notsohard}, we have 
\begin{displaymath}
S_M(t) = \sum_{i=0}^{n-d+1}(-1)^i\beta_it^i. \end{displaymath}
Thus the system of equations becomes 
\begin{eqnarray*} 
\left[ \begin{array}{ccccccc}
\frac{1!}{1!} & \frac{2!}{2!} & \frac{3!}{3!} & & \cdots & & \frac{(n-d+1)!}{(n-d+1)!} \\ 
\frac{1!}{0!} & \frac{2!}{1!} & \frac{3!}{2!} & & & & \\
& \frac{2!}{0!} & \frac{3!}{1!} & & & & \vdots \\
& & \frac{3!}{0!}& & & & \\
& & & \ddots & & & \\
& & & & \frac{(n-d-1)!}{0!} & \frac{(n-d)!}{1!} & \frac{(n-d+1)!}{2!} \\
& & & & & \frac{(n-d)!}{0!} & \frac{(n-d+1)!}{1!} 
\end{array} \right]  & \left[ \begin{array}{c}
-\beta_1 \\ \beta_2 \\ -\beta_3 \\ \beta_4 \\ \vdots \\ (-1)^{n-d}\beta_{n-d} \\ (-1)^{(n-d+1)}\beta_{n-d+1} \end{array} \right]
=  & \\
& \left[ \begin{array}{c}
-\mu \\ 0 \\ 0 \\ 0 \\ \vdots \\ 0 \\ (-1)^{n-d}(n-d)!e \end{array} \right] 
& \end{eqnarray*}

\noindent One easily calculates that this system of equations is equivalent to:
\begin{eqnarray*} \left[ \begin{array}{cccccc}
\frac{1!}{1!} & \frac{2!}{2!} & \frac{3!}{3!} &&& \\
& \frac{1!}{0!} & \frac{2!}{1!} &&& \\
&& \frac{2!}{0!} &&& \\
&&& \ddots && \\
&&&& \frac{(n-d-1)!}{0!} & \frac{(n-d)!}{1!} \\
&&&&& \frac{(n-d)!}{0!} \end{array} \right] &
\left[ \begin{array}{c}
-\beta_1 \\ \beta_2 \\ -\beta_3 \\ \beta_4 \\ \vdots \\ (-1)^{(n-d+1)}\beta_{n-d+1} \end{array} \right] = & \\   
& \left[ \begin{array}{c} -\mu \\ \mu \\ -2\mu \\ \vdots \\ (-1)^{n-d-2}(n-d-1)!\mu \\ (-1)^{n-d}(n-d)!(e-\mu) \end{array} \right] & \end{eqnarray*}
The claim follows from the final two rows. \QED

\begin{Proposition}\label{finlen} Let $S=k[x_1,\ldots,x_n]$ and $\phi: S^m \lar S^d$ such that $\phi$ is given by a matrix of linear forms in $S$ and $\length(\coker \phi) < \infty$.  Then \begin{displaymath}
\length(\coker \phi) \leq { n+d-1 \choose d-1} \end{displaymath}
\end{Proposition}

\proof  We may assume that $k$ is infinite.  We let $\m$ denote the irrelevant maximal ideal of $S$.

Let $I_d(\phi)$ be the ideal generated by the maximal minors of $\phi$.  Then $\length(\coker \phi) < \infty$ is equivalent to $I_d(\phi)$ is $\m$-primary.  Thus, by \cite[Theorem 1]{EN} we must have $m \geq n+d-1$.

If $m=n+d-1$, then by \cite{BV} we have $\length(\coker \phi) = \length(S/I_d(\phi))$.  On the other hand, by \cite[Theorem 1]{EN} the Eagon-Northcott complex is a minimal free resolution of $I_d(\phi)$; in particular, the ${ n+d-1 \choose d}$ maximal minors of $\phi$ minimally generate $I_d(\phi)$.  Therefore, $I_d(\phi) = \m^d$ and $\length(\coker \phi) = { n+d-1 \choose d-1}$.

If $m \geq n+d$, it suffices to show that $\phi$ can, after a change of basis of $S^m$, be written in the form \begin{displaymath}
\phi = \left[ \;\;\;\; \psi \;\;\;\; \left| \begin{array}{c} \alpha_1 \\ \vdots \\ \alpha_d \end{array} \right. \right], \end{displaymath}
where $\alpha_1,\ldots,\alpha_d$ are linear forms in $S$ and $\psi$ is a $d \times(m-1)$ matrix of linear forms in $S$ such that $I_d(\psi)$ is again $\m$-primary.  Indeed, if this is the case, then $\coker \phi$ is a homomorphic image of $\coker \psi$ which, by induction on $m$, satisfies $\length(\coker \psi) \leq {n+d-1 \choose d-1}$, which completes the proof.

For any ${\bf w} \in k^n$, write $\phi_{\bf w}$ for the $d \times m$ matrix of scalars obtained by evaluating the entries of $\phi$ at the point ${\bf w}$.  Let $A$ be the ring of polynomials $k[{\bf y}] = k[y_1,\ldots,y_m]$, and $\phi'$ the $d \times n$ matrix of linear forms in $A$ whose $i^{th}$ column is given by $\phi_{{\bf e}_i}\cdot {\bf y}$, where ${\bf e}_1,\ldots,{\bf e}_n$ is the standard basis for $k^n$.  For any generalized column \begin{displaymath}
a = \sum \beta_i (\phi_{{\bf e}_i}\cdot {\bf y}) = (\phi_{\sum \beta_i {\bf e}_i})\cdot {\bf y} \end{displaymath} of $\phi'$, we write $(a)$ for the ideal in $A$ generated by the entries of $a$; in particular, we write $(a)_{\langle 1\rangle}$ for the degree 1 component of this ideal.  We have the following equivalences:
\begin{eqnarray*} I_d(\phi) \mbox{ is $\m$-primary } & \Longleftrightarrow & \mbox{Every $\phi_{\bf w}$ (${\bf w} \neq 0$) is surjective} \\
 & \Longleftrightarrow & \height (a) = d \mbox{ for every generalized column $a$ of $\phi'$} \end{eqnarray*}

Thus, recasting the problem in terms of $\phi'$, we must show that there is a hyperplane $H$ such that for every generalized column $a$ in $\phi'|_H$, we have $\height (a) = d$; in other words, that there is a linear form $h \in A_{\langle 1\rangle}$ such that $h \not \in (a)_{\langle 1\rangle}$ for every generalized column $a$ of $\phi'$.  Each $(a)_{\langle 1\rangle}$ is a vector space of dimension $d$ in $A_{\langle 1\rangle}$ which has dimension $m$.  Since the vector spaces $(a)_{\langle 1\rangle}$ are parameterized by ${\bf P}^{n-1}$, their union is a variety of dimension $n+d-1$ in $A_{\langle 1\rangle}$.  Since $m \geq n+d$, we can find such
an $h$. \QED

\bigskip

\noindent {\bf Proof of Theorem \ref{hdeg-dim2}.}  Let $L=\gm(R)$.  Then $\max\{\, i \; | \; L_{\langle i \rangle} \neq 0\,\} \leq r$.  Thus, since $L$ is a subquotient of $S$, we may bound $\length(L)$ in terms of $n$ and $r$.  Now, if we let $\overline{R} = R/L$, then we have $\reg(\overline{R}) \leq r,\,\deg(\overline{R}) = e$ and $\dim(\overline{R}) = d$.  Since $\hdeg(R)=\hdeg(\overline{R})+\length(L)$, we may replace $R$ by $\overline{R}$ and assume that $\depth(R)>0$.

If $d \leq 1$, we are done, so we assume $d=2$.  Now we show that the claim is true if $R$ is unmixed; that is, if every associated prime of $R$ has dimension $2$.

Let $M=R_{\geq r}$ and let $K=R/M$.  We have $\hdeg(R) \leq \hdeg(M) + \length(K)$ by Proposition \ref{hdeg-prop}, and, since $K$ is a subquotient of $S$ such that $K_{\langle i \rangle}=0$ for $i\geq r$, we can bound $\length(K)$ in terms of $n$ and $r$.  Thus is suffices to bound $\hdeg(M)$.

Now, $\Ass(M) \subseteq \Ass(R)$, so $M$ is also unmixed.  Hence $\length(M_1) < \infty$ and 
$$
\hdeg(M) = \deg(M) + \hdeg(M_1) = e+\length(M_1).$$
It remains to bound $\length(M_1)$.

$M$ has a linear resolution,
$$0 \lar S^{\mu-e}\stackrel{\phi}{\lar}S^{(n-d+1)\mu-(n-d)e} \lar \cdots \lar S^{\mu} \lar M \lar 0,$$
where $\mu=\nu(M)$ and $\phi$ is a matrix of linear forms.  Hence $M_1$ has a presentation
$$S^{(n-d+1)\mu-(n-d)e}\stackrel{\phi^t}{\lar}S^{\mu-e}\lar M_1 \lar 0.$$
Since $\length(M_1) < \infty$, by Proposition \ref{finlen} we have
$$\length(M_1) \leq {n+\mu-e-1 \choose \mu-e-1}.$$
But $\mu = \nu(M) = H(R;r)$ can also be bounded in terms of $n$ and $r$, so we are done in the case where $R$ is unmixed.

Now suppose $R$ is not unmixed.  Say $R=S/I$ and $I=I_1\cap I_2$, where $I_1$ and $I_2$ are the intersections of the primary components of $I$ of dimensions $1$ and $2$, respectively.  We have an exact sequence
\begin{equation}\label{I1I2}
0 \lar I_2/I \lar R \lar S/I_2 \lar 0. \end{equation}
Since $I_2/I \cong (I_1+I_2)/I_1$, it has dimension $1$, and, since it is embedded in $R$, it must be Cohen-Macaulay.  From (\ref{I1I2}), we obtain the exact sequence
$$0 \lar (S/I_2)_1 \lar R_1 \lar (I_2/I)_1 \lar 0,$$
and since $\length(S/I_2)_1<\infty$ and $I_2/I$ is Cohen-Macaulay, we have
\begin{eqnarray*}
\hdeg (R) & = & \deg(R) + \hdeg(R_1) \\
& = & \deg(S/I_2) + \length((S/I_2)_1)+\hdeg((I_2/I)_1) \\
& = & \hdeg(S/I_2) + \deg(I_2/I). \end{eqnarray*}

Since every module in (\ref{I1I2}) has positive depth, we may choose a linear form $h$, which is regular on all three modules.  Thus we have an inclusion
$$0 \lar I_2/(hI_2+I) \lar R/hR.$$
So $\deg(I_2/I) = \length(I_2/(hI_2+I))$ can be bounded in terms of $n$ and $\reg(R/hR) = r$.  Now it only remains to bound $\hdeg(S/I_2)$.  It is unmixed, so it is bounded in terms of $n$ and $\reg(S/I_2)$, but $\reg(S/I_2) \leq r$.  Indeed, by Proposition \ref{reg-ses}, we have
$$\reg(S/I_2) \leq \max\{\,r,\,\reg(I_2/I)-1\,\}.$$
But $I_2/I$ is a Cohen-Macaulay submodule of $R$, so $\reg(I_2/I) \leq \reg(R)$, and we are done. \QED

\bigskip

By following the above proof carefully, one may extract an explicit bound for $\hdeg(R)$.  The present author has in fact done so, but to record it here would violate standards of mathematical aesthetics.  (The sensitive reader may feel that such standards have already been violated.)  In light of the proof of Theorem \ref{hdeg-dim2}, it seems that Conjecture \ref{hdeg-reg} may be rather difficult, unless one can prove Conjecture \ref{hdeg-gin}.

\bigskip

By now, we have compiled a list of grievances against what is at this point our only friend, the only cohomological degree which we know how to compute, $\hdeg(-)$.  The converse to Theorem \ref{Yosh} does not hold for $\hdeg(-)$, and it is difficult to relate $\hdeg(-)$ to Castelnuovo-Mumford regularity.  The homological degree tends to give numbers which are too large, due to the deadly combination of the recursive nature of the definition with the binomial coefficients contained therein.  This is reflected in the fact that the inequality $$\hdeg(M/xM) \leq \hdeg(M)$$ for generic $x$ is frequently strict, even when $x$ is regular on $M$.

It turns out that the extremal cohomological degree, $\bdeg(-)$, about which we have so far said little, is free of these difficulties.  This is the subject of the next chapter.

%% file: thesis-chap3.tex
\chapter{The Extremal Cohomological Degree, $\bdeg$}\label{bdeg}

\section{Axiomatic Characterization of $\bdeg$} \label{bdeg-axiom}

In this section, $(R,\m)$ will be a local Noetherian ring with infinite residue field $k$, or a homogeneous algebra over an infinite field, $k$, with irrelevant maximal ideal $\m$.  We recall the definition in the previous chapter of the extremal cohomological degree.

\begin{Definition} Let $U(-)$ be a notion of genericity for ${\cal M}(R)$.  Then for all $M$ in ${\cal M}(R)$,
$$\bdeg_U(M)=\min \{\,\Deg(M)\; | \; \Deg(-) \mbox{ is a cohomological degree relative to } U(-)\,\}.$$
\end{Definition}

The main theorem of this section is the following.

\begin{Theorem}\label{b-prop} There exists a notion of genericity $U(-)$ such that for all $M$ in ${\cal M}(R)$ with $\depth(M)>0$,
$$\bdeg_U(M) = \bdeg_U(M/xM)$$
for $x \in U(M)$. \end{Theorem}

An easy induction argument shows that if $U(-)$ and $V(-)$ are two notions of genericity for which the conclusion to Theorem \ref{b-prop} holds, then we have $\bdeg_U(-) = \bdeg_V(-)$.  Therefore, once we have proved Theorem \ref{b-prop}, we will have established the existence of a cohomological degree -- written simply $\bdeg(-)$, without any reference to a notion of genericity -- which satisfies the following axiomatic description.

\begin{Definition} \label{bdeg-def} $\bdeg(-)$ is the unique numerical function $$\bdeg(-) : \{\mbox{isomorphism classes of }{\cal M}(R)\,\} \lar {\bf N}$$ satisfying the following conditions.\begin{enumerate} 
\item $\bdeg(M) = \deg(M)$ if $M$ is Cohen-Macaulay.
\item If $L = \gm(M)$, then $\bdeg(M) = \bdeg(M/L) + \length(L)$.
\item If $\depth(M)>0$ and $x$ is a generic hyperplane on $M$, then $\bdeg(M) = \bdeg(M/xM)$. \end{enumerate} \end{Definition}

We remark that the above definition is redundant.  In place of condition {\em 1.}, one only need require
\begin{center} {\em 1'.\quad  $\bdeg(0) = 0.$} \end{center}
Then, if $\dim M = 0$, from the second condition we obtain
$$\bdeg(M) = \bdeg(0) + \length(M) = \length(M).$$
Applying the third condition, we recover
$$\bdeg(M) = \deg(M)$$
if $M$ is Cohen-Macaulay.

In order to prove Theorem \ref{b-prop}, we first establish a series of lemmas.  The first shows that, by passing to the completion if necessary, we may assume that $R$ is the homomorphic image of a Gorenstein ring.

\begin{Lemma} Let $(R,\m)$ denote a Noetherian local ring and let $(\hat{R},\m\hat{R})$ denote the $\m$-adic completion of $R$.  If $M$ is an $R$-module, write $\hat{M}$ for $M\otimes_R\hat{R}$.  Suppose $U(-)$ is a notion of genericity on ${\cal M}(\hat{R})$ with respect to which there exists at least one cohomological degree.  Let $\overline{U(M)}$ denote the image of $U(M)$ in $\m\hat{R}/(\m\hat{R})^2 \cong \m/\m^2$.  For $M$ in ${\cal M}(R)$, we define $U(M)$ by $$U(M) = \{\, x \;|\; \overline{x} \in \overline{U(\hat{M})}\,\}.$$ Then $U(-)$ is a notion of genericity on ${\cal M}(R)$, and, for all $M$ in ${\cal M}(R)$, $$\bdeg_U(M) = \bdeg_U(\hat{M}).$$ \end{Lemma}

\proof That $U(-)$ is a notion of genericity on ${\cal M}(R)$ is clear.  To prove that $\bdeg_U(M) = \bdeg_U(\hat{M})$, we observe that $(M/xM)\,\hat{} \cong \hat{M}/x\hat{M},\,\gm(M) \cong \gm(\hat{M}),$ and $(M/\gm(M))\,\hat{} \cong \hat{M}/\gm(\hat{M}),$ and apply induction on the dimenstion of $M$ to the characterization of $\bdeg_U(-)$ given in Equation (\ref{bdeg-characterization}). \QED

\bigskip

For the remainder of this section, we assume that $R$ is the homomorphic image of a Gorenstein ring.  The next lemma is adapted from the appendix of Watanabe \cite{Watanabe}.

\begin{Lemma}\label{watanabe} Let $M$ be an $R$-module of finite length.  Then the function 
$$x \mapsto \length(M/xM)$$ achieves its minimum generically among hyperplanes. \end{Lemma}

\begin{Corollary}\label{watanabe-dim1} Let $M$ be a one dimensional finitely generated $R$-module.  Then the function $$x \mapsto \length(M/xM)$$ achieves its minimum generically among hyperplanes $x$ such that $\length([0:x]_M) < \infty$. \end{Corollary}

\proof Let $L= \gm(M)$ and let $x$ be a regular hyperplane on $M/L$.  Then the short exact sequence
$$ 0 \lar L \lar M \lar M/L \lar 0$$ induces
$$ 0 \lar L/xL \lar M/xM \lar M/(xM + L) \lar 0.$$
Hence $\length(M/xM) = \length(M/(xM + L)) + \length(L/xL)$. We have $\length(M/(xM+L)) \geq \deg(M)$, with equality generically (when $x$ is a degree system for $M$).  Since   $\length(L/xL)$ achieves its minimum generically, $\length(M/xM)$ must also achieve its minimum generically. \QED

%
%

%
%

\begin{Lemma} \label{silly} Suppose $\depth (M) > 0$.  Then for $x,y \in R$, regular on $M$,  we have
$$ \length\left( \frac{\gm(M/xM)}{y \gm(M/xM)} \right) = \length \left( \frac{\gm(M/yM)}{x \gm(M/yM)}\right).$$ \end{Lemma}

\proof $R$ is the homomorphic image of a Gorenstein ring, $S$.  Let $M_1 = \Ext_S^{n-1}(M,S)$, where $n=\dim(S)$.  We will show that $$\length\left( \frac{\gm(M/xM)}{y \gm(M/xM)} \right) = \length \left(\frac{M_1}{(x,y)M_1}\right),$$ which will suffice, since the roles of $x$ and $y$ may be inverted.  First, note that
\begin{eqnarray*}
\length\left( \frac{\gm(M/xM)}{y \gm(M/xM)} \right) & = & \length( 0 :_{\gm(M/xM)}y) \\
& = & \length \left( \frac{\Ext_S^n(M/xM,S)}{y\Ext_S^n(M/xM,S)} \right). \end{eqnarray*}
Since $\Ext_S^n(M/xM,S) \cong M_1/xM_1$, we are done. \QED

\begin{Lemma} \label{depth0} Let $U(-)$ be a notion of genericity, $x \in U(M)$ and let $L= \gm(M)$.  Then
$$\bdeg_U(M) \geq  \bdeg_U(M/xM) -\length(L/xL)+\length(L).$$
Equality holds if and only if 
$$ \bdeg_U(\overline{M}) = \bdeg_U(\overline{M}/x\overline{M}),$$
where $\overline{M} = M/L$. \end{Lemma}

\proof 
\begin{eqnarray*}
\bdeg_U(M/xM) & = & \length(L/xL) + \bdeg_U(\overline{M}/x\overline{M})\\
& \leq & \length(L/xL) + \bdeg_U(\overline{M}) \\
& = & \length(L/xL) + \bdeg_U(M) - \length(L). \end{eqnarray*} 
\QED

\begin{Lemma}\label{degree-system} Given an exact sequence of modules $$0 \lar A \lar B \lar C \lar 0,$$ such that $\length(C) < \infty$ and a sequence ${\bf x}$ of elements in $\m$, then ${\bf x}$ is a degree system for $A$ if and only if it is a degree system for $B$. \end{Lemma}

\proof This is trivial if $\dim(A)=\dim(B)=0$, so assume that $A$ and $B$ have positive dimension.  Thus if ${\bf x}$ is a multiplicity system on $A$ or on $B$ then it is a multiplicity system on both $A$ and $B$, and we have $e({\bf x};A) = e({\bf x};B)$.  Since the property of being a degree system is a generic condition, there exists a sequence which is a degree system for both $A$ and $B$.  Thus $\deg(A)=\deg(B)$, so ${\bf x}$ is a degree system for $A$ if and only if it is a degree system for $B$. \QED

\newcommand{\E}{{\cal E}}
\begin{Definition} If $\dim(M)=0$, we let $\E(M) = \emptyset$, and if $\dim(M)=d>0$, we define $\E(M)$ inductively as follows. $$\E(M) = \left(\bigcup_{i=1}^{d-1}\E(M_i)\right)\cup \{M\}.$$
\end{Definition}

\begin{Definition} If $M$ is a module of dimension at most $1$, then we define the {\em Rees number} of $M$ to be the integer $$\rees(M) = \min\{\, \length(M/xM) \;|\; x \in \m\,\}.$$ \end{Definition}

By Lemma \ref{watanabe} and Corollary \ref{watanabe-dim1}, the Rees number is attained by a generic element.

\begin{Definition} We say that a sequence ${\bf x} = x_1,\ldots,x_t$ in $\m$ is a {\em $t$-Rees sequence} for $M$ if $$\length([0:x_{i+1}]_{M/(x_1,\ldots,x_i)M}) < \infty$$ for $i=0,\ldots,t-1$ and, for all such sequences ${\bf y}$, we have $$\rees((M/{\bf x}M)_1) \leq \rees((M/{\bf y}M)_1).$$ \end{Definition}

\begin{Lemma}\label{colonlemma} If $x$ and $y$ are regular on $M$, then $[0:x]_{M/yM} \cong [0:y]_{M/xM}$. \end{Lemma}

\proof We have an exact sequence $$0 \lar yM/xyM \lar M/xyM \lar M/yM \lar 0.$$
Tensoring with $R/xR$, we obtain 
$$0 \to yM/xyM \stackrel{\sim}{\to} yM/xyM \to [0:x]_{M/yM} \to yM/xyM \to M/xM \to M/(x,y)M \to 0.$$  But there is an isomorphism $M/xM \stackrel{\sim}{\lar}yM/xyM$, and the kernel of the composition $$M/xM \lar yM/xyM \lar M/xM$$ is $[0:y]_{M/xM}$. \QED 

\begin{Lemma} If ${\bf x}$ is a $t$-Rees sequence for $M$, then any permutation of ${\bf x}$ is a $t$-Rees sequence for $M$. \end{Lemma}

\proof This follows from Lemma \ref{colonlemma}.

\begin{Lemma} \label{preswitch} Suppose $\dim(M)\geq 1$.  Let $x,y \in \m$ be such that $[0:x]_M,\,[0:y]_M,$ and $[0:y]_{M/xM}$ each have finite length.  Then
$$\length\left(\frac{(M/xM)_1}{y(M/xM)_1}\right) = \length\left(\frac{(M/yM)_1}{x(M/yM)_1}\right).$$ \end{Lemma}

\proof We may assume that $M$ has positive depth, in which case $x$ and $y$ are regular on $M$.  The map $$(M/xM)\stackrel{y}{\lar}(M/xM)$$ may be decomposed as 
\begin{equation}\label{ONE}M/xM \stackrel{\sim}{\to} yM/xyM \onto (x,y)M/xM \hookrightarrow M/xM.\end{equation}  From (\ref{ONE}), we derive exact sequences
\begin{equation}
\label{TWO} 0 \lar [0:y]_{M/xM} \lar yM/xyM \lar (x,y)M/xM \lar 0, \mbox{ and }\end{equation} \begin{equation}
\label{THREE} 0 \lar (x,y)M/xM \lar M/xM \lar M/(x,y)M \lar 0. \end{equation}
The sequence (\ref{TWO}) yields
\begin{equation} 
\label{FOUR} ((x,y)M/xM)_1 \stackrel{\sim}{\lar}(yM/xyM)_1, \mbox{ and }
\end{equation} \begin{equation}
\label{FIVE} 0 \lar ((x,y)M/xM)_0 \lar (yM/xyM)_0 \lar ([0:y]_{M/xM})^{\vee} \lar 0. \end{equation}
From (\ref{THREE}) we obtain
\begin{equation} \label{SIX}
(M/xM)_1 \to ((x,y)M/xM)_1 \to (M/(x,y)M)_0 \to (M/xM)_0 \to ((x,y)M/xM)_0 \to 0. \end{equation}
Splicing (\ref{FOUR}) and (\ref{FIVE}) into (\ref{SIX}), and bearing in mind (\ref{ONE}), we have
$$(M/xM)_1 \stackrel{y}{\lar} (M/xM)_1 \lar (M/(x,y)M)_0 \lar$$ 
$$ \lar (M/xM)_0 \stackrel{y}{\lar} (M/xM)_0 \lar ([0:y]_{M/xM})^{\vee}\lar 0.$$
Thus $$\length\left(\frac{(M/xM)_1}{y(M/xM)_1}\right) = \length((M/(x,y)M)_0) - \length([0:y]_{M/xM}).$$
Since $[0:y]_{M/xM} \cong [0:x]_{M/yM}$ by Lemma \ref{colonlemma}, we may reverse the roles of $x$ and $y$, whence $$\length\left(\frac{(M/xM)_1}{y(M/xM)_1}\right) = \length\left(\frac{(M/yM)_1}{x(M/yM)_1}\right).$$ \QED

\begin{Lemma} Given a finitely generated module $M$, any $t$ sufficiently generic linear combinations of the generators of $\m$ form a $t$-Rees sequence for $M$. \end{Lemma}

\proof Let ${\bf x} = x_1,\ldots,x_t$ be a $t$-Rees sequence for $M$.  For a generic $z_1$, we have that $[0:z_1]_M$ and $[0:z_1]_{M/(x_2,\ldots,x_t)M}$ have finite length and $$\rees((M/{\bf x}M)_1) = \length\left( \frac{(M/{\bf x}M)_1}{z_1(M/{\bf x}M)_1}\right).$$
Thus, by Lemma \ref{preswitch}, we have that $z_1,x_2,\ldots,x_t$ is a $t$-Rees sequence.  Continuing in the fashion, we may replace $x_2$ by a generic $z_2$, and so on. \QED

\begin{Definition} An element $x$ of $\m$ is called a Rees element for $M$ if $x$ is part of a $t$-Rees sequence for all $t\leq \dim(M)-1$. \end{Definition}

\begin{Definition} Let $M$ be a finitely generated module.  Let
$$U(M) = \left\{ x \;\left|\; \begin{array}{c}\mbox{for each } N\in\E(M),\, \length([0:x]_N)<\infty \\ \mbox{and } x \mbox{ is part of a degree system for $N$}\end{array}\right.\right\},$$
and $$V(M) = U(M)\cap\{\,x\;|\; x \mbox{ is a Rees element for $M$ and } \length(M_1/xM_1) = \rees(M_1)\;\}.$$
\end{Definition}

Observe that $U(-)$ and $V(-)$ are notions of genericity.

\begin{Lemma} \label{switcheroo} If $x \in V(M)$ and $y \in V(M)\cap\left( \bigcap_{N\in\E(M)}V(N/xN)\right)$, then $x \in V(M/yM)$. \end{Lemma}

\proof We begin by showing that if $x\in U(M)$ and $y \in U(M)\cap\left( \bigcap_{N\in\E(M)}U(N/xN)\right)$ then $x \in U(M/yM)$.  We proceed by induction on $\dim(M)$.  Observe that we may assume that $\depth(M)>0$.

Let $N \in \E(M/yM)$.  We wish to show that $[0:x]_N$ has finite length and that $x$ is part of a degree system on $N$.
\begin{itemize}
\item Case 1: $N=M/yM$.  Since $y \in U(M/xM)$, we have that $[0:y]_{M/xM}$ has finite length and that $y$ is part of a degree system for
 $M/xM$.  It follows that $[0:x]_{M/yM}$ has finite length and that $x$ is part of a degree system for $M/yM$.
\item Case 2: $N=(M/yM)_i$, where $i<\dim(M/yM) = \dim(M)-1$.  The exact sequence $$0 \lar M \stackrel{y}{\lar}M \lar M/yM \lar 0$$ yields the exact sequence
\begin{equation} \label{armyguy} 0 \lar M_{i+1}/yM_{i+1} \lar (M/yM)_i \lar [0:y]_{M_i} \lar 0.\end{equation}
Since $\dim(M_{i+1}) < \dim(M)$, we have, by induction, that $x\in U(M_{i+1}/yM_{i+1})$.  Since $[0:y]_{M_i}$ and $[0:x]_{M_{i+1}/yM_{i+1}}$ have finite length, $[0:x]_{(M/yM)_i}$ has finite length as well.  Furthermore, since $x$ is part of a degree system for $M_{i+1}/yM_{i+1}$, it follows that $x$ is part of a degree system for $(M/yM)_i$, by Lemma \ref{degree-system}.
\item Case 3: $N \in \E((M/yM)_i), N \neq (M/yM)_i$, where $i<\dim(M/yM) = \dim(M)-1$.  From (\ref{armyguy}), we have
$$(M_{i+1}/yM_{i+1})_j \cong ((M/yM)_i)_j$$ if $j \geq 2$.  Therefore, either $N  = ((M/yM)_j)_1$ or 
$N$ is isomorphic to a module in $\E(M_{i+1}/yM_{i+1})$.  Since $x \in U(M_{i+1})$ and $$y\in U(M_{i+1})\cap\left(\bigcap_{N\in\E(M_{i+1})}U(N/xN)\right),$$ we have $x \in U(M_{i+1}/yM_{i+1})$, so we are done in the former case.

If $N = ((M/yM)_i)_1$, then we have, from (\ref{armyguy}),
$$0 \lar N \lar (M_{i+1}/yM_{i+1})_1 \lar ([0:y]_{M_i})^\vee.$$
Since $x \in U((M_{i+1}/yM_{i+1})_1)$, it follows that $x \in U(N)$ by Lemma \ref{degree-system}.\end{itemize}

Now suppose $x,y \in V(M)$ and $y\in V(N/xN)$ for each $N \in \E(M)$.  By the above reasoning, we have that $x \in U(M/yM)$.  It remains to show that $x$ is a Rees element for $M/yM$ and that $$\length\left(\frac{(M/yM)_1}{x(M/yM)_1}\right) = \rees((M/yM)_1).$$
Since $x \in V(M)$ and $y\in V(M/xM)$, we have that $x$ is a Rees element for $M$ and that $y$ is a Rees element for $M/xM$.  Thus $x,y$ can be extended to a $t$-Rees sequence for $M$, for every $t\leq \dim(M)-1$.  Since any permutation of $t$-Rees sequence is again a $t$-Rees sequence, $x$ is a Rees element for $M/yM$.

Finally, since $x$ and $y$ are Rees elements for $M$, we have $$\rees((M/xM)_1)=\rees((M/yM)_1),$$ and, by Lemma \ref{preswitch}, we have $$\length\left(\frac{(M/yM)_1}{x(M/yM)_1}\right)=\length\left(\frac{(M/xM)_1}{y(M/xM)_1}\right).$$  Since $y \in V(M/xM)$, we have $$\length\left(\frac{(M/xM)_1}{y(M/xM)_1}\right) = \rees((M/xM)_1),$$ and this forces $$\length\left(\frac{(M/yM)_1}{x(M/yM)_1}\right) = \rees((M/yM)_1).$$ \QED

\bigskip

\noindent {\bf Proof of Theorem \ref{b-prop}.} Observe that $U(-)$ is the notion of genericity with respect to which $\hdeg(-)$ is a cohomological degree (Proposition \ref{hdeg}).  Since $V(-)$ is a refinement of $U(-)$, $\bdeg(-)$ is well defined. We will show that if $\depth(M)>0$ then $$\bdeg_V(M) = \bdeg_V(M/xM)$$ for all $x \in V(M)$.  We proceed by induction on $\dim(M)$.  Let $x_1$ and $x_2$ lie in $V(M)$.  It suffices to show that $\bdeg_V(M/x_1M) = \bdeg_V(M/x_2M)$.  Let $$y \in V(M)\cap\left( \bigcap_{N\in\E(M)}V(N/x_1N)\cap V(M/x_2M)\right).$$  By Lemma \ref{switcheroo}, we have $x_1,x_2 \in V(M/yM)$.  Thus, by induction and Lemma \ref{depth0},
$$
\bdeg_V(M/x_1M) =  \bdeg_V(M/(x_1,y)M) - \length\left(\frac{\gm(M/x_1M)}{y\gm(M/x_1M)}\right) + \length(\gm(M/x_1M)).$$
 Now, $$(\gm(M/x_1M))^{\vee} \cong (M/x_1M)_0 \cong M_1/x_1M_1,$$
so $\length(\gm(M/x_1M)) = \rees(M_1)$ and, similarly, $\length(\gm(M/yM)) = \rees(M_1)$.  Thus, by Lemma \ref{silly}, \begin{eqnarray*}
\bdeg_V(M/x_1M) & = & \bdeg_V(M/(x_1,y)M) - \length\left(\frac{\gm(M/yM)}{x_1\gm(M/yM)}\right) + \length(\gm(M/yM))\\
& = & \bdeg_V(M/yM). \end{eqnarray*}
Similarly, $\bdeg_V(M/x_2M) = \bdeg_V(M/yM)$, so we are done. \QED

\bigskip

We remark that the above proof is not constructive, in the sense that it is not clear how one could determine, using a program such as {\sc Macaulay}, if a given hyperplane section is generic for a given module.  Nonetheless, a ``random'' hyperplane will, with probability one, be generic.  Therefore, with the help of a random number generator, one has an algorithm for the computation of $\bdeg(-)$.  Though the theoretical status of this algorithm may be dubious, from a practical standpoint it is quite satisfactory:  the script that the present author has written for {\sc Macaulay} works, and works quickly (much more quickly than the script for $\hdeg(-)$, for example).

\section{Basic Properties}

We preserve the assumptions on $(R,\m)$ of the previous section.

\noindent Lemma \ref{depth0} can be recast in terms of $\bdeg(-)$.  We record the reformulation for reference.

\begin{Proposition} Let $M$ be in ${\cal M}(R)$ and let $L=\gm(M)$.  Then for a generic hyperplane $x$,
$$\bdeg(M) = \bdeg(M/xM) - \length(L/xL) + \length(L).$$
\end{Proposition}

\begin{Proposition} \label{submod} Given an exact sequence $$0 \lar A \lar B \lar C \lar 0,$$ where $\length(C) < \infty$, we have
$$\bdeg(B) \leq \bdeg(A) + \length(\gm(B)) - \length(\gm(A)).$$  In particular, we have $\bdeg(B) \leq \bdeg(A) + \length(C)$, and, if $\depth (B) > 0$, $\bdeg(B) \leq \bdeg(A)$. \end{Proposition}

\proof We induct on the dimension of $B$.  If $\dim B = 0$, the statement is trivial, so assume $\dim B > 0$.  First, we treat the case in which $\depth (B) > 0$.  Choose $x$ to be generic on $A$ and $B$.  Since $x$ is regular on $B$, we obtain the exact sequence

\setlength{\unitlength}{.02cm}
\begin{picture}(650,130)\thicklines
\put(100,90){$0$}
\put(135,95){\vector(1,0){30}}
\put(185,90){$[0:x]_C$}
\put(245,95){\vector(1,0){30}}
\put(285,90){$A/xA$}
\put(345,95){\vector(1,0){30}}
\put(385,90){$B/xB$}
\put(445,95){\vector(1,0){30}}
\put(485,90){$C/xC$}
\put(545,95){\vector(1,0){30}}
\put(595,90){$0$}
\put(325,80){\vector(1,-1){21.2}}
\put(367,58.8){\vector(1,1){21.2}}
\put(346.5,43){$K$}
\put(325,16.8){\vector(1,1){21.2}}
\put(367,38){\vector(1,-1){21.2}}
\put(310,2){$0$}
\put(392,2){$0$}

\end{picture}

\noindent where $K$ is the kernel of the map $B/xB \lar C/xC$.  Now, $x$ is generic on $B$, so $\bdeg(B) = \bdeg(B/xB)$, and, by induction, we have
\begin{eqnarray*} 
\bdeg(B/xB) & \leq & \bdeg(K) + \length(\gm(B/xB)) - \length(\gm(K)) \\
& = & \bdeg(A/xA) - \length(0:_Cx) + \length(\gm(B/xB)) - \length(\gm(K)) \\
& = & \bdeg(A) - \length(C/xC) + \length(\gm(B/xB)) - \length(\gm(K)) \\
& \leq & \bdeg(A). \end{eqnarray*}
The last inequality is from the exact sequence
$$ 0 \lar \gm(K) \lar \gm(B/xB) \lar C/xC$$
which is induced by
$$0 \lar K \lar B/xB \lar C/xC \lar 0.$$  This completes the proof in the case where $\depth (B) > 0$.  If $\depth (B) = 0$, let $\overline{B} = B/\gm(B)$ and $\overline{A}=A/\gm(A)$.  Since $\overline{A} \hookrightarrow \overline{B}$ we have $\bdeg(\overline{A}) \leq \bdeg(\overline{B})$, from which the proposition follows. 
\QED

\begin{Proposition} \label{qBuchs} Let $M$ in ${\cal M}(R)$ have dimension $d$.  Then
$$ \bdeg(M) \leq \deg(M) + \sum_{i=0}^{d-1} {d-1 \choose i} \length( H_{\m}^i(M) ).$$
Equality holds if $M$ is Buchsbaum\footnote{Originally, this read, "...if 
$\m H_{\m}^i(M) = 0$ for $i=0,\ldots,d-1$; that is, if $M$ is 'quasi-Buchsbaum.'"  The proof 
incorrectly assumed that the property of being quasi-Buchsbaum was preserved under generic
hyperplane sections.  This error was pointed out in {\rm {\sc U. Nagel and T. R\"{o}mer}, 
Extended degree functions and monomial modules, {\em Trans. Amer. Math. Soc.}  
{\bf 358}  (2006), 3571-3589.}}.
\end{Proposition}

\proof Clearly, we may assume $\depth (M) > 0$.  Let $x$ be generic on $M$.  We have the exact sequence
$$0 \lar M \stackrel{x}{\lar} M \lar M/xM \lar 0.$$
The long exact sequence in cohomology breaks up into short exact sequences
\begin{equation}\label{ses} 0 \lar \frac{H_{\m}^i(M)}{xH_{\m}^i(M)} \lar {H_{\m}^i(M/xM)} \lar ([0:x]_{H_{\m}^{i+1}(M)}) \lar 0.\end{equation}
Thus $\length(H_{\m}^i(M/xM)) \leq \length({H_{\m}^i(M)}) + \length({H_{\m}^{i+1}(M)}),$ so the claim follows from Pascal's triangle and induction on the dimension of $M$.

If $M$ is Buchsbaum, then so is $M/xM$. Since $x {H_{\m}^i(M)} = 0$ for $0 \leq i \leq d-1$, (\ref{ses}) becomes  
$$0 \lar {H_{\m}^i(M)} \lar {H_{\m}^i(M/xM)} \lar {H_{\m}^{i+1}(M)} \lar 0,$$ and we finish as before. \QED

\begin{Exercise} Let $A$ be a finitely generated graded $R$-algebra,
$$A = R\oplus A_{\langle 1 \rangle} \oplus A_{\langle 2 \rangle} \oplus \cdots,$$ such that $A_{\langle 0 \rangle} = R$ and $A$ is generated in degree $1$ as an $R$-algebra.  Let $M$ be a finitely generated graded $A$-module,
$$M = \bigoplus_{i\in{\bf Z}}M_{\langle i \rangle}.$$
Show that the function $$n \mapsto \bdeg(M_{\langle n \rangle})$$ agrees with a polynomial for $n \gg 0$. \end{Exercise}

\hint This may be proven as in \cite[Proposition 3.3.4]{Doering}.

\section{Linear Resolution}

In this section, we preserve the assumptions and notation of Section \ref{bdeg-axiom}.  Recall that in Section \ref{Deg-linres} we showed that if $\Deg(-)$ is any cohomological degree on ${\cal M}(R)$, then
$$ \nu(M) = \Deg(M) \; \; \Longrightarrow \; \; \grm(M) \mbox{ has linear resolution.}$$
Furthermore, we saw, using the example of $\Deg(-) = \hdeg(-)$, that the converse does not hold.  In this section, we prove that the converse does hold for $\bdeg(-)$; that is, we have:

\begin{Theorem} \label{linresthm} The following are equivalent.
\begin{itemize} 
\item $\grm(M)$ has linear resolution.
\item $\nu(M) = \bdeg(M).$
\item $\nu(M) = \bdeg(\grm(M)).$
\end{itemize}
\end{Theorem}

We begin by establishing some tools for dealing with associated graded modules.

\begin{Lemma} Suppose $x$ is superficial on $M$.  Then we have an exact sequence
$$ 0 \lar K \lar \frac{\grm(M)}{x^*\grm(M)} \lar \grm\left( \frac{M}{xM} \right) \lar 0$$ where $\length(K) < \infty$. \end{Lemma}
\proof $$ \left[ \frac{\grm(M)}{x^*\grm(M)} \right]_i \cong \frac{\m^iM}{x\m^{i-1}M+\m^{i+1}M},$$ and 
$$\left[ \grm\left(\frac{M}{xM}\right)\right]_i \cong \frac{\m^iM}{xM\cap\m^iM + \m^{i+1}M}.$$
Since $x\m^{i-1}M \subseteq xM\cap\m^iM$, there is a natural surjection
$$ \frac{\grm(M)}{x^*\grm(M)} \onto \grm\left(\frac{M}{xM}\right).$$
The kernel is 
$$K = \bigoplus_{i\geq 0} \frac{xM\cap\m^iM + \m^{i+1}M}{x\m^{i-1}M+\m^{i+1}M},$$ and since by Lemma \ref{superlemma}, $xM\cap \m^iM = x\m^{i-1}M$ for $i \gg 0$, $\,\length(K) < \infty$. \QED

\begin{Corollary} If $x$ is superficial on $M$, then
$$\Deg\left(\grm\left( \frac{M}{xM} \right)\right) \leq \Deg\left( \frac{\grm(M)}{x^*\grm(M)}\right)$$ for any cohomological degree $\Deg(-)$. \end{Corollary}

\proof This follows from Proposition \ref{basics}(2). \QED

\begin{Theorem} \label{grim} $\bdeg(M) \leq \bdeg(\grm(M))$. \end{Theorem}

\proof We argue by induction on $\dim M$.  Suppose $\depth (M) > 0$.  Choose $x$ generic on $M$ such that $x^{*}$ is generic on $\grm(M)$.  In particular, we choose $x$ to be superficial on $M$.  Then
\begin{eqnarray*} \bdeg(M) & = & \bdeg(M/xM) \\
& \leq & \bdeg(\grm(M/xM))\\
& \leq & \bdeg(\grm(M)/x^*\grm(M))\\
& \leq & \bdeg(\grm(M)). \end{eqnarray*}

\noindent So assume $\depth (M) = 0$.  Let $L = \gm(M)$ and let
$$L' = \bigoplus_{i\geq 0} \frac{L\cap \m^iM}{L\cap \m^{i+1}M}.$$ Note that $\length(L) = \length(L')$ and that $\grm(M/L) \cong \grm(M)/L'$.  So
\begin{eqnarray*}
\bdeg(M) & = & \bdeg(M/L) + \length(L)\\
& \leq  & \bdeg(\grm(M/L)) + \length(L)\\
& = & \bdeg(\grm(M)/L')+\length(L') = \bdeg(\grm(M)). \end{eqnarray*}
\QED

\bigskip

\noindent {\bf Proof of Theorem \ref{linresthm}}.  
We have
\begin{eqnarray*}
\nu(M) = \bdeg(\grm(M)) & \Longrightarrow & \nu(M) = \bdeg(M) \\
& \Longrightarrow & \grm(M) \mbox{ has linear resolution } \end{eqnarray*}
by Theorems \ref{Yosh} and \ref{grim}.  Hence, replacing $M$ by $\grm(M)$ and $R$ by $\grm(R)$, we may assume $M \cong \grm(M)$ and need only show
\begin{eqnarray*}
M \mbox{ has linear resolution } & \Longrightarrow & \nu(M) = \bdeg(M). \end{eqnarray*}

We argue by induction on $\dim M$.  The case $\dim M = 0$ is trivial, since such a module has linear resolution if and only if it is a vector space concentrated in degree $0$.  So suppose $\dim M > 0$.  Since $M$ has linear resolution, $M \cong L \oplus N$, where $\m L = 0$ and $N$ is a module with linear resolution and positive depth.  Let $x$ be generic on $N$.  Then $N/xN$ has linear resolution also, so, by induction
\begin{eqnarray*}
\bdeg(M) = \length(L) + \bdeg(N) & = & \nu(L) + \bdeg(N/xN)\\
& = & \nu(L) + \nu(N/xN)\\
& = & \nu(M). \end{eqnarray*}\QED

\begin{Corollary} Let $M$ be a finitely generated graded module over $S=k[x_1,\ldots,x_n]$.  If $M$ has positive depth and $\reg (M) = r$, then
$$\bdeg(M) \leq H(M;r).$$ \end{Corollary}

\proof Let $N = M_{\geq r}$.  By Lemma \ref{submod}, $\bdeg(M) \leq \bdeg(N)$.  But $N$ has a linear resolution, so $\bdeg(N) = \nu(N) = H(M;r)$. \QED

\bigskip

This corollary gives an indication of the connection between $\bdeg(-)$ and Castelnuovo-Mumford regularity.  We shall have more to say about this in the next section.

\section{Castelnuovo-Mumford Regularity}\label{bdeg-reg}

In this section, we assume that $R$ is a homogeneous $k$-algebra with presentation $R=S/I$, where $S=k[x_1,\ldots,x_n]$ and $I$ is a homogeneous ideal.  We further assume that the characteristic of the field $k$ is zero.  Under these assumptions, we will solve Problem \ref{Deg-reg-prob} in its entirety for $\Deg(-)= \bdeg(-)$.  

Recall the $r^{th}$ Macaulay representation of an integer $e$:
$$e = {k(r) \choose r} + \cdots + {k(t) \choose t},$$
where $k(r) > \cdots > k(t) \geq t$.  Then we let
$$e^{(r,d)} = {k(r) + d \choose r} + \cdots + {k(t)+ d \choose t}.$$
Our main theorem is the following.

\begin{Theorem}\label{homog-bound} If $\embdim(R)=n,\, \deg(R) = e,\, \reg(R) = r,\, \dim R = d$, and $\depth (R) = g$, then $$\bdeg(R) \leq e + {n-g+r \choose r} - e^{(r,d-g)}. $$ \end{Theorem}

The method of proof is, briefly, as follows.  First, we show that $\bdeg(R)$ is unchanged when we pass to the generic initial ideal.  Generic initial ideals are Borel-fixed, so our next step is to give a combinatorial description of $\bdeg(R)$ when $I$ is Borel-fixed.  Finally, we exhibit the Borel-fixed ideal which maximizes $\bdeg(R)$ and calculate.  We remark that Theorem \ref{homog-bound} is sharp, since we will give an ideal for which the bound is obtained.

\begin{Theorem} \label{bdeg-gin} $\bdeg(R) = \bdeg(\gin(R))$. \end{Theorem} 

\proof We induct on the dimension of $R=S/I$, the 0-dimensional case being trivial.  First, suppose $\depth(S/I) > 0$ and let $x$ be a generic hyperplane.  Recall from Section \ref{GIN-intro} that $S_x$ denotes $S/(x)$ and $I_x$ denotes $(I+(x))/(x) \subseteq S_x$.  Then 
\begin{eqnarray*}
\bdeg(S/I) & = & \bdeg(S/I+(x)) = \bdeg(S_x/\gin(I_x))\\
&  = & \bdeg(S/\gin(I)+(x_n)) = \bdeg(S/\gin(I)),\end{eqnarray*}
by induction and 
Propositions \ref{BF-summary}(2) and \ref{gin-gen-hyp}. 

Now, if $\depth(S/I) = 0$,
 we have 
\begin{eqnarray*}
\bdeg(S/I) & = & \bdeg(S/I^{sat}) + \length(I^{sat}/I) \\
& = & \bdeg(S/\gin(I^{sat})) + \length(\gin(I^{sat})/\gin(I)) \\
& = & \bdeg(S/\gin(I)^{sat}) + \length(\gin(I)^{sat}/\gin(I)) \\
& = & \bdeg(S/\gin(I)). \end{eqnarray*} \QED 

\begin{Definition} Let $I$ be a Borel-fixed monomial ideal.  For a monomial $m$, 
let $$\range(m) = \max \{ t : x_t \mbox{ divides } m \}.$$  A monomial $m$ is 
{\em standard} with  respect to $I$ if $m \not\in I$ but 
$m \in [I : (x_{1},\ldots,x_t)^{\infty}] = [I : x_t^{\infty}] $, where $t = \range(m)$; that is, if
$$ \left(\prod_{i=1}^{t}x_i^{\alpha_i}\right) \cdot m \in I,$$
for $\alpha_{1}\ldots,\alpha_t \gg 0$. \end{Definition}

\begin{Proposition} \label{standard} If $I$ is Borel-fixed, $\bdeg (S/I)$ is equal to the number of monomials standard with respect to $I$. \end{Proposition}

\proof We proceed by induction on the dimension of $S/I$.  If $\dim S/I = 0$, then $\bdeg(S/I) = \length(S/I)$ is equal to the number of monomials not in $I$.  It is easy to see that, in this case, every such monomial is standard.  So assume $\dim S/I > 0$.  Then we have 
$$ \bdeg(S/I) = \bdeg(S/I^{sat}) + \length(I^{sat}/I),$$
with $I^{sat} = (I:x_n^{\infty})$.  So
$$\length(I^{sat}/I) =  |\{ \mbox{ monomials } m \not\in I, \mbox{ but } m \in (I:x_n^{\infty}) \, \}|, $$
and, by induction,
\begin{eqnarray*}
\bdeg(S/I^{sat}) & = & \bdeg(S/I^{sat}+(x_n)) \\
& = & \bdeg \left( \frac{k[x_1,\ldots,x_{n-1}]}{(I:x_n^{\infty})\cap k[x_1,\ldots,x_{n-1}]} \right) \\
& = & |\{ m \in (I:x_t^{\infty}), \mbox{ where } t=\range(m)<n, \mbox{ but } m \not\in (I:x_n^{\infty}) \}|. \end{eqnarray*}
Clearly, these two sets are disjoint, and their union is the set of all monomials standard with respect to $I$. \QED

\bigskip

Note that if $m$ is standard and $t = \max \{\, i \; | \; m \in (I:x_i^{\infty})\,\}$, then $(m,\{x_{t+1},\ldots,x_n\})$ agrees with what Sturmfels, Trung and Vogel \cite{STV} call a {\em standard pair}.  It follows that if $I$ is Borel-fixed, $\adeg(S/I) = \bdeg(S/I)$, where $\adeg(-)$ is the arithmetic degree function. See \cite[Lemma 3.3]{STV}.

In what follows, we write $\Mon\{x_{i1},\ldots,x_{ij}\}$ for the set of all monomials in the variables  \linebreak $x_{i1},\ldots,x_{ij}$ and $\Mon_k\{x_{i1},\ldots,x_{ij}\}$ for the set of all monomials in $\Mon\{x_{i1},\ldots,x_{ij}\}$ of degree $k$.  Let $\succ_{lex}$ denote the {\em lexicographic order} on monomials in $S$:  $m_1 \succ_{lex} m_2$ if $\deg_{x_i}m_1 = \deg_{x_i}m_2$ for $i<k$ but $\deg_{x_k}m_1 > \deg_{x_k}m_2$ for some $1 \leq k \leq n$.  Note that we do {\em not} sort by total degree here, as we did for the reverse lexicographic order. A set of the form $M(m) = \{\, m' \in \Mon_k\{x_1,\ldots,x_c\} \; | \; m' \prec_{lex} m \, \}$, where $m \in \Mon_k\{x_1\ldots,x_c\}$, is called a {\em lex-segment} in $\Mon_k\{x_1,\ldots,x_c\}$.
We say that $M \subseteq \Mon_k\{x_1,\ldots,x_c\}$ satisfies the {\em shift axiom} if
$$ m \in M, \, x_i \mbox{ divides } m,\mbox{ and } i<j\leq c \Longrightarrow \left( \frac{x_j}{x_i}\right)m \in M.$$

\begin{Lemma}\label{monk} Out of all subsets $M$ of $\Mon_k\{x_1,\ldots,x_c\}$ of size $\varepsilon$ which satisfy the shift axiom, the one which is a lex-segment minimizes the function
$$\range(M) = \sum_{m\in M} \range(m).$$ \end{Lemma}

\proof Let $t = \min\{\, \range(m) \; | \; m\in M \,\}$.  Then all monomials $m$ of degree $k$ with support in $\{x_{t+1},\ldots,x_c\}$ lie in $M$.  Thus, to minimize $\range(M)$, it suffices to minimize $t$.  A lex-segment accomplishes this. \QED

\begin{Definition} Say $M\subseteq \Mon\{x_1,\ldots,x_c\}$ is $(r,e)$-compact if \begin{enumerate}
\item $|M|=e$,
\item $\deg(m) \leq r$ for all $m \in M$,
\item If $m \in M$ and $m'$ divides $m$, then $m' \in M$, and
\item (The shift axiom) If $m\in M, \; x_i$ divides $m$, and $i<j\leq c$, then $(x_j/x_i)m \in M$. \end{enumerate} \end{Definition}

Let $M_{(r,e)} \subseteq \Mon\{x_1,\ldots,x_c\}$ consist of the $e$ lexicographically least monomials of degree $\leq r$.  Observe that $M_{(r,e)}$ is $(r,e)$-compact.

\begin{Lemma} \label{omega} Let $\omega : {\bf N} \to {\bf Z}$ be such that if $a<b$, then $\omega(a)<\omega(b)$.  For all $M \subseteq \Mon\{x_1,\ldots,x_c\}$, let $$\omega(M) = \sum_{m\in M}\omega(\deg(m)).$$  If $M \subseteq \Mon\{x_1,\ldots,x_c\}$ is $(r,e)$-compact, then $\omega(M) \leq \omega(M_{(r,e)})$. \end{Lemma}

\proof For $M \subseteq \Mon_k\{x_1,\ldots,x_c\}$, write $$\sigma(M) = \{\, m \in \Mon_{k+1}\{x_1,\ldots,x_c\} \; | \; m/x_i \in M \mbox{ if } x_i \mbox{ divides } m \, \},$$ and for $m \in \Mon_k\{x_1,\ldots,x_c\}$, let 
\begin{eqnarray*}
\sigma(m) & = & \{\, x_i m \; | \; i \geq \range(m)\,\} \\
& = & \{\, m' \in \Mon_{k+1}\{x_1\ldots,x_c\} \; | \; m' \in \sigma(M) \Leftrightarrow m\in M \mbox{ for all }\\ 
&  & \quad M \subseteq \Mon_k\{x_1,\ldots,x_c\} \mbox{ satisfying the shift axiom }\} \end{eqnarray*}
Note that if $M$ satisfies the shift axiom, $\sigma(M)$ can be expressed as the disjoint union $\sigma(M) = \coprod_{m\in M} \sigma(m)$.

Now, it follows from Lemma \ref{monk} that if $M,M' \subseteq \Mon_k\{x_1,\ldots, x_c\}$, both satisfying the shift axiom and both of size $\varepsilon$, and $M$ is a lex-segment, then $|\sigma(M)| \geq |\sigma(M')|$, since for all $M$ satisfying the shift axiom we have
$$|\sigma(M)| = \sum_{m\in M}|\sigma(m)| = \sum_{m\in M}(c-\range(m)+1) = (c+1)|M|-\range(M).$$

Now, say $M \subseteq \Mon\{x_1,\ldots,x_c\}$ is $(r,e)$-compact.  Let $M_k = M\cap \Mon_k\{x_1,\ldots,x_c\}$, and let $M_k' \subseteq \Mon_k\{x_1,\ldots,x_c\}$ be the lex-segment of size $|M_k|$.

\claim $M' := \cup_0^r M_k'$ is $(r,e)$-compact.

\poc It is clear that $|M'|=e$ and that $M'$ satisfies the shift axiom.  To show that if $m\in M$ and $x_i$ divides $m$ then $m/x_i \in M$, it suffices to show that $|M_{k+1}'| \leq |\sigma(M_k')|$.  But $|M_{k+1}'| = |M_{k+1}| \leq |\sigma(M_k)| \leq |\sigma(M_k')|$. \PQED

Since $\omega(M) = \omega(M')$, we may assume that $M_k$ is a lex-segment for each $k$.  Let $m_1$ be the lexicographically least monomial $m$ such that $m \not\in M$ and $\deg(m) \leq r$, and let $m_2$ be the lexicographically greatest monomial $m$ such that $m \in M$.  Then we have $\deg(m_2) < \deg (m_1)$ because each $M_k$ is a lex-segment.  If $M$ maximizes $\omega(M)$, then $m_2 \prec_{lex} m_1$, since otherwise $M\cup\{m_1\}\backslash\{m_2\}$ is still compact, but with $\omega(M\cup\{m_1\}\backslash\{m_2\}) > \omega(M)$.  Hence $M = \{\, m \in \Mon\{x_1,\ldots,x_c\} \; | \; m \preceq_{lex} m_2,\, \deg (m) \leq r\,\}$. \QED

\bigskip

Suppose $I$ is Borel-fixed.  Then $\Ass(R) \subseteq \{ P_1,\ldots,P_n \}$, where $P_i = (x_1,\ldots,x_i)$.  Thus, if $I$ has height $c$, then it can be written $$I = Q_c\cap \cdots \cap Q_n,$$ where $Q_i$ is $P_i$-primary.

\begin{Lemma} \label{sweetheart} Suppose $I$ is Borel-fixed of height $c$ and regularity $r$.  Let $I' = Q_c\cap\m^{r+1}$.  Then
$$ \bdeg(R) \leq \bdeg(S/I').$$ \end{Lemma}

\proof Let $J = I \cap \m^{r+1}$.  Then $J$ is generated in degree $r+1$ and $\length(I/J) < \infty$.  Thus $\bdeg(S/I) \leq \bdeg(S/J)$ by Proposition \ref{basics} and so we may assume that $I$ is generated in degree $r+1$.

\claim A monomial $m$ is standard with respect to $I$ if and only if $\deg (m) \leq r$ and $m \in [I:x_t^{\infty}]$, where $t = \range(m)$.  

\poc Observe that we need only show that $m$ is standard implies that $\deg (m) \leq r$, since $I$ is generated in degree $r+1$.  Suppose $m = x_1^{\alpha_1}\cdots x_t^{\alpha_t}, \alpha_t>0,$ is standard.  Then there exists $m' = x_1^{\beta_1}\cdots x_t^{\beta_t} \in I$, where $\sum_{i=1}^t \beta_i = r+1$ and $\beta_i\leq \alpha_i$ for $1\leq i \leq t-1$ but $\beta_t > \alpha_t$.  Since $I$ is Borel-fixed, we have
$$m'' = x_1^{\alpha_1}\cdots x_{t-1}^{\alpha_{t-1}} \cdot x_t^{\beta_t - \sum_1^{t-1}(\alpha_i - \beta_i)} \in I.$$
Since $m$ is standard, it does not lie in $I$, and so $\beta_t - \sum_1^{t-1}(\alpha_i - \beta_i)>\alpha_t$, whence $\deg (m) < r+1$. \PQED

Now suppose that $m$ is standard with respect to $I$.  Then $\deg (m) \leq r$ and $$m \in [I:x_t^{\infty}] \subseteq [Q_c\cap\m^{r+1}:x_t^{\infty}]$$ and so $m$ is standard with respect to $Q_c\cap \m^{r+1}$.  By Proposition \ref{standard}, we are done. \QED

\bigskip
Note that the ideal $Q_c\cap\m^{r+1}$ has the same height, degree and regularity as $I$ and is also Borel-fixed.

\bigskip

\noindent {\bf Proof of Theorem \ref{homog-bound}.} First, by factoring out a maximal regular sequence of hyperplanes, we may assume that $g=0$.  By passing to the generic initial ideal, we may assume that $I$ is a Borel-fixed monomial ideal.  Furthermore, by Lemma \ref{sweetheart}, we may assume that $I$ is of the form $I = J \cap \m^{r+1}$, where $J$ is $P_c$-primary, where $c = \height I = n-d$.

Let us establish some notation.  As above, let $c=n-d$.  For any $m \in \Mon\{x_1,\ldots,x_n\}$, write $m = m_c\cdot m_{n-c}$, where $m \in \Mon\{x_1,\ldots,x_c\}$ and $m_{n-c} \in \Mon\{x_{c+1},\ldots,x_n\}$.  Let $M_J = \{ \, m \in \Mon\{x_1,\ldots,x_c\} \; | \; m \not\in J \, \}$.  Observe that $M_J$ is $(r,e)$-compact.

Now, $\bdeg(S/I) = \bdeg(S/I^{sat}) + \length(I^{sat}/I)$, but $I^{sat}=J$, so
\begin{eqnarray*}
\bdeg(S/I) & = & |\{\, m\in\Mon\{x_1,\ldots,x_n\} \; | \; \deg (m) \leq r \mbox{ and either } m \in M_J \mbox{ or } m\in J \, \}| \\ 
& = & {n+r \choose r} - \left| \left\{ \, m \in \Mon\{x_1,\ldots,x_n\} \;
\left| \; \begin{array}{c} \deg(m) \leq r,\, m_c \in M \\ \mbox{ and } m_{n-c}\neq 1\end{array}\,\right.\right\}\right| \\
& = & {n+r \choose r} - \sum_{m\in M_J}{n-c+r-\deg (m) \choose r-\deg (m)}. \end{eqnarray*}
By Lemma \ref{omega}, this is maximized when $M_J=M_{(r,e)}$.  

Let $\widetilde{M}_{(r,e)} = \{ m \in \Mon\{x_1,\ldots,x_n\} \; | \; \deg(m)\leq r \mbox{ and }m_c \in M_{(r,e)} \, \}$.  Then $$\bdeg(S/I) \leq {n+r \choose r} + e - |\widetilde{M}_{(r,e)}|.$$  Now, if $$e = {s_0 + r \choose r} + \cdots + {s_t + r-t \choose r-t},$$ where $s_0 \geq s_1 \geq \cdots \geq s_t \geq 0$, then $M_{(r,e)}$ consists of
\begin{quote}
The ${s_0+r \choose r}$ monomials $m \in \Mon\{x_{c-s_0+1},\ldots,x_c\}$ of degree $\leq r$; and

The ${s_1+r-1 \choose r-1}$ monomials of the form $x_{c-s_0} \cdot m$, where 
$\deg (m) \leq r-1$ and $m \in\Mon\{x_{c-s_1+1},\ldots,x_c\}$; and

$\cdots$

The ${s_t+r-t \choose r-t}$ monomials of the form $x_{c-s_0} \cdots x_{c-s_{t-1}} \cdot m$, where $\deg (m) \leq r-t$ and 
$m \in  \Mon\{x_{c-s_t+1},\ldots,x_c\}$.
\end{quote}
Therefore, $\widetilde{M}_{(r,e)}$ consists of
\begin{quote}
The ${s_0+n-c+r \choose r}$ monomials $m \in \Mon\{x_{c-s_0+1},\ldots,x_n\}$ of degree $\leq r$; and

The ${s_1+n-c+r-1 \choose r-1}$ monomials of the form $x_{c-s_0} \cdot m$, where $\deg (m) \leq r-1$ and $m \in\Mon\{x_{c-s_1+1},\ldots,x_n\}$; and

$\cdots$

The ${s_t+n-c+r-t \choose r-t}$ monomials of the form $x_{c-s_0} \cdots x_{c-s_{t-1}} \cdot m$, where $\deg (m) \leq r-t$ and $m \in  \Mon\{x_{c-s_t+1},\ldots,x_n\}$.
\end{quote}
This completes the proof. \QED

\bigskip

Of course, combining Theorem \ref{homog-bound} with the propositions in Section \ref{gen-props} yields bounds on numbers of generators of ideals, etc., in which the cohomological degrees have been banished, leaving the multiplicity, Castelnuovo-Mumford regularity, and more elementary invariants.

As further application of Theorem \ref{homog-bound}, we will make some remarks on minimal multiplicity, as promised in Section \ref{Deg-elem}.  Let $R$ be a homogeneous $k$-algebra with minimal presentation $R=S/I$, so that $\height I = \embcod(R)$.  $R$ is said to have $r$-linear resolution if $I$ is concentrated in degree $r$ and $I$ has a linear $S$-resolution (after a shift of degree $r$, of course).  Said  another way, $I$ is generated in degree $r$ and $\reg(R)=r-1$.  The following theorem (which is valid even if $k$ has nonzero characteristic) is from \cite[Theorems 2.1 and 3.1]{EG}.

\begin{Proposition}\label{EGprop} Let $R$ be a homogeneous ring of dimension $d$, multiplicity $e$ and embedding codimension $c$.  Then
$$ c+1 \leq e + \sum_{i=0}^{d-1}{d-1 \choose i}\dim \left[ H_{\m}^i(R) \right]_{\langle 1-i \rangle}.$$ Furthermore, if equality holds, then $R$ has $2$-linear resolution.
\end{Proposition} 


\begin{Corollary} Suppose $R$ is a homogeneous ring with embedding codimension $c$.  Then
$$R \mbox{ has $2$-linear resolution } \Longrightarrow \bdeg(R) = c+1.$$  The converse 
is true if $R$ is Buchsbaum \footnote{Originally, if $R$ is quasi-Buchsbaum.  See
the note to Proposition \ref{qBuchs}.}.
\end{Corollary}

\proof By Proposition \ref{inequalities}, we have $c+1 \leq \bdeg(R)$.  If $R$ has 
$2$-linear resolution, then $\reg(R)=1$ and by Theorem \ref{homog-bound}, 
$\bdeg(R) \leq c+1$; hence equality holds.  If, on the other hand, $\bdeg(R) = c+1$ and
$R$ is Buchsbaum, we have
$$c+1 \leq e + \sum_{i=0}^{d-1}{d-1 \choose i}\dim \left[ H_{\m}^i(R) \right]_{\langle 1-i \rangle} = \bdeg(R) =c+1,$$ by Proposition \ref{qBuchs}.  By Proposition \ref{EGprop}, $R$ has a $2$-linear resolution. \QED

%% file: thesis-chap4.tex
\chapter{The Homological Dilworth Number, $\hdil$} \label{hdil-chap}

\section{Motivation}

Once again, we let $(R,\m)$ be a local Noetherian ring with residue field $k$, or a homogeneous $k$ algebra with irrelevant maximal ideal $\m$.  Given an exact sequence of $R$-modules of finite length $$0 \lar A \lar B \lar C \lar 0,$$ the length function satisfies the following equality: \begin{equation}\label{star} \length(B) = \length(A)+\length(C). \end{equation} Usefully, the length function is also an estimator of the number of generators of a module of finite length, \begin{equation}\label{2star} \nu(M) \leq \length(M),\end{equation} though this estimate can be quite crude.  Equation (\ref{star}) and inequality (\ref{2star}) are at the heart of the proofs of the statements (Propositions \ref{inequalities}, \ref{gen-a} and \ref{gen-b}) which bound the number of generators of a module in terms of $\hdeg(-)$ (and other cohomological degrees).  More precisely:  the homological degree
\begin{equation}\label{dagger} \hdeg (M) = \left\{ \begin{array}{ll}
                      \deg M + \displaystyle\sum_{i=0}^{d-1}{d-1 \choose i}\hdeg(M_i) & d = \dim M >0 \\
                      \ell(M) & \dim M =0 \end{array} \right. \end{equation}
satisfies \begin{equation} \label{3star} \hdeg(M/xM) \leq \hdeg(M) \end{equation} for a general hyperplane $x$, a fact whose proof relies crucially on (\ref{star}); inequality  (\ref{3star}) allows us to reduce to dimension $0$, where $\hdeg(-)$ and $\length(-)$ coincide; and finally, inequality (\ref{2star}) provides the desired bound.

Unfortunately, the homological degree inherits the inefficiency of the inequality (\ref{2star}), and therefore bounds using $\hdeg(-)$ can be significantly wasteful.  In order to sharpen $\hdeg(-)$, one may be tempted to replace $\length(-)$ with $\nu(-)$ in (\ref{dagger}).  But $\nu(-)$ does not interact well with exact sequences; instead of (\ref{star}), we have $$\nu(C) \leq \nu(B) \leq \nu(A)+\nu(C),$$ which is not strong enough to pass through the proof of (\ref{3star}).  Our solution is to use another function, the Dilworth number, $\dil(-)$ (Definition
\ref{dil}) which is intermediate to $\nu(-)$ and $\length(-)$, and satisfies (Lemma \ref{diltriangle}) $$\max\{\, \dil(A),\dil(C)\,\} \leq \dil(B) \leq \dil(A)+\dil(C).$$ Though this is still weaker than (\ref{star}), it turns out that it is {\em just enough} to prove an inequality similar to (\ref{3star}).

\section{Dilworth Numbers, Elementary Properties}

Let $(R,\m)$ be a local Noetherian ring with residue field $k$, or a homogeneous $k$-algebra with irrelevant maximal ideal $\m$.

\begin{Definition}\label{dil} Let $M$ be a finitely generated module over $R$.  Then we define the Dilworth number of $M$, $\dil(M)$ as follows:
$$\dil(M) = \sup\{\nu(N)\;|\; N \subseteq M\}.$$\end{Definition}

The Dilworth number was first studied systematically by Watanabe \cite{Wat1}, though it had already been considered in passing by Sally \cite{Sally}.  The terminology is due to Watanabe.  Note that $\dil(M) < \infty$ if and only if $\dim M \leq 1$.  See, for example, \cite[Chapter 3, Theorem 1.2]{Sally}.  Clearly, we have $\nu(M)\leq\dil(M)$.  If equality holds, we say that $M$ is {\em Dilworth}.

Recall (Definition \ref{type-def}) that if $\depth(M)=t$ then the type of $M,\, \type(M),$ is the dimension of $\Ext^t(k,M)$ as a vector space over $k$.  Therefore, if $M$ has finite length, $\type(M) = \dim_k(\socle(M))$.

\begin{Lemma}\label{onto}If $\length(M) < \infty$ then we have:
$$\dil(M) = \sup\{\type(U) \;|\; M \onto U \}.$$\end{Lemma}

\proof 
If $N \subseteq M$ then $$\nu(N) \leq \type(M/\m N),$$ since $N/\m N \hookrightarrow \socle(M/ \m N)$.  Similarly, if $M \onto U$ then $$\type (U) \leq \nu(N),$$ where $$N=[(\ker M \to U):\m]_M,$$ since $N/\m N \onto \socle(U)$. \QED

\bigskip

Recall that if $M$ is a module of finite length, then we write $M^{\vee} = \Hom_R(M,E_R(k))$ for the Matlis dual of $M$; see Section \ref{intro-local-duality}.  As a corollary of Lemma \ref{onto}, we obtain the following result of Ikeda \cite{Ik2}.

\begin{Corollary}\label{mcheck} If $\length(M) < \infty$, then $\dil(M) = \dil(M^{\vee})$.
\end{Corollary} 
\proof This is Matlis duality; see Proposition \ref{matlis-duality}. \QED

\bigskip

On the other hand, from \cite{Ik2} we have the following proposition.

\begin{Proposition}\label{dil-dim1}{\em (Ikeda)} If $\dim M = 1$, we have 
$$\dil (M) = \deg(M) + \dil( \gm(M)).$$ \end{Proposition}

Now, $\gm(M)^{\vee} \cong M_0$, by Local Duality (Proposition \ref{local-duality}).  Thus, combining Corollary \ref{mcheck} and Proposition \ref{dil-dim1}, we have that if $\dim M =1$, then \begin{equation}\label{dil-dim1-eq}\dil(M)=\deg(M)+\dil(M_0),\end{equation} which is a suggestive formula!  In Section \ref{hdil-sec}, we will pursue this analogy with the homological degree, but in the meantime, we need to develop some more tools.

\begin{Lemma}\label{diltriangle} If $\dim M \leq 1$ and $$ 0 \lar N \lar M \lar U \lar 0$$ is exact, then
$$\max \{ \dil(N), \dil(U) \} \leq \dil(M) \leq \dil(N)+\dil(U).$$\end{Lemma}

\proof Trivial. \QED

\bigskip

We remark that $\dil(-)$ is the unique minimal function satisfying the conclusion to Lemma \ref{diltriangle} which also satisfies the inequality $\nu(M)\leq \dil(M)$.

\begin{Lemma}\label{lattice} Let $M$ be a finitely generated module over $R$ such that $\dim M \leq 1$.  Let $d = \dil (M)$.  If $N_1$ and $N_2$ are submodules of $M$ such that $\nu(N_1) = \nu(N_2) = d$, then $\nu(N_1+N_2) = \nu(N_1 \cap N_2) = d$ as well. \end{Lemma}
\proof 
We have the exact sequence,
$$0 \lar N_1 \cap N_2 \lar N_1 \oplus N_2 \lar N_1+N_2 \lar 0.$$
Tensoring with $R/\m$ yields
$$ \overline{N_1 \cap N_2} \lar \overline{N_1} \oplus \overline{N_2} \lar \overline{N_1+N_2} \lar 0.$$
The vector space in the middle has
dimension $2d$.  Since each of the other two vector spaces has dimension less than or equal to $d$, we must have $$\dim \overline{N_1 \cap N_2} = \dim \overline{N_1 + N_2} = d.$$ \QED

\bigskip Thus, the collection of submodules $N$ of $M$ satisfying $\nu(N) = \dil(M)$ forms a lattice under $\cap$ and $+$.  As a consequence, we have the following.

\begin{Corollary}\label{unique} Let $M$ be a finitely generated module over $R$ such that $\dim M \leq 1$.  Let $d = \dil (M)$.  There exists a unique maximal submodule $N \subseteq M$ such that $\nu(M) = d$.  If $\length(M) < \infty$, then there is also a unique minimal submodule with this property.\end{Corollary}

\begin{Definition} For $\dim(M)\leq 1$, let $D(M)$ be the unique maximal submodule of $M$ satisfying $\nu(D(M)) = \dil(M)$.  \end{Definition}

Though we will not require this fact in the sequel, it turns out that $D(-)$ is a functor, as we shall see in the following sequence of propositions.

\begin{Lemma}\label{functorlemma} If $M$ is Dilworth, then all homomorphic images of $M$ are Dilworth. \end{Lemma}

\proof Let $N\subseteq M$.  We wish to show that $M/N$ is Dilworth.  From the
exact sequence
$$0 \lar M/(N\cap \m M) \lar M/N \oplus M/\m M \lar M/(N+\m M) \lar 0,$$
we have
$$\dil(M/N)+\dil(M/\m M) \leq \dil(M/(N\cap\m M)) + \dil(M/(N+\m M)).$$
But $\dil(M/\m M) = \dil(M/(N\cap\m M)) = \nu(M)$, so $\dil(M/N) \leq \dil(M/(N+\m M)) = \nu(M/N)$.  Thus $M/N$ is Dilworth. \QED

\begin{Proposition} Suppose $M$ and $N$ are modules of dimension at most $1$, and let $\phi : M \to N$.  Then $$\phi(D(M)) \subseteq D(N).$$ \end{Proposition}

\proof By replacing $M$ by $D(M)$, we may assume that $M$ is Dilworth.  By Lemma \ref{functorlemma}, $\phi(M)$ is also Dilworth, so we may replace $M$ by $\phi(M)$ and assume $M\subseteq N$.  Thus we need only show that if $M\subseteq N$ and $M$ is Dilworth then $M \subseteq D(N)$.

We have an exact sequence
$$0 \lar M\cap D(N) \lar M\oplus D(N) \lar M + D(N) \lar 0.$$
$M$ is Dilworth, and so $\nu(M\cap D(N)) \leq \nu(M)$.  But $\nu(M\oplus D(N)) = \nu(M)+\dil(N)$, so $\nu(M+D(N)) \geq \dil(N)$.  Hence $M+D(N) \subseteq D(N)$, which finishes the proof. \QED

\begin{Corollary} $D(-)$ is a covariant functor on ${\cal M}_1(R)$, where ${\cal M}_1(R)$ denotes the category of finitely generated $R$-modules of dimension at most $1$. \end{Corollary}

Before leaving this topic, we conclude with some remarks whose proofs are left as an exercise for  the reader.

\begin{Exercise} $D(-)$ is neither right exact nor left exact on ${\cal M}_1(R)$, though it obviously preserves inclusions.  On the other hand, on exact sequences of the form $$0 \to \gm(M) \to M \stackrel{\pi}{\to} \overline{M} \to 0,$$ $D(-)$ is left exact.  Let $$E(M)={D(\overline{M})}/{\pi (D(M))}.$$  Then $E(-)$ is a functor on ${\cal M}_1(R)$, with $E(M)=0$ if $\dim M = 0$ or if $M$ is Cohen-Macaulay of dimension $1$.  $E(-)$ is neither left exact nor right exact. \end{Exercise}
\hint Use Proposition \ref{dil-dim1}. \QED

\section{Behavior Under Generic Hyperplane Sections}

In this section, we explore the behavior of the Dilworth number under generic hyperplanes sections.  We preserve the notation and assumptions of the previous section, and in addition assume that the field $k$ is infinite.

\begin{Lemma}\label{MovxM}  Suppose $\dim M \leq 1$.  If $\dil (M) = \dil (M/xM)$ for a generic hyperplane $x$, then $M$ is Dilworth, {\em i.e.},  $\nu(M) = \dil (M)$. \end{Lemma}

\proof Let $N$ be the unique maximal submodule such that $\nu(N) = \dil(M) = d$.  Let $x$ be generic.  Let $N_x$ be the preimage of $D(M/xM)$ in $M$.  Since $\dil(M) = \dil(M/xM),\, \nu(N_x) = d$.  Since $xM \subseteq N_x$, we must have $xM \subseteq \m N_x$.  On the other hand, since $N_x \subseteq N$, we have $xM \subseteq \m N$.  Since this is true for generic $x$, it follows that $\m M \subseteq \m N \subseteq N \subseteq M$.  Since $\nu(M) \leq \nu(N)$, this forces $M=N$. \QED

\begin{Corollary}\label{0colx} Suppose $\dim M = 0$.  If $\dil (M) = \dil [0: x]_M$ for a generic hyperplane $x$ then $\dil (M) = \type (M)$ \end{Corollary}

\proof Matlis duality; see Proposition \ref{matlis-duality}. \QED

\bigskip

The following ugly orphan of a lemma will be the crucial step in the proof of Proposition \ref{dil-structure}

\begin{Lemma}\label{orphan} If $M$ is a module over $R$ such that $\length(M) < \infty$ and $$\type(M) = \dil(M) = \nu([0:x]_M) = d$$ for a generic hyperplane $x$, then $[0:x]_M = \socle(M)$ for such $x$.  \end{Lemma}

\proof We proceed by induction on $\length(M)$.  For a hyperplane $x$, write $L_x = [0:x]_M$.  If $L_x = M$ for generic $x$, then $\m M = 0$ and the assertion is trivial.  So we may assume $L_x \neq M$ for generic $x$.  Let $x$ and $y$ be generic.  Since $\nu(L_x) = \nu(L_y) = d$, we have $\nu(L_x\cap L_y) = d$.  But since $[0:y]_{L_x} = L_x \cap L_y$, we have that the module $L_x$ satisfies the hypotheses of the lemma, and $\length(L_x) < \length(M)$. So, by induction, we have $L_x \cap L_y = \socle(M)$ for generic $x,y$.

Choose hyperplanes $x$ and $y$ such that $\nu(L_y) = \nu(L_{x+\alpha y}) = d$ for infinitely many units $\alpha \in R$.  For the remainder of this argument, we will call such $\alpha$ ``good''.  We shall also employ the following language: if $A$ is a finitely generated module and $a \in A$, then we say that $a$ is a {\em generator} of $A$ if $a$ is part of a minimal generating set for $A$; that is, if $a \in A \setminus \m A$.

\noindent {\bf Claim.} For good $\alpha$, if $L_{x+\alpha y} \neq \socle(M)$ then there exists a module $K_{\alpha}$ such that \begin{itemize}
\item $\socle(M) \subset K_{\alpha} \subseteq L_{x+\alpha y}$,
\item $yK_{\alpha} \subseteq \socle(M)$, and
\item $\nu(K_{\alpha}) = d$.
\end{itemize}
We will prove the claim presently.  For the moment, assume it is true.  For brevity, write $L_{\alpha}$ for $L_{x+\alpha y}$ and write $K_{\alpha} = A_{\alpha} + \socle(M)$, where $A_{\alpha}$ has no generators in $\socle(M)$.  (The choice of $A_{\alpha}$ is not unique, but no matter.)  We wish to show, by induction on $m$, that $\nu(A_{\alpha_1} + \cdots + A_{\alpha_m}) = \nu(A_{\alpha_1}) + \cdots + \nu(A_{\alpha_m})$ for distinct, good $\alpha_i$.  This is a contradiction, and so will complete the proof.

The case $m=1$ is trivial.  By induction, we may assume $\nu(A_{\alpha_1} + \cdots + A_{\alpha_{m-1}}) = \nu(A_{\alpha_1}) + \cdots + \nu(A_{\alpha_{m-1}})$. 
For good $\alpha,\, L_\alpha \cap L_y = \socle(M)$, so $\length(L_\alpha/y L_\alpha) = \length(L_\alpha \cap L_y) = d$.  Since $\nu(L_\alpha) = d$, this forces $yL_\alpha = \m L_\alpha$.  Now, we also have $\socle(M) \subseteq (L_{\alpha_1}+\ldots + L_{\alpha_{m-1}})\cap L_y$ and \begin{eqnarray*} 
\length((L_{\alpha_1}+\ldots + L_{\alpha_{m-1}})\cap L_y) & = & \length\left(\frac{L_{\alpha_1} + \cdots + L_{\alpha_{m-1}}}{y(L_{\alpha_1} + \cdots + L_{\alpha_{m-1}})}\right) \\
& = & \length\left(\frac{L_{\alpha_1} + \cdots + L_{\alpha_{m-1}}}{\m(L_{\alpha_1} + \cdots + L_{\alpha_{m-1}})}\right) = d.\end{eqnarray*}
So $(L_{\alpha_1}+\ldots + L_{\alpha_{m-1}})\cap L_y = \socle(M)$.  In particular, $y$ does not annihilate any generator of $A_{\alpha_1}+\cdots +A_{\alpha_{m-1}}$.  Therefore, if for each $1\leq i \leq m-1,\,v_i$ is a generator for $A_{\alpha_i}$, then $\{yv_1,\ldots,yv_{m-1}\}$ is a linearly independent set in the vector space $\socle(M)$.



Now, to show $\nu(A_{\alpha_1} + \cdots + A_{\alpha_m}) = \nu(A_{\alpha_1}) + 
\cdots + \nu(A_{\alpha_m})$, it suffices to show that if $v_i$ is a generator
 for $A_{\alpha_i}$, $1 \leq i \leq m$, and $$\sum_{i=1}^m \beta_i v_i 
\in \m (A_{\alpha_1} + \cdots + A_{\alpha_m})$$ then each $\beta_i \in \m$.

Multiplying by $x+\alpha_m y$, we obtain
\begin{eqnarray*}
(x+\alpha_m y)\sum_{i=1}^m \beta_i v_i & \in & \m (x+\alpha_m y)(A_{\alpha_1} + \cdots + A_{\alpha_m}) \\
\sum_{i=1}^{m-1}(\alpha_m - \alpha_i) \beta_i y v_i & \in & \m ((\alpha_m - \alpha_1)yA_{\alpha_1} + \cdots + (\alpha_m - \alpha_{m-1})yA_{\alpha_{m-1}})\\
& = & 0.\end{eqnarray*}

Since the $yv_i$ are linearly independent, $1\leq i \leq m-1$, and the $\alpha_i$ are distinct, this implies that $\beta_i \in \m$ for $1\leq i \leq m-1$.  Similarly, one shows that $\beta_m \in \m$.  Thus, we are done.

\noindent {\bf Proof of claim.}  It suffices to show that if $\socle(M) \subseteq K \subseteq L_{x+\alpha y}$, with $\nu(K) = d$ and $K' := yK+\socle(M)$, then $\nu(K') = d$ as well.  Indeed, if this is true, then we can find our $K_\alpha$ among the modules $L_{x+\alpha y}, L'_{x+\alpha y},L''_{x+\alpha y},L'''_{x+\alpha y}, \ldots$

Thus, let $K=A+\socle(M)$, where no generator of $A$ lies in $\socle(M)$.  Then $K' = yA + \socle(M)$. Now, we have
$$\nu(K) = \nu(A) + d-\dim(A\cap \socle(M))$$
and
$$\nu(K') = \nu(yA) + d - \dim(yA \cap \socle(M)).$$
Since $L_{\alpha}\cap L_y = \socle(M)$ and since no generator of $A$ lies in $\socle(M)$, we have $\nu(yA) = \nu(A)$.  Since, clearly, $yA\cap \socle(M) \subseteq A \cap \socle(M)$, this yields $\nu(K') \geq \nu(K)$.  Since $\nu(K) = \dil(M) = d,\, \nu(K') = d$ and we are done. \QED


\bigskip

We shall also require the dual formulation of this lemma, and so we record it here.
\begin{Lemma}\label{orphandual} If $M$ is a module over $R$ such that $\length(M) < \infty$ and $$\nu(M) = \dil(M) = \type(M/xM) = d$$ for a generic hyperplane $x$, then $M/xM = M/\m M$ for such $x$.  \end{Lemma}

\section{The Definition and Basic Properties of the Homological Dilworth Number} \label{hdil-sec}

As previously  remarked, Equation (\ref{dil-dim1-eq}) suggests an analogy between the Dilworth number and the homological degree.  To extend the Dilworth number to modules of dimension greater than or equal to $2$, we are therefore led to the following definition.

\begin{Definition}\label{hdil-def} If $R$ is the homomorphic image of a Gorenstein ring and
$\dim M = d$, then the homological Dilworth number is defined recursively as follows,
$$ \hdil (M) = \left\{ \begin{array}{ll}
                      \deg M + \displaystyle\sum_{i=0}^{d-1}{d-1 \choose i}\hdil(M_i) & d>0 \\
                      \dil(M) & d=0 \end{array} \right. . $$ 
If $R$ is not the homomorphic image of a Gorenstein ring, then we define 
$$ \hdil(M) = \hdil(M \otimes_R \hat{R}),$$ where $\hat{R}$ is the $\m$-adic completion of $R$.
\end{Definition}

By Equation (\ref{dil-dim1-eq}) we have that
$$\hdil(M) = \dil(M)$$ whenever $\dil(M)<\infty$; that is, whenever $\dim M \leq 1$.  

The homological Dilworth number is not a cohomological degree, as it violates the condition $$\Deg(M) = \length(M),$$ whenever $\dim M = 0$.  However, if $\dim M > 0$ and $M$ is Cohen-Macaulay, then $\hdil(M) = \deg(M)$.  Suspecting that the difficulty is confined to dimension $0$, one may bravely speculate that there exists some cohomological degree, $\Deg(-)$, which agrees with the homological Dilworth number, at least for modules of positive depth;  but this too is false, as we shall see.  Nonetheless, many of the statements which hold for $\hdeg(-)$ may also be proven for $\hdil(-)$, a task to which we now dedicate ourselves.  Generally speaking, a statement about $\hdeg(-)$ and the corresponding statement about $\hdil(-)$ have similar proofs,
except on the level of details (the home of the Devil;  see, for example, Lemma \ref{orphan}), where the proofs of statements about $\hdil(-)$ tend to be more complicated.

\begin{Lemma}\label{triangle}If $$ 0 \lar N \lar M \lar U \lar 0$$ is exact, then \begin{enumerate}
\item If $\length(N) < \infty$ or $\length(U) < \infty$, then $\hdil(M) \leq \hdil(N) + \hdil(U)$.
\item If $N=\gm(M)$, then $\hdil(M) = \dil(N) + \hdil(U)$. \end{enumerate}\end{Lemma}

\proof By taking completions if necessary, we may assume that $R$ is the homomorphic image of a Gorenstein ring.  Suppose $\length(U)<\infty$.  It suffices to show that $\hdil(M_i) \leq \hdil(N_i)$ for $i \geq 1$ and $\dil(M_0) \leq \dil(N_0) + \dil(U_0)$.  But $M_i \cong N_i$ for $i>1$ and we have
$$ 0 \lar M_1 \lar N_1 \lar U_0 \lar M_0 \lar N_0 \lar 0.$$
Since, by Proposition \ref{h-Spec}, $\dim N_1 = \dim M_1 \leq 1$, we have that $\hdil (-) = \dil (-)$ for these modules, and since $M_1 \subseteq N_1$, we have that $\dil (M_1) \leq \dil (N_1)$.  

Meanwhile, if we let $A = (\ker M_0 \onto N_0)$, then we have, by Lemma~\ref{diltriangle}, $\dil(M_0) \leq \dil(N_0) + \dil(A)$ and $\dil(A) \leq \dil(U_0)$.  Hence we are done in this case.

Now suppose that $\length(N)<\infty$.  Then $M_i \cong U_i$ for $i \geq 1$ and we have
$$0 \lar U_0 \lar M_0 \lar N_0 \lar 0.$$  Thus $\dil(M_0) \leq \dil(U_0)+\dil(N_0)$.  Furthermore, if $N=\gm(M)$, we have $U_0 = 0$.  Hence we are done in all cases. \QED

\bigskip The following proposition was proved for the homological degree by Vasconcelos \cite{V96}.  The proof given here follows the argument of Vasconcelos closely; we merely observe that the properties of the Dilworth number are just enough to carry through.

\begin{Proposition}\label{inductive}  If $x$ is a generic hyperplane on $M$, then, for $r \geq 1$, we have
$$r\cdot\hdil(M)  \geq  \hdil(M/x^rM). $$ In particular,
$$\hdil(M)  \geq  \hdil(M/xM). $$

\end{Proposition}

\proof We may assume that $R$ is the homomorphic image of a Gorenstein ring.  We proceed by induction on $d=\dim M$.  

First, observe that we may assume that $\depth (M) > 0$.  Indeed, assume the result has been proven for modules of positive depth and suppose $\depth(M)=0$.

If not, then the exact sequence 
$$ 0 \lar \gm(M) \lar M \lar M' \lar 0$$
yields
$$ 0 \lar  \gm(M)/x^r\gm(M) \lar M/x^rM \lar M'/x^rM' \lar 0.$$
Hence, by Proposition \ref{triangle}, \begin{eqnarray*}
\hdil(M/x^rM) & \leq & \dil(\gm(M)/x^r\gm(M)) + \hdil(M'/x^rM') \\
& \leq & \dil(\gm(M)) + r \cdot  \hdil(M') \\
& \leq & r \cdot  \hdil(M). \end{eqnarray*}

So we assume that $x^r$ is regular on $M$.  We have the exact sequence
$$0 \lar M \stackrel{x^r}{\lar} M \lar \overline{M} \lar 0,$$
from which we obtain the long exact sequence
$$ 0 \lar M_d \stackrel{x^r}{\lar} M_d \lar \overline{M}_{d-1} \lar \cdots \lar M_1 \stackrel{x^r}{\lar} M_1 \lar \overline{M}_0 \lar 0.$$
Since $x$ is generic, we may take it to be part of a degree system for $M$, and so, by Proposition~\ref{lastone?}, $r \cdot  \deg(M) = \deg(\overline{M})$. Thus it suffices to show the inequality
\begin{equation}\label{ineq-goal} \hdil(\overline{M}_{i-1}) \leq r(\hdil(M_i) + \hdil(M_{i-1})).\end{equation}

With this in mind, we break the long exact sequence above into short exact sequences as follows.
$$ 0 \lar L_i \lar M_i \lar \widetilde{M}_i \lar 0 $$
\begin{equation}\label{b} 0 \lar \widetilde{M}_i \lar M_i \lar G_i \lar 0 \end{equation}
\begin{equation}\label{c} 0 \lar G_i \lar \overline{M}_{i-1} \lar L_{i-1} \lar 0 \end{equation}
Since $x$ is generic, we have that $\dim L_i = 0$.
If $\dim M_i = 0$, then we have \begin{eqnarray*} 
\hdil(\overline{M}_{i-1}) = \dil(\overline{M}_{i-1}) & \leq & \dil(G_i) + \dil(L_{i-1}) \\
& \leq & \dil(M_i) + \dil(M_{i-1}) \end{eqnarray*}
Since $r \geq 1$, this immediately implies (\ref{ineq-goal}).

If $\dim M_i > 0$, we begin by applying the functor $\gm(-)$ to (\ref{b}), from which we obtain
\begin{equation}\label{diagram} \begin{array}{ccccccccc}
0 & \lar & \widetilde{M}_i & \lar & M_i & \lar & G_i & \lar & 0 \\
&& \uparrow & & \uparrow & & \uparrow \\
0 & \lar & \gm(\widetilde{M}_i) & \lar & \gm(M_i) & \lar & \gm(G_i) \end{array} \end{equation}
Let $H_i = \image(\gm(M_i) \lar \gm(G_i))$.  Then applying the snake lemma to (\ref{diagram}) yields
\begin{equation} \label{d}0 \lar \widetilde{M}_i/\gm(\widetilde{M}_i) \lar M_i/\gm(M_i) \lar G_i/H_i \lar 0. \end{equation}
Now, $\widetilde{M}_i/\gm(\widetilde{M}_i) \simeq M_i/\gm(M_i)$, and the composite map
$$M_i/\gm(M_i)\stackrel{\sim}{\lar}\widetilde{M}_i/\gm(\widetilde{M}_i) \lar M_i/\gm(M_i)$$ is induced by multiplication by $x^r$.  Hence we may rewrite (\ref{d}) as
$$0 \lar M_i/\gm(M_i) \stackrel{x^r}{\lar} M_i/\gm(M_i) \lar G_i/H_i \lar 0.$$
Since $\dim (M_i/\gm(M_i)) \leq i < \dim M$, by Proposition \ref{h-Spec}, we may apply induction.  Hence, by Proposition \ref{triangle}, we have 
\begin{eqnarray*} \hdil(G_i/H_i) & \leq & r \cdot  \hdil(M_i/\gm(M_i))\\
& = & r\cdot  \hdil(M_i) - r\cdot \dil(\gm(M_i)). \end{eqnarray*}
Now, from (\ref{c}), we have
\begin{eqnarray*}\hdil(\overline{M}_{i-1}) & \leq & \dil(L_{i-1})+\hdil(G_i) \\
& \leq & \hdil(M_{i-1}) + \hdil(G_i/H_i) + \dil(H_i) \\
& \leq & \hdil(M_{i-1}) + \dil(H_i) + r \cdot  \hdil(M_i) - r \cdot  \dil(\gm(M_i)). \end{eqnarray*}
Since $\dil(H_i) \leq \dil(\gm(M_i))$, we have $\hdil(\overline{M}_{i-1}) \leq r\Big(\hdil(M_{i-1}) + \hdil(M_i)\Big).$ \QED

\begin{Corollary}\label{moo} $\nu(M) \leq \hdil(M)$ \end{Corollary}

\begin{Definition} If $\nu(M) = \hdil(M)$ then we say that $M$ is Dilworth. \end{Definition}

Since, clearly, $\hdil(M) \leq \hdeg(M)$, a Yoshida module (that is, a module which satisfies $\nu(M) = \hdeg(M)$) is necessarily Dilworth.  The converse is not true, even if $\depth(M)>0$, as the following example shows.

\begin{Example} A module which is Dilworth but not Yoshida.\end{Example}
Let $R = k [\![ x,y,z ] \! ]$ or $k[x,y,z]_{(x,y,z)}$, and let $\m = (x,y,z)$.  Let $M = \m^tR/x^tR$.  Then $M$ is Dilworth but not Yoshida, with
$$\nu(M) = \hdil(M) = {t+2 \choose 2} - 1$$
and
$$\hdeg(M) = {t+2 \choose 3} + t.$$
(We leave the calculation as an exercise for the reader.)  Observe that $\grm(M)$ has a linear resolution (which foreshadows the results in section \ref{linres-dil}) and hence, in this case, $\hdil(M)=\bdeg(M)$.

\begin{Exercise} Given an exact sequence of finitely generated modules
$$0 \lar A \lar B \lar C \lar 0$$ with $\m C = 0$ and $\hdil(B) = \hdil(A) + \length(C)$, then
$$0 \lar \gm(A) \lar \gm(B) \lar C \lar 0$$ is also exact, and $\dil(\gm(B)) = \dil(\gm(A)) + \length(C)$. \end{Exercise}

\section{Bounding Numbers of Generators and Hilbert Functions}

The following three propositions are the analogues of Propositions \ref{gen-a}, \ref{gen-b} and \ref{gen-c}.  The proofs of Propositions \ref{hdil-gen-a} and \ref{hdil-gen-c} are almost identical to the proofs of Propositions \ref{gen-a} and \ref{gen-c}, and so we omit them.  The proof of \ref{hdil-gen-b} is slightly different from that of \ref{gen-b} (in fact, Proposition \ref{hdil-gen-b} is marginally simpler to prove than Proposition \ref{gen-b}, an exception to the rule stated in the remarks following Definition \ref{hdil-def}) so we will include the proof in this case.

\begin{Proposition} \label{hdil-gen-a} Suppose $R$ is Cohen-Macaulay of dimension $d$ and let $I$ be an ideal of codimension $g > 0$.  If $\depth (R/I) = r$, then
$$\nu(I) \leq \deg(R) + (g-1)\deg(R/I) + (d-r-1)(\hdil(R/I) - \deg(R/I)).$$\end{Proposition}
\proof Similar to \cite[Theorem 3.1]{DGV}, cited as Proposition \ref{gen-a}. \QED

\begin{Proposition} \label{hdil-gen-b} Let $\dim R = d > 0$ and let $I$ be an $\m$-primary ideal.  Then
$$ \nu(I) \leq \hdil(R) {s+d-2 \choose d-1} + {s+d-2 \choose d-2} ,$$
where $s$ is the index of nilpotency of $R/I$. \end{Proposition}

\proof Let $x_1,\ldots,x_d$ be a generic system of parameters for $R$ and let $L=(x_1,\ldots,x_{d-1})$.  Since $L^s \subseteq I$, we have $\nu(I) \leq \nu(L^s) + \nu(I/L^s)$.  Clearly, $\nu(L^s) \leq {s+d-2 \choose d-2}$.  Since $I/L^s$ is an ideal of $R/L^s$, a one dimensional ring, it suffices to
show that $\dil(R/L^s) \leq \hdil(R) {s+d-2 \choose d-1}$.

We use induction on $s$.  If $s = 1$, the claim follows from Proposition \ref{inductive}.  Suppose $s > 1$.  The exact sequence
$$ 0 \lar L^{s-1}/L^s \lar R/L^s \lar R/L^{s-1} \lar 0$$
shows that 
$$\dil(R/L^s) \leq \dil(L^{s-1}/L^s) + \dil(R/L^{s-1}).$$
Let $r = {s+d-3 \choose d-3}$.  Since $\nu(L^{s-1}/L^s) \leq r$, we have $(R/L)^r \onto L^{s-1}/L^s$, whence $$\dil(L^{s-1}/L^s) \leq r\cdot \dil(R/L) \leq r \cdot \hdil(R).$$  On the other hand, by the induction hypothesis we have $$\dil(R/L^{s-1}) \leq \hdil(R){s+d-3 \choose d-2}.$$  Putting this all
together, we have
$$\dil(R/L^s) \leq r \cdot \hdil(R) + \hdil(R){s+d-3 \choose d-2} \leq \hdil(R) {s+d-2 \choose d-1}.$$ \QED

\begin{Proposition} \label{hdil-gen-c} The reduction number $\red(\m)$ of $\m$ is bounded as follows. $$\red(\m) \leq d\cdot \hdil(R) -2d + 1.$$ \end{Proposition}

\proof Similar to \cite[Corollary 4.7]{DGV}, cited as Proposition \ref{gen-c}. \QED

\section{Dilworth Modules and Linear Resolution}\label{linres-dil}

In this section, we study Dilworth modules: modules which satisfy $$\nu(M) = \hdil(M).$$  The analogous statement for a cohomological degree $\Deg(-)$ is $$\nu(M) = \Deg(M),$$ and, as we saw in Chapter 2, this equality forces $\grm(M)$ to have a linear resolution.  Therefore, one is moved to ask:
\begin{quote} If $M$ is Dilworth, does $\grm(M)$ necessarily have a linear resolution?\end{quote} The answer is no, for in dimension $0$, the question becomes
\begin{quote} Does $\nu(M) = \dil(M)$ imply $\m M = 0$?\end{quote} (The reader may satisfy herself that counterexamples are myriad.)  What is perhaps surprising is that misbehavior in dimension $0$ is the {\em only} obstruction to linear resolution.  Thus, we will prove:

\begin{Theorem}\label{main} Let $M$ be a Dilworth module such that $\m \gm(M) = 0$.  Then $\grm(M)$ has a linear resolution. \end{Theorem}

Before we give a proof, we need to establish some preliminary results.

\begin{Proposition}\label{tootie} If $M$ is Dilworth and $x$ is a generic hyperplane on $M$, then $\overline{M} = M/xM$ is Dilworth also. \end{Proposition}

\proof Squeeze:  $\nu(M) = \hdil(M) \geq \hdil(\overline{M}) \geq \nu(\overline{M}) = \nu(M)$. \QED

\begin{Proposition}\label{agl} Assume $R$ is the homomorphic image of a Gorenstein ring.  Let $M$ be a module such that $\dim M = d$ and $ \depth (M) > 0$.  If $M$ is Dilworth, then $M_{d-1}$ is Dilworth and, for $i \leq d-2$, we have $\dim M_i = 0$ and $\nu(M_i) = \type(M_i) = \dil(M_i)$. \end{Proposition}

\proof Let $x$ be generic on $M$.  We write $\overline{M}$ for $M/xM$ and $\overline{M}_i$ for $\Ext_S^{n-i}(\overline{M},S)$.  The exact sequence
$$0 \lar M \stackrel{x}{\lar} M \lar \overline{M} \lar 0$$
yields
$$M_{i+1} \stackrel{x}{\lar} M_{i+1} \lar \overline{M}_i \lar M_i \stackrel{x}{\lar} M_i,$$
or
$$ 0 \lar \frac{M_{i+1}}{xM_{i+1}} \lar \overline{M}_i \lar [0:x]_{M_i}) \lar 0,$$
for each $i$.  Since $x$ is generic, the module  $[0:x]_{M_i}$ has finite length, and so, by Lemma \ref{triangle}, we have
\begin{eqnarray*}
\hdil(\overline{M}_i) & \leq & \hdil\left(\frac{M_{i+1}}{xM_{i+1}}\right) + \dil([0:x]_{M_i})\\
                      & \leq & \hdil(M_{i+1}) + \hdil(M_i),\end{eqnarray*}
for each $i$.  On the other hand, we have
\begin{eqnarray*}
\hdil (\overline{M}) & = & \deg (\overline{M}) + \sum_{i=0}^{d-2}{d-2 \choose i}\hdil(\overline{M}_i) \\
                     & \leq & \deg (M) + \sum_{i=0}^{d-2}{d-2 \choose i}(\hdil(M_{i+1}) + \hdil(M_i)) \\
                     & =    & \deg (M) + \hdil(M_{d-1}) + \sum_{i=0}^{d-2}\left[{d-2 \choose i}+{d-2 \choose i-1}\right]\hdil(M_i) \\
                     & =    & \deg (M) + \sum_{i=0}^{d-1}{d-1 \choose i}\hdil(M_i) \\
                     & =    & \hdil(M) = \hdil(\overline{M}).\end{eqnarray*}
Hence we must have equality all the way through, and so
\begin{eqnarray*}
\hdil(\overline{M}_i) & = & \hdil\left(\frac{M_{i+1}}{xM_{i+1}}\right) + \dil([0:X]_{M_i}), \\
\hdil\left(\frac{M_{i+1}}{xM_{i+1}}\right) & = & \hdil(M_{i+1}), \mbox{ and } \\
\dil([0:x]_{M_i}) & = & \hdil(M_i), \end{eqnarray*}
for $i\leq d-2$.  This forces $\dim M_i = 0$ (since $[0:x]_{M_i} \subseteq \gm(M_i)$) and, by Lemma~\ref{MovxM} and Corollary~\ref{0colx}, $\nu(M_i) = \dil(M_i) = \type(M_i)$, for $i \leq d-2.$  To show the final claim, that $\nu(M_{d-1}) = \hdil(M_{d-1})$, we proceed by induction on $d$.

If $d=2$, in which case $M_{d-1} = M_1$, we have $\hdil(M_1) = \hdil({M_1}/{xM_1})$ for generic $x$, and since $\dim M_1 \leq 1$, we are done by Lemma~\ref{MovxM}.

If $d>2$, then we have, by induction, \begin{eqnarray*}
\hdil(\overline{M}_{d-2}) & = & \nu(\overline{M}_{d-2}) \\
                          & \leq & \nu\left(\frac{M_{d-1}}{xM_{d-1}}\right) + \nu([0:x]_{M_{d-2}}) \\
                          & \leq & \nu\left(\frac{M_{d-1}}{xM_{d-1}}\right) + \nu(M_{d-2}) \\
                          & \leq & \hdil\left(\frac{M_{d-1}}{xM_{d-1}}\right) + \dil(M_{d-2}) \\
                          & \leq & \hdil(M_{d-1}) + \dil(M_{d-2}) \\
                          & =    & \hdil(\overline{M}_{d-2}). \end{eqnarray*}
Hence we have equality all the way through, whence
$$ \nu(M_{d-1}) = \nu\left(\frac{M_{d-1}}{xM_{d-1}}\right) = \hdil\left(\frac{M_{d-1}}{xM_{d-1}}\right) = \hdil(M_{d-1}),$$
which finishes the proof. \QED

\begin{Lemma} \label{h0} Let $M$ be a Dilworth module with $L = \gm(M)$ and $N=M/L$.  Then $L$ and $N$ are both Dilworth modules with $\nu(M) = \nu(L) + \nu(N)$.  Furthermore, if $\m L = 0$, then $M \cong L \oplus N$. \end{Lemma}

\proof  It is harmless to assume that $R$ is the homomorphic image of a Gorenstein ring.  From the exact sequence
\begin{equation}\label{tosplit} 0 \lar L \lar M \lar N \lar 0 \end{equation}
we have
\begin{equation}\label{e} \nu(L) + \nu(N) \leq \dil(L) + \hdil(N) = \hdil(M) = \nu(M) \leq \nu(L) + \nu(N). \end{equation}
Hence $\nu(M) = \nu(L) + \nu(N)$ and $L$ and $N$ are again Dilworth.  On the other hand, if $\m L = 0$, tensoring (\ref{tosplit}) with $R/ \m$ yields
\begin{equation}\label{f} L \lar \overline{M} \lar \overline{N} \lar 0. \end{equation}
By (\ref{e}), the map $L \lar \overline{M}$ must be an injection, so (\ref{f}) becomes
$$ 0 \lar L \lar \overline{M} \lar \overline{N} \lar 0.$$
This induces a splitting of (\ref{tosplit}); {\em i.e.} $M \cong L \oplus N$. \QED

\begin{Proposition}\label{dil-structure} Let $M$ be Dilworth of dimension $d$ and of positive depth.  Let $x$ be a generic, regular hyperplane on $M$, and write $\overline{M} = M/xM$.  Then $\m H_{\m}^i(\overline{M}) = 0$ for $i \leq d-3$.  Furthermore, if $d=2$, then $\m \gm(\overline{M}) = 0$. \end{Proposition}

\proof We may assume that $R$ is the homomorphic image of a Gorenstein ring. As in the proof of Proposition \ref{agl}, we have, for $i \leq d-3$, short exact sequences of modules of finite length
$$ 0 \lar M_{i+1}/xM_{i+1} \lar \overline{M}_i \lar [0:x]_{M_i} \lar 0 $$
subject to
\begin{eqnarray}
\label{ff} \dil(\overline{M}_i) & = & \dil(M_{i+1}/xM_{i+1}) + \dil([0:x]_{M_i}) \\
\label{g} \dil(M_{i+1}/xM_{i+1}) & = & \dil(M_{i+1})\\
\label{h} \dil([0:x]_{M_i}) &  = & \dil(M_i) \end{eqnarray}
Since $\nu(M_i) = \type(M_i) = \dil(M_i)$ for $i\leq d-2$, equations (\ref{g}) and (\ref{h}), along with Lemmas \ref{orphan} and \ref{orphandual}, force $\m[0:x]_{M_i} = \m(M_{i+1}/xM_{i+1}) = 0$ for $i \leq d-3$.  By virtue of (\ref{ff}), we have
$$\overline{M}_i \cong \left( \frac{M_{i+1}}{xM_{i+1}} \right) \bigoplus [0:x]_{M_i},$$
whence $\m \overline{M}_i = 0$, for $i \leq d-3$.

To show the second claim, first note that $\overline{M}_0 = M_1/xM_1$.  If $\dim M_1 = 0$, then we may proceed as above, so we assume $\dim M_1 = 1$.  $\overline{M}$ is Dilworth, so $\nu(\gm(\overline{M})) = \dil(\gm(\overline{M}))$, by Lemma \ref{h0}.  Since $(\overline{M}_0)^{\vee} \cong \gm(\overline{M})$, we have $\type(M_1/xM_1) = \dil (M_1/xM_1)$.

Let $L = \gm(M_1)$ and $M_1' = M/L$.  Note that since $M_1$ is Dilworth, $L$ is Dilworth and $M_1'$ is a one dimensional linear Cohen-Macaulay module.  The exact sequence
$$0 \lar L \lar M_1 \lar M_1' \lar 0$$ yields
$$0 \lar L/xL \lar M_1/xM_1 \lar M_1'/xM_1' \lar 0.$$
Since $M_1'$ is linear Cohen-Macaulay, $\m (M_1'/xM_1') = 0.$  Further,
\begin{eqnarray*}
\nu(M_1) & \leq & \dil(M_1/xM_1) = \type (M_1/xM_1) \\
& \leq & \type(L/xL) + \type(M_1'/xM_1') \\
& \leq & \dil(L) + \deg(M_1') \\
& \leq & \dil(M_1) = \nu(M_1), \end{eqnarray*}
so $\type(L/xL) = \dil(L) = \nu(L)$.  Consequently, by Lemma \ref{orphandual}, $\m (L/xL) = 0$.  Meanwhile, the above calculation also forces $\dil(M_1/xM_1) = \ell(L/xL) + \ell(M_1'/xM_1')$; hence $$\frac{M_1'}{xM_1'} \cong \left( \frac{L}{xL} \right) \bigoplus \left( \frac{M_1'}{xM_1'} \right),$$ which finishes the proof.  \QED

\bigskip

By replacing the condition that $\nu(M) = \Deg(M)$ by the conditions $\m \gm(M) = 0$ and $\nu(M) = \hdil(M)$, we may prove the following three Lemmas exactly as Lemmas \ref{QM}, \ref{SV-lemma} and \ref{exact}.

\begin{Lemma}\label{QMdil} Let $M$ be Dilworth with $\m \gm(M) = 0$ and $Q = (x_1,\ldots,x_d)$, where $x_1,\ldots,x_d$ is a generic system of parameters for $M$.  Then $\m M = QM$. \end{Lemma}

\begin{Lemma}\label{SVlemma} Let $I$ be an ideal of $R$ and $M$ a finitely generated $R$-module such that $M/IM$ is Dilworth and $\m \gm(M/IM)=0$.  Let $Q = (x_1,\ldots,x_r)$, where $x_1,\ldots,x_r$ is a generic partial system of parameters for $M/IM$.  Then 
$$ [IM:\m]\cap Q^kM \subseteq IQ^{k-1}M \qquad \mbox{ for } k\geq 1.$$ \end{Lemma}

\begin{Lemma} Let $M$ be Dilworth with $\depth(M)>0$ and let $x$ be a regular, generic hyperplane on $M$.  Set $\overline{M} = M/xM$.  Then 
$$ 0 \lar \grm(M) \stackrel{x^*}{\lar} \grm(M) \lar \grm(\overline{M}) \lar 0$$ is exact, where $x^{*}$ is an initial form of $x$ in $\grm(R)$. \end{Lemma}

\noindent {\bf Proof of Theorem \ref{main}.} We have shown (Propositions \ref{tootie} and \ref{dil-structure}) that if $M$ is Dilworth and $\m \gm(M) = 0$ and $x$ is generic on $M$, then $M/xM$ is Dilworth and $\m \gm(M/xM) = 0$; that is, that the hypotheses of Theorem \ref{main} are preserved under a generic hyperplane section.  Thus we may prove Theorem \ref{main} exactly as we proved Theorem \ref{Yosh}. \QED

\bigskip

\noindent The converse to Theorem \ref{main} is not true, as the following example shows.

\begin{Example} A module which has a linear resolution but is not Dilworth.\label{example}\end{Example} Let $(R,\m)$ be a regular local ring of dimension $n \geq 2$.  Let $M= \m^{t+1}$.  Then $\nu(M) = {n+t \choose n-1}$.  Since $M_1 \cong (R/M)^{\vee}$ and $M_i=0$ for $i=0$ and for $2\leq i \leq n-1$,  $$\hdil(M) = 1 + (n-1)\dil(R/M) = 1 + (n-1){n+t-1 \choose n-1}.$$  If $n=2$, then $M$ is Dilworth and if $n > 2$, the difference $\hdil(M) - \nu(M)$ is a polynomial in $t$ of degree $n-1$.  In either case, $\grm(M)$ has a linear resolution, by Proposition \ref{trunc-prop}.

\section{Comparison Between $\bdeg(-)$ and $\hdil(-)$}

In this section we compare the homological Dilworth number with the extremal cohomological degree.  We preserve the notation and assumptions of the previous sections, including the assumption that $k$ is infinite. Our first result is:

\begin{Proposition} \label{hdil-bdeg-dim2} Let $M$ be in ${\cal M}(R)$ and suppose $\dim M \leq 2$.  Then $$\hdil(M) \leq \bdeg(M).$$ \end{Proposition}

\noindent We shall require the following result of Ikeda \cite[Theorem 1.1]{Ik2}.

\begin{Proposition} \label{rees} If $\dim M \leq 1$ then $$\dil(M) \leq \length(M/xM),$$ for all $x \in \m$. \end{Proposition}

\noindent {\bf Proof of Proposition \ref{hdil-bdeg-dim2}.}  By taking completions, we may assume that $R$ is the homomorphic image of a Gorenstein ring.  The proposition is almost trivial is $\dim(M) \leq 1$, so we assume that $\dim(M)=2$.  Let $L = \gm(M)$.  Then, by Lemma \ref{triangle}(2) and by Definition \ref{bdeg-def},
$$\hdil(M) = \dil(L) + \hdil(M/L)$$ and $$\bdeg(M) = \length(L) + \bdeg(M/L).$$ Since $\dil(L) \leq \length(L)$, we may assume $\depth(M) > 0$.  So $$\hdil(M) = \deg(M) + \dil(M_1).$$
Let $x$ be generic on $M$ and let $\overline{M} = M/xM$.  Then, since $\dim \overline{M} = 1$, $$\bdeg(M) = \bdeg(\overline{M}) =  \deg(\overline{M})+\length(\gm(\overline{M})).$$  But $ \gm(\overline{M}) \cong (\overline{M}_0)^{\vee}$ and $\overline{M}_0 \cong M_1/xM_1$.  Thus $$\bdeg(M) =  \deg(M) + \length(M_1/xM_1).$$ By Proposition \ref{rees}, $\length(M_1/xM_1) \geq \dil(M_1)$, and so we are done. \QED

\begin{Corollary} \label{dil-dim2} Suppose $M$ in ${\cal M}(R)$ has dimension at most $2$.  Then
$\grm(M)$ has a linear resolution if and only if $M$ is Dilworth and $\m \gm(M) = 0$. \end{Corollary}

\proof This follows from Theorems \ref{linresthm} and \ref{main} and Proposition \ref{hdil-bdeg-dim2}. \QED

\bigskip

\bigskip For another proof (in a special case) of the above theorem which does not make reference to $\bdeg(-)$, see Appendix \ref{App}.

We now turn our attention to the modules of dimension $3$ or greater.  It turns out that in this case, either inequality
$$\hdil(M) < \bdeg(M)$$ or $$\hdil(M) > \bdeg(M)$$ is possible (though, obviously, not at the same time).  First, observe that if $L=\gm(M)$ then
$$\hdil(M) = \dil(L) + \hdil(M/L)$$ and $$\bdeg(M) = \length(L) + \bdeg(M/L).$$ Since $\length(L) \geq \dil(L)$, and since the difference can be made as large as desired, it is in some sense {\em cheating} to admit modules with non-trivial $\gm(M)$.  Therefore, for the remainder of this section, the module $M$ is assumed to have dimension at least $3$ and positive depth.

Under these assumptions, the inequality $$\hdil(M) > \bdeg(M)$$ is typical.  Any module with a linear resolution which is not Dilworth, as in Example \ref{example}, must satisfy this inequality.  In fact, as the reader may verify for himself, randomly\footnote{In the sense of human activity, not in the mathematical sense.} constructing non-Cohen-Macaulay modules on, for example, {\sc Macaulay}, will almost invariably produce modules satisfying this inequality.

Thus, to produce a module satisfying $\hdil(M) < \bdeg(M)$ a systematic approach is required.  Suppose $R$ is Cohen-Macaulay and let $L$ be a module of finite length with $\nu(L) = m$.  Let $M$ be the module of syzygies of $L$: $$0 \lar M \lar R^m \lar L \lar 0.$$  Then $\dim M = 3, \, \deg(M)=m, \, H_{\m}^1(M) = L$ and $\gm(M) = H_{\m}^2(M) = 0$.  Therefore $$\hdil(M) = m + \dil(L)$$ and, as the reader may verify,
$$\bdeg(M) = m + \length(L/xL) + \length(L/(x,y)L),$$ where $x,y$ is a generic partial system of parameters.  Thus, we wish to find a module of finite length satisfying
$$\dil(L) < \length(L/xL)$$ for generic $x$.

We begin with some remarks concerning modules $L$ satisfying $\m^2L=0$.  Let $\nu(L)=m$ and $\nu(\m L) = n$.  Assume that no generator of $L$ is killed by $\m$; that is, that $$[0:\m]_L \subseteq \m L.$$  Thus $\length(L) = m+n$.  Let $V=L/\m L$ and $W=\m L$. Then we have a homomorphism $$\phi_L \in \Hom(\m / \m^2, \Hom(V,W))$$ given by $$\phi_L(\overline{x}) = (\overline{v} \mapsto xv),$$
where $x \in \m,\, v\in L,$ and $\overline{x}$ and $\overline{v}$ denote their residue classes in $\m / \m^2$ and $V$, respectively.  By fixing bases of $\m / \m^2, \, V$ and $W,\, \phi_L$ corresponds to an $m\times n$ matrix of linear forms, as follows.  Let $x_1,\ldots,x_t$ be a minimal generating set of $\m$. Then for each $i,\, \phi_L(\overline{x}_i) \in \Hom(V,W)$ corresponds to a matrix, and if $\overline{x} = \hat{x}_1 \overline{x}_1 + \cdots + \hat{x}_t \overline{x}_t$, then $$\phi_L(\overline{x}) = \sum_{i=1}^t \hat{x}_i \cdot \phi_L(\overline{x}_i),$$ an $m \times n$ matrix of linear forms in the variables $\hat{x}_1,\ldots,\hat{x}_t$.

Conversely, any matrix of linear forms $\phi \in \Hom(\m / \m^2, \Hom(V,W))$ determines a module $L_\phi$ with $\m^2 L_{\phi} = 0, \, L_{\phi}/\m L_{\phi} \cong V$ and $\m L_{\phi} \cong W$, subject to relations $$ x \cdot v = \phi(\overline{x})(\overline{v}).$$

Now let us use this correspondence to construct a module $L$ which satisfies $\m L^2 = 0$ and $\dil(L) < \length(L/xL)$ for generic $x$.  For convenience, let us assume that $\nu(L) = \length(\m L) = n$; that is, that $\phi = \phi_L$ corresponds to a square matrix.  We have the following proposition.

\begin{Proposition} With notation as above,
\begin{enumerate} \item $\length(L/xL) > \nu(L)$ for generic $x$ if and only if $\det(\phi) = 0$; and 
\item $\dil(L) = \nu(L)$ if and only if there is no choice of bases of $V$ and $W$ such that, with respect to these bases, the matrix of linear forms corresponding to $\phi$ has a rectangular block of zeroes in the lower left hand corner which is large enough to intersect the main diagonal; that is, of size $a \times b$ where $a+b > n$.
\end{enumerate}
\end{Proposition}

We leave the proof as an exercise to the reader.  Now, if $\m = (x,y,z)$, then the matrix $$ \phi = \left[ \begin{array}{ccc} 
0 & \hat{x} & \hat{y} \\
-\hat{x} & 0 & \hat{z} \\
-\hat{y} & -\hat{z} & 0 \end{array} \right]$$
satisfies $\det(\phi) = 0$, and it is clear that no sequence of row and column operations can get $\phi$ into either of the following forms:
$$
\phi = \left[ \begin{array}{ccc} 0 & \ast & \ast \\ 0 & \ast & \ast \\ 0 & \ast & \ast \end{array} \right],\,
\phi = \left[ \begin{array}{ccc} \ast & \ast & \ast \\ 0 & 0 & \ast \\ 0 & 0 & \ast \end{array} \right],\, \mbox{ or } \quad
\phi = \left[ \begin{array}{ccc} \ast & \ast & \ast \\ \ast & \ast & \ast \\ 0 & 0 & 0 \end{array} \right].$$
Thus $\nu(L_{\phi}) = \dil(L_{\phi}) = 3$ and $\length(L_{\phi}/xL_{\phi}) > 3$ for generic $x$.  (In fact, $I_2(\phi) = (\hat{x},\hat{y},\hat{z})^2$, so $\length(L_{\phi}/xL_{\phi}) = 4$ for generic $x$.)

%% file: thesis-appendix.tex
\chapter{Almost Cohen-Macaulay Modules with Linear Resolution}\label{App}
\setcounter{Th}{0}
In this appendix, we derive some consequences of Proposition \ref{linres}, reproduced here for convenience.

\begin{Pro}\label{App-linres} Let $M$ be a finitely generated graded module over $S=k[x_1,...,x_n]$ such that $M$ has a linear resolution, $\dim M=d$ and $\depth M \geq d-1$.  Then \begin{eqnarray*} \beta_{n-d+1} & = & \mu - e \\ \beta_{n-d} & = & (n-d+1)\mu - (n-d)e, \end{eqnarray*}
where $\beta_i = \beta_i(M)$ is the $i^{th}$ Betti number of $M$, $\mu = \nu(M)$ is the minimal number of generators of $M$, and $e = e(M)$ is the degree of $M$.
\end{Pro}

In Chapter \ref{hdil-chap} we established (Corollary \ref{dil-dim2}) the following theorem.

\begin{Th} \label{star-app} Let $M$ be a finitely generated graded $S$-module of dimension at most $2$ such that $M\cong \grm(M)$, where $S=k[x_1,\ldots,x_n]$ is a polynomial ring over the infinite field $k$.  Then $M$ has linear resolution if and only if $M$ is Dilworth and $\m \gm(M) = 0$. \end{Th}

The condition $M \cong \grm(M)$ is equivalent to the condition that $M$ be generated in degree $0$.  One direction of Theorem \ref{star-app}
\begin{quote} If $M$ is Dilworth and $\m \gm(M) = 0$ then $M$ has a linear resolution \end{quote}
is just a special case of Theorem \ref{main}, and is valid without the assumption that $\dim M \leq 2$.  To prove the converse
\begin{quote} If $M$ has a linear resolution then $M$ is Dilworth and $\m \gm(M) = 0$\end{quote}
we appealed to Theorem \ref{linresthm}.  Here, we shall give a second proof based on Proposition \ref{App-linres}.

\begin{Le}\label{lastminute} Suppose $\dim M \leq 1$ and $M/xM = M/ \m M$ for some $x$.  Then $M$ is Dilworth. \end{Le}

\proof Let $\mu = \nu(M)$ and $d = \dil(M)$.  We induct on $\mu$, the case $\mu=1$ being straightforward; thus suppose $\mu  > 1$.  Let $N = D(M)$.  We may assume $N \subseteq \m M$.  Indeed, if not, then let $a \in N \setminus \m M$ and let $\overline{M}=M/\langle a \rangle$.    Then $\nu(\overline{M}) = \mu - 1$ and $\dil(\overline{M}) \geq d-1$.  Since $\overline{M}/x \overline{M} = \overline{M}/ \m  \overline{M}$, by induction we have $\dil(\overline{M}) = \nu(\overline{M}).$  It follows that $d=\mu$.

Thus we assume $N \subseteq \m M$.  Let $N' = x[N:x]$.  We claim that $N'=N$.  It is clear that $N' \subseteq N$.  To show the opposite inclusion, let $a \in N$.  Then $a \in \m M = xM$, so $a=xb$ for some $b \in [N:x]$.  So $a \in N'$, which proves the claim.

Now $\nu([N:x]) < \nu(N)$, hence $\nu(N') \leq \nu([N:x]) < \nu(N)$, a contradiction.  \QED

\bigskip

\noindent {\bf Proof of Theorem \ref{star-app}.}  By Theorem \ref{main}, we have that if $M$ is Dilworth and $\m \gm(M) = 0$, then $M$ has linear resolution.  We need only show the converse.

Suppose $M$ has linear resolution.  Then $\m \gm(M) = 0$ and $M \cong \gm(M) \oplus M/\gm(M)$, so we may assume $\depth(M)>0$.  If $\dim(M) = 1$, then $M$ is Cohen-Macaulay, so the result follows from Proposition \ref{BHU}.  Thus, we may assume that $M$ has dimension $2$.  Let $x$ be generic on $M$ and let $\overline{M} = M/xM$.  Then, from the long exact sequence derived from
$$0 \lar M \stackrel{x}{\lar} M \lar M/xM \lar 0,$$ we have $\overline{M}_0 = M_1/xM_1$.  $\overline{M}$ has linear resolution and is of dimension $1$, so it is a Yoshida module.  In particular, $\m \overline{M}_0 = 0$.  Hence $M_1/xM_1 = M_1/\m M_1$, whence $M_1$ is Dilworth, by Lemma \ref{lastminute}.  Now $M_1$ is presented by
$$S^{\beta_{n-2}} \lar S^{\beta_{n-1}} \lar M_1 \lar 0,$$ where $\beta_i$ is the $i^{th}$ Betti number of $M$. Thus, by Proposition \ref{App-linres}, $\nu(M_1) = \mu-e$, where $\mu = \nu(M)$ and $e=\deg(M)$.  Therefore 
\begin{eqnarray*}\mu = \nu(M) \leq \hdil(M) & = & e+\dil(M_1)\\
& = & e+\nu(M_1) = e+\mu-e = \mu.\end{eqnarray*}
Thus $M$ is Dilworth. \QED

\begin{Ex} Let $S=k[x,y]$ be a polynomial ring in two variables.  Yoshino \cite{Yoshino} has classified $S$-modules with linear resolution.  Using Yoshino's classification, verify that all such modules are Dilworth, thus giving a {\em third} proof of (this special case of) Theorem \ref{star-app}. \end{Ex}

\begin{Pro} Let $M$ be a finitely generated graded $S$-module of dimension at most $2$ such that $M\cong \grm(M)$, where $S=k[x_1,\ldots,x_n]$ is a polynomial ring over the infinite field $k$.  Suppose that $M$ is Yoshida and unmixed.  Let $e=\deg(M)$.  then 
$$\hdeg(M) \in \{e,e+1,e+2,\ldots,2e\},$$
and each of these values may be obtained. \end{Pro}

\proof $\hdeg(M) = e + \length(M_1)$ and we must have $\m M_1 = 0$.  $M_1$ is presented by 
$$S^{\beta_{n-2}} \stackrel{\phi}{\lar} S^{\beta_{n-1}} \lar M_1 \lar 0,$$
where $\phi$ is a matrix of linear forms.  We have $\beta_{n-2} = (n-1)\mu - (n-2)e$ and $\beta_{n-1} = \mu-e$.  Since $\m M_1  = 0$, we have $\beta_{n-2} \geq n\beta_{n-1}$.  Hence 
$$(n-1)\mu-(n-2)e \geq n(\mu-e),$$ whence $e \leq \mu \leq 2e$.  

To see that each of these value may be obtained, let $e$ and $t$ be integers such that $e \geq 1$ and $0 \leq t \leq e$; let $\overline{S}=S/(x_3,\ldots,x_n)$; let $\mathfrak{n} = \m \overline{S}$; and let $$M = (\overline{S})^{e-t}\oplus\left(\bigoplus_1^t \mathfrak{n}\right).$$  Then $M$ is unmixed with $\deg(M)=e$ and $\hdeg(M) = \nu(M)=e+t.$ \QED

\bigskip

Similarly, we may derive:

\begin{Pro} If $M \cong \grm(M)$ and $\dim M = d$ and $\depth(M) \geq d-1$ and $M$ is linear Buchsbaum, then $$\length(M_{d-1}) = I(M) \leq e/(d-1).$$ \end{Pro}

%% file: thesis-vita.tex
\begin{vita}  
\heading{Tor Kenneth Gunston}
\vspace{15pt}


\begin{descriptionlist}{19xx xxxx} 

\item[1988]	Graduated from Foxborough High School, Foxborough, Massachusetts.

\item[1988-91]	Attended Stevens Institute of Technology, Hoboken, New Jersey.  Majored in Mathematics.

\item[1992]	B.S., Stevens Institute of Technology.

\item[1992-98]	Graduate work in Mathematics, Rutgers, The State University of New Jersey, New Brunswick, New Jersey.

\item[1993-97]	Teaching Assistant, Department of Mathematics.

\item[1996]	L.R. Doering and T. Gunston, Algebras Arising from Bipartite Planar Graphs, {\em Comm. Algebra} {\bf 24} (1996), 3589-3598.

\item[1996]	L.R. Doering, T. Gunston and W.V. Vasconcelos, Cohomological degrees and Hilbert functions of graded modules, {\em American J. Math.}, to appear.

\item[1998]	Ph.D. in Mathematics, Rutgers, The State University of New Jersey, New Brunswick, New Jersey.

\end{descriptionlist} 
\end{vita}